\documentclass[11pt,reqno]{amsart}
\usepackage{amsmath,amssymb,amsthm,mathtools,mathrsfs,upgreek}
\usepackage[margin=1.08in]{geometry}
\usepackage{enumitem}
\usepackage{graphicx}
\usepackage{tikz}
\usepackage[export]{adjustbox}
\usetikzlibrary{calc,arrows.meta,decorations.markings,decorations.pathreplacing}
\tikzset{arf diagram/.style={line width=0.5pt,>={Latex[length=4pt,width=3pt]}}}

\newcommand{\figOctagon}{%
\begin{tikzpicture}[arf diagram,x=0.75pt,y=0.75pt,yscale=-1,xscale=1]
\draw   (395.66,183.41) -- (348.41,230.66) -- (281.59,230.66) -- (234.34,183.41) -- (234.34,116.59) -- (281.59,69.34) -- (348.41,69.34) -- (395.66,116.59) -- cycle ;
\draw  [dash pattern={on 3pt off 3pt on 3pt off 3pt}]  (314.41,68.92) -- (314.91,229.42) ;
\draw  [fill={rgb, 255:red, 0; green, 0; blue, 0 }  ,fill opacity=1 ] (318.46,150) .. controls (318.46,148.09) and (316.91,146.54) .. (315,146.54) .. controls (313.09,146.54) and (311.54,148.09) .. (311.54,150) .. controls (311.54,151.91) and (313.09,153.46) .. (315,153.46) .. controls (316.91,153.46) and (318.46,151.91) .. (318.46,150) -- cycle ;
\draw (310,55) node [anchor=north west][inner sep=0.75pt]  {$\mathsf{t}$};
\draw (312,232) node [anchor=north west][inner sep=0.75pt]    {$\mathsf{t}'$};
\draw (318,110) node [anchor=north west][inner sep=0.75pt]    {$a_{O}$};
\draw (318,150) node [anchor=north west][inner sep=0.75pt]  {$o$};
\draw (370,76) node [anchor=north west][inner sep=0.75pt]   {$\mathsf{f}_{1+}$};
\draw (240,205) node [anchor=north west][inner sep=0.75pt]   {$\mathsf{f}_{2+}$};
\draw (376,205) node [anchor=north west][inner sep=0.75pt]  {$\mathsf{f}_{2-}$};
\draw (231,76) node [anchor=north west][inner sep=0.75pt]   {$\mathsf{f}_{1-}$};
\draw (399,142) node [anchor=north west][inner sep=0.75pt]   {$\mathsf{b}$};
\draw (218,140) node [anchor=north west][inner sep=0.75pt]   {$\mathsf{b}'$};
\end{tikzpicture}
}

\newcommand{\Gndiagram}
{\tikzset{every picture/.style={line width=0.75pt}} 
\begin{tikzpicture}[x=0.75pt,y=0.75pt,yscale=-1,xscale=1]

\draw   (10,175) .. controls (10,83.87) and (83.87,10) .. (175,10) .. controls (266.13,10) and (340,83.87) .. (340,175) .. controls (340,266.13) and (266.13,340) .. (175,340) .. controls (83.87,340) and (10,266.13) .. (10,175) -- cycle ;
\draw   (77.5,180) .. controls (77.5,163.43) and (90.93,150) .. (107.5,150) .. controls (124.07,150) and (137.5,163.43) .. (137.5,180) .. controls (137.5,196.57) and (124.07,210) .. (107.5,210) .. controls (90.93,210) and (77.5,196.57) .. (77.5,180) -- cycle ;
\draw   (160,180) .. controls (160,163.43) and (173.43,150) .. (190,150) .. controls (206.57,150) and (220,163.43) .. (220,180) .. controls (220,196.57) and (206.57,210) .. (190,210) .. controls (173.43,210) and (160,196.57) .. (160,180) -- cycle ;
\draw   (270,180) .. controls (270,163.43) and (283.43,150) .. (300,150) .. controls (316.57,150) and (330,163.43) .. (330,180) .. controls (330,196.57) and (316.57,210) .. (300,210) .. controls (283.43,210) and (270,196.57) .. (270,180) -- cycle ;
\draw  [color={rgb, 255:red, 0; green, 0; blue, 0 }  ,draw opacity=1 ][fill={rgb, 255:red, 0; green, 0; blue, 0 }  ,fill opacity=1 ] (7.5,180) .. controls (7.5,178.62) and (8.62,177.5) .. (10,177.5) .. controls (11.38,177.5) and (12.5,178.62) .. (12.5,180) .. controls (12.5,181.38) and (11.38,182.5) .. (10,182.5) .. controls (8.62,182.5) and (7.5,181.38) .. (7.5,180) -- cycle ;
\draw  [color={rgb, 255:red, 0; green, 0; blue, 0 }  ,draw opacity=1 ][fill={rgb, 255:red, 0; green, 0; blue, 0 }  ,fill opacity=1 ] (75,180) .. controls (75,178.62) and (76.12,177.5) .. (77.5,177.5) .. controls (78.88,177.5) and (80,178.62) .. (80,180) .. controls (80,181.38) and (78.88,182.5) .. (77.5,182.5) .. controls (76.12,182.5) and (75,181.38) .. (75,180) -- cycle ;
\draw  [color={rgb, 255:red, 0; green, 0; blue, 0 }  ,draw opacity=1 ][fill={rgb, 255:red, 0; green, 0; blue, 0 }  ,fill opacity=1 ] (157.5,180) .. controls (157.5,178.62) and (158.62,177.5) .. (160,177.5) .. controls (161.38,177.5) and (162.5,178.62) .. (162.5,180) .. controls (162.5,181.38) and (161.38,182.5) .. (160,182.5) .. controls (158.62,182.5) and (157.5,181.38) .. (157.5,180) -- cycle ;
\draw  [color={rgb, 255:red, 0; green, 0; blue, 0 }  ,draw opacity=1 ][fill={rgb, 255:red, 0; green, 0; blue, 0 }  ,fill opacity=1 ] (267.5,180) .. controls (267.5,178.62) and (268.62,177.5) .. (270,177.5) .. controls (271.38,177.5) and (272.5,178.62) .. (272.5,180) .. controls (272.5,181.38) and (271.38,182.5) .. (270,182.5) .. controls (268.62,182.5) and (267.5,181.38) .. (267.5,180) -- cycle ;
\draw    (150,180) -- (160,180) ;
\draw    (260,180) -- (270,180) ;
\draw    (40,180) .. controls (62.33,122.17) and (140.67,123.5) .. (150,180) ;
\draw    (40,180) .. controls (65,70.83) and (247.67,106.5) .. (260,180) ;
\draw [line width=2.25]    (10,180) -- (40,180) ;
\draw    (40,180) -- (80,180) ;
\draw (161,31.4) node [anchor=north west][inner sep=0.75pt]    {$G_{n}$};
\draw (12,180.9) node [anchor=north west][inner sep=0.75pt]    {$p_{1}$};
\draw (79.5,180.9) node [anchor=north west][inner sep=0.75pt]    {$p_{2}$};
\draw (162,181.4) node [anchor=north west][inner sep=0.75pt]    {$p_{3}$};
\draw (272,181.4) node [anchor=north west][inner sep=0.75pt]    {$p_{n}$};
\draw (29,101.4) node [anchor=north west][inner sep=0.75pt]    {$C_{1}$};
\draw (96.5,211.4) node [anchor=north west][inner sep=0.75pt]    {$C_{2}$};
\draw (179,211.4) node [anchor=north west][inner sep=0.75pt]    {$C_{3}$};
\draw (289,211.4) node [anchor=north west][inner sep=0.75pt]    {$C_{n}$};
\draw (231,173.4) node [anchor=north west][inner sep=0.75pt]    {$...$};
\draw (41,181.4) node [anchor=north west][inner sep=0.75pt]    {$\alpha _{2}$};
\draw (130,132.4) node [anchor=north west][inner sep=0.75pt]    {$\alpha _{3}$};
\draw (221,111.4) node [anchor=north west][inner sep=0.75pt]    {$\alpha _{n}$};
\draw (170,127.13) node [anchor=north west][inner sep=0.75pt]   {$...$};
\end{tikzpicture}}

\newcommand{\figQuotient}{%
\begin{tikzpicture}[arf diagram,x=0.75pt,y=0.75pt,yscale=-1,xscale=1,scale=0.9]
\draw     (140,70) -- (140,250) ;
\draw     (140,70) -- (450,70) ;
\draw     (140,250) -- (446,250) ;
\draw     (400,70) -- (400,130) ;
\draw     (446,190) .. controls (306.41,194.42) and (305.41,182.92) .. (306,160) .. controls (306.59,137.08) and (309.91,127.42) .. (446,130) ;
\draw     (446,70) .. controls (525.91,69.42) and (525.41,129.92) .. (446,130) ;
\draw     (400,190) -- (400,250) ;
\draw     (446,190) .. controls (525.91,189.42) and (525.41,249.92) .. (446,250) ;
\draw [color={rgb, 255:red, 208; green, 2; blue, 27 }  ,draw opacity=1 ]    (140,120) .. controls (250,120) and (306,130) .. (306,160) .. controls (306,190) and (250,200) .. (140,200) ;
\draw [<->] (516,90) .. controls (567.66,95.28) and (556.69,221.88) .. (516,230) ;
\draw     (140,330) -- (140,530) ;
\draw     (140,330) -- (450,330) ;
\draw     (140,530) -- (446,530) ;
\draw     (400,330) -- (400,530) ;
\draw     (450,330) .. controls (529.91,329.42) and (530.09,529.92) .. (446,530) ;
\draw [->] (90,240) -- (90,350) ;
\draw [color={rgb, 255:red, 208; green, 2; blue, 27 }  ,draw opacity=1 ]    (140,380) -- (270,380) ;
\draw [color={rgb, 255:red, 208; green, 2; blue, 27 }  ,draw opacity=1 ]    (140,480) -- (270,480) ;
\draw [color={rgb, 255:red, 208; green, 2; blue, 27 }  ,draw opacity=1 ]    (270,380) -- (270,480) ;
\draw (431,82) node [anchor=north west][inner sep=0.75pt]   [align=left] {$\displaystyle M$};
\draw (435,202) node [anchor=north west][inner sep=0.75pt]  [align=left] {$\displaystyle M$};
\draw (285,46.4) node [anchor=north west][inner sep=0.75pt]  {$W_{F,M}$};
\draw (261,231.4) node [anchor=north west][inner sep=0.75pt]   {$W_{F}$};
\draw (380,92.4) node [anchor=north west][inner sep=0.75pt]    {$Y'$};
\draw (380,212.4) node [anchor=north west][inner sep=0.75pt]    {$Y'$};
\draw (121,231.4) node [anchor=north west][inner sep=0.75pt]  {$Y$};
\draw (116,112.4) node [anchor=north west][inner sep=0.75pt]  [color={rgb, 255:red, 208; green, 2; blue, 27 }  ,opacity=1 ]  {$C_{1}$};
\draw (557,142.4) node [anchor=north west][inner sep=0.75pt]    {$\tau $};
\draw (217,105.4) node [anchor=north west][inner sep=0.75pt]  [opacity=1 ]  {\textcolor[RGB]{208,2,27}{$S$}};
\draw (116,193.4) node [anchor=north west][inner sep=0.75pt]  [color={rgb, 255:red, 208; green, 2; blue, 27 }  ,opacity=1 ]  {$C_{2}$};
\draw (280,305) node [anchor=north west][inner sep=0.75pt]    {$([ 0,1] \ \times \ Y') \cup M$};
\draw (374,418.4) node [anchor=north west][inner sep=0.75pt]   {$Y'$};
\draw (116,418.4) node [anchor=north west][inner sep=0.75pt]   {$Y'$};
\draw (116,370.4) node [anchor=north west][inner sep=0.75pt]  [color={rgb, 255:red, 208; green, 2; blue, 27 }  ,opacity=1 ]  {$C_{1}$};
\draw (116,472.4) node [anchor=north west][inner sep=0.75pt]  [color={rgb, 255:red, 208; green, 2; blue, 27 }  ,opacity=1 ]  {$C_{2}$};
\draw (435,410) node [anchor=north west][inner sep=0.75pt]  [align=left] {$\displaystyle M$};
\draw (94,281.4) node [anchor=north west][inner sep=0.75pt]    {$\pi _{W}$};
\draw (164,358) node [anchor=north west][inner sep=0.75pt]  [opacity=1 ]  {\textcolor[RGB]{208,2,27}{$[ 0,\frac{1}{2}] \times C_{1}$}};
\draw (164,482.4) node [anchor=north west][inner sep=0.75pt]  [opacity=1 ]  {\textcolor[RGB]{208,2,27}{$[ 0,\frac{1}{2}] \times C_{2}$}};
\draw (271,411.4) node [anchor=north west][inner sep=0.75pt]  [opacity=1 ]  {\textcolor[RGB]{208,2,27}{$\left\{\tfrac{1}{2}\right\} \times F$}};
\draw (122,532.4) node [anchor=north west][inner sep=0.75pt]    {$t=0$};
\draw (381,532.4) node [anchor=north west][inner sep=0.75pt]    {$t=1$};
\end{tikzpicture}
}

\newcommand{\figCharts}{%
\colorlet{bcol}{black}
\colorlet{tcol}{black}
\colorlet{wcol}{black}
\definecolor{ufill}{RGB}{210,228,248}
\definecolor{vfill}{RGB}{242,220,220}
\begin{tikzpicture}[arf diagram,x=0.75pt,y=0.75pt,yscale=-1,every node/.style={font=\small}]
\coordinate (v1) at (395.66,183.41);
\coordinate (v2) at (395.66,116.59);
\coordinate (v3) at (348.41,69.34);
\coordinate (v4) at (281.59,69.34);
\coordinate (v5) at (234.34,116.59);
\coordinate (v6) at (234.34,183.41);
\coordinate (v7) at (281.59,230.66);
\coordinate (v8) at (348.41,230.66);
\coordinate (tm) at ($(v3)!0.5!(v4)$);
\coordinate (tpm) at ($(v7)!0.5!(v8)$);
\fill[ufill] (tm)--(v3)--(v2)--(v1)--(v8)--(tpm)--cycle;
\fill[vfill] (tm)--(v4)--(v5)--(v6)--(v7)--(tpm)--cycle;
\fill[wcol,opacity=0.28] ($(tm)+(-6,0)$) rectangle ($(tpm)+(6,0)$);
\draw[bcol] (v1)--(v2);
\draw[wcol] (v2)--(v3);
\draw[tcol] (v3)--(v4);
\draw[wcol] (v4)--(v5);
\draw[bcol] (v5)--(v6);
\draw[wcol] (v6)--(v7);
\draw[tcol] (v7)--(v8);
\draw[wcol] (v8)--(v1);
\draw[dashed] (tm)--(tpm);
\fill (315,150) circle (2.6pt); \node[right=1pt] at (318,150) {\(o\)};
\node at (361,126) {\(U\)};
\node at (269,174) {\(V\)};
\node[bcol,right] at ($(v1)!0.5!(v2)+(-1,0)$) {\(\textsf{b}\)};
\node[bcol,left] at ($(v5)!0.5!(v6)+(0,0)$)  {\(\textsf{b}'\)};
\node[tcol,above] at (tm) {\(\textsf{t}\)};
\node[tcol,below] at (tpm) {\(\textsf{t}'\)};
\node at (345,190) {\(N(a_O)\)};
\end{tikzpicture}%
}

\newcommand{\figUpstairsConnected}{%
\begin{tikzpicture}[arf diagram,x=0.75pt,y=0.75pt,yscale=-1,xscale=1]
\draw     (140,70) -- (140,250) ;
\draw     (140,70) -- (450,70) ;
\draw     (140,250) -- (446,250) ;
\draw     (396,70) -- (396,130) ;
\draw     (446,190) .. controls (306.41,194.42) and (305.41,182.92) .. (306,160) .. controls (306.59,137.08) and (309.91,127.42) .. (446,130) ;
\draw     (446,70) .. controls (525.91,69.42) and (525.41,129.92) .. (446,130) ;
\draw     (396,190) -- (396,250) ;
\draw     (446,190) .. controls (525.91,189.42) and (525.41,249.92) .. (446,250) ;
\draw [color={rgb,255:red,208;green,2;blue,27},draw opacity=1]
  (140,120)
  .. controls (361.3333,120) and (361.3333,200) ..
  (140,200);
\draw [<->] (516,90) .. controls (567.66,95.28) and (556.69,221.88) .. (516,230) ;
\draw (428,90) node [anchor=north west][inner sep=0.75pt] [align=left] {$\displaystyle M$};
\draw (428,212) node [anchor=north west][inner sep=0.75pt]  [align=left] {$\displaystyle M$};
\draw (280,40) node [anchor=north west][inner sep=0.75pt]   {$W_{F,M}$};
\draw (211,211.4) node [anchor=north west][inner sep=0.75pt]   {$W_{F}$};
\draw (377,92.4) node [anchor=north west][inner sep=0.75pt]    {$Y'$};
\draw (378,212.4) node [anchor=north west][inner sep=0.75pt]    {$Y'$};
\draw (121,212.4) node [anchor=north west][inner sep=0.75pt]   {$Y$};
\draw (123,113.4) node [anchor=north west][inner sep=0.75pt]  [color={rgb, 255:red, 208; green, 2; blue, 27 }  ,opacity=1 ]  {$C$};
\draw (557,142.4) node [anchor=north west][inner sep=0.75pt]    {$\tau $};
\draw (167,102.4) node [anchor=north west][inner sep=0.75pt]  [color={rgb, 255:red, 208; green, 2; blue, 27 }  ,opacity=1 ]  {$\text{Fix}( \tau ) \cong F$};
\end{tikzpicture}
}

\usepackage[colorlinks=true,linkcolor=blue!55!black,citecolor=blue!55!black,urlcolor=blue!55!black]{hyperref}

\theoremstyle{plain}
\newtheorem{theorem}{Theorem}[section]
\newtheorem{proposition}[theorem]{Proposition}
\newtheorem{lemma}[theorem]{Lemma}
\newtheorem{corollary}[theorem]{Corollary}
\newtheorem{conjecture}[theorem]{Conjecture}

\theoremstyle{definition}
\newtheorem{definition}[theorem]{Definition}

\newtheorem{remark}[theorem]{Remark}

\newcommand{\F}{\mathbb F_2}
\newcommand{\Z}{\mathbb Z}
\newcommand{\R}{\mathbb R}
\newcommand{\Q}{\mathbb Q}
\newcommand{\C}{\mathbb C}

\newcommand{\Spinc}{\operatorname{spin}^{c}}
\newcommand{\spinc}{\ensuremath{\operatorname{spin}^{c}}}
\newcommand{\Arf}{\operatorname{Arf}}
\newcommand{\PD}{\operatorname{PD}}
\newcommand{\rad}{\operatorname{rad}}
\newcommand{\RRSC}{\textsc{rrsc}}
\newcommand{\RSpin}{\textsc{rspin}}
\newcommand{\tsr}{\widetilde{\text{\normalfont\textsc{r}}}}
\newcommand{\sr}{\text{\normalfont\textsc{r}}}

\newcommand{\Fix}{\operatorname{Fix}}
\newcommand{\gr}{\operatorname{gr}}
\newcommand{\ind}{\operatorname{ind}}

\newcommand{\Fh}{\widehat F}
\newcommand{\Tcal}{\mathcal{T}}

\newcommand{\ord}{\mathsf{ord}}
\newcommand{\dbcv}{\mathsf{\Sigma}}
\newcommand{\sF}{\mathfrak s}
\newcommand{\tf}{\mathfrak t}

\newcommand{\hHFR}{\text{\smaller{$\widehat{\textit{HFR}}$}}}

\newcommand{\Lift}{\mathscr L}

\newcommand{\Bh}{\widehat B}

\newcommand{\Bin}{B_{\mathrm{in}}}
\def\s{\mathfrak{s}}

\DeclareMathOperator{\tHMR}{\text{\smaller{$\widetilde{\textit{HMR}}$}}}
\DeclareMathOperator{\thmr}{\text{\smaller{$\widetilde{\textit{HMR}}$}}}

\DeclareMathOperator{\tcmr}{\text{\smaller{$\widetilde{\textit{CMR}}$}}}

\DeclareMathOperator{\Rot}{\mathsf{Rot}}

\title[Absolute mod-two gradings on framed real monopole Floer homology]
{Absolute mod-two gradings\\on framed real monopole Floer homology}

\author[Jiakai Li]{Jiakai Li}
\address{Department of Mathematics, Stanford University\\
450 Jane Stanford Way, Building 380, Stanford, CA 94305, USA}
\email{lijk@stanford.edu}

\author[Ciprian Manolescu]{Ciprian Manolescu}
\address{Department of Mathematics, Stanford University\\
450 Jane Stanford Way, Building 380, Stanford, CA 94305, USA}
\email{cm5@stanford.edu}
\begin{document}

\begin{abstract}
We construct absolute mod-two gradings on framed real monopole Floer homology, for real 3-manifolds such that the image of the fixed point set in the quotient is nullhomologous. This allows us to fix the sign of Miyazawa's degree invariant. Our construction uses a mod-two index-theoretic correction term, the Arf invariant.
 \end{abstract}

\maketitle

\section{Introduction}

In a series of papers~\cite{Li22,Li23,Li24,Li26}, the first author developed {\em real monopole Floer homology}, an invariant of real three-manifolds (i.e., three-manifolds $Y$ equipped with an involution $\tau$ with fixed point set of codimension two). An alternative construction was given by Konno--Miyazawa--Taniguchi~\cite{KMT:homology}; later, Xiao~\cite{Xiao} proved that these two invariants are isomorphic. Real monopole Floer homology is an essential piece of real Seiberg-Witten theory, whose applications include new results about surfaces in four-manifolds \cite{miyazawa,Baraglia:surface} and knot concordance \cite{KPT:cables}.

A typical way to extract additional information from Floer homologies is through their gradings. These can be relative (i.e., given two generators, only the difference of their gradings is well-defined as an element in an abelian group $G$), or absolute (i.e., to each generator we can assign an element of $G$ in a natural way). Real monopole Floer homology has a relative $\Z/n$ grading for some $n$ (which could be odd, and depends on the three-manifold and its $\text{spin\textsuperscript{c}}$ structure); see \cite[Section 11.2]{Li22}. On the other hand, ordinary monopole Floer homology, as developed by Kronheimer and Mrowka, has a relative $\Z/2n$ grading, and hence a relative $\Z/2$ grading. Furthermore, the $\Z/2$-grading can be made absolute. There is also an absolute  $\Q$-grading in case the $\text{spin\textsuperscript{c}}$ structure is torsion. See \cite[Sections 22.3, 22.4, 28]{KM}.

There is a version of real monopole Floer homology, the framed theory $\tHMR$, which was shown in \cite{Li26} to admit  a relative $\Z/2$-grading. This version is constructed from a branched covering of the usual real monopole configuration space, by using a set $\mathbf{p}$ of markings (basepoints) on the three-manifold. While the real monopole Floer homologies defined in \cite{Li22} depend on a real \text{spin\textsuperscript{c}} structure, the framed version depends on more data: Specifically, it is associated to a relative real \text{spin\textsuperscript{c}} ($\RRSC$) structure $\widetilde{\sr}$ on $(Y,\tau,\mathbf{p})$. An $\RRSC$ is a triple consisting of a \text{spin\textsuperscript{c}} structure, a trivialization $T$ of the spinor bundle at $\mathbf{p}$, and a Clifford compatible conjugate-linear involution $\sr$ lifting $\tau$, for which $T\circ\sr\circ T^{-1}$ is the standard complex conjugation.

The purpose of this paper is to define absolute $\Z/2$-gradings on $\tHMR$, lifting the relative one from \cite{Li26}. (As we shall see, the lift may depend on additional topological data: the homology class of a Seifert surface, and the ordering of the components of the fixed point set.) The absolute $\Z/2$-gradings have several benefits. For example, they allow us to compare the gradings between generators in different \text{spin\textsuperscript{c}} structures. More importantly, they will be helpful for getting more information from exact triangles \cite{LiNew}, and for more general cut-and-paste techniques for the four-dimensional real Seiberg--Witten invariants.

\subsection{Setup and notation}
Let $Y$ be a closed oriented 3-manifold, and let $\tau \colon Y \to Y$ be an orientation-preserving involution.
Let $Y'$ denote the quotient manifold $Y/\tau$,
and $\pi$ denote the projection map $Y \to Y'$.
Write the fixed-point set $C = \Fix(\tau)=C_1\sqcup\cdots\sqcup C_n$, which is assumed to be a non-empty link. Fix also an orientation on $C$.

Suppose $C$ is nullhomologous in $Y'$, and let $F \subset Y'$ be a connected oriented Seifert
surface for $C$, with $\partial F=C$. (Note that, unlike in the real Heegaard Floer setting considered in \cite{Srivastava}, the complement of $F$ does not have to be a handlebody.) Assume the branched cover on
$Y\setminus C$ is the cover determined by mod-two intersection with $F$.
Let $\Sigma$ be the double of $F$ along its boundary,
which can be identified with the
preimage $\Sigma=\pi^{-1}(F)$.
Then the manifold $Y'$-cut-along-$F$ has boundary $\partial (Y' \setminus N'(F)) \cong \Sigma$.
Concretely, cutting $Y$ along $\Sigma$ produces the two copies of $Y'$ cut along $F$.
In particular, $\Sigma$ separates $Y$.
Choose a side of $\Sigma$, equivalently one of the two lifts of a small
Seifert push-off of $F$, and choose one basepoint $p_k\in C_k$.

\subsection{The knot case}
\label{sec:knot_idea}
Suppose that $C$ is connected,
and $(\mathfrak{s},\sr)$ is a real $\Spinc$ structure on $(Y,\tau)$.
(Here, $\s$ is the $\Spinc$ structure and $\sr$ is the real structure on the spinor bundle.)
Choose a real $\Spinc$ cap
$(W,\tf)$ with $\partial W=-Y$ that extends $(\mathfrak s,\sr)$.
Fix a framed real configuration $\mathfrak a$ on $Y$, and use the corresponding Atiyah--Patodi--Singer (APS) boundary condition on $W$.
Let $\gr_{SW}(\mathfrak{a},W,\tf)$ be the formal dimension of the Seiberg--Witten moduli space with this APS condition.
Write $D^+_{\tf}$ for the Dirac operator on the cap and $D^{+,\sr}_{\tf}$ for its real part.
Given such a $W$,
define
\begin{equation*}
 \iota_R(W):= b^+_R(W) - b^1_R(W).
\end{equation*}
where $b^i_R(W) = \dim_{\mathbb{R}} H^i(W;\mathbb{R})^{-\tau^*}$ is the anti-invariant part of the $i$-th Betti number, cf.~\cite[Definition 25.4.1]{KM}.
Write
\[
\mathcal I_{\tf}(\mathfrak a) = \gr_{SW}(\mathfrak{a},W,\mathfrak{t}) + \iota_R(W).
\]

The quantity $\mathcal I_{\tf}(\mathfrak a)$ should be viewed as a preliminary absolute mod-two grading which is, however, not independent of $\mathfrak{t}$.
Indeed, when $\mathfrak{\mathfrak{a}}$ is a reducible configuration with invertible Dirac operator,
\[\mathcal I_{\tf}(\mathfrak a) = \ind^{APS} D^{+,\sr}_{\tf}.\]
In the ordinary case, a complex Dirac operator contributes to the Seiberg--Witten dimension formula twice its complex index.
In the real case, the real Dirac operator contributes the complex index by conjugate-linearity.
Let $L$ be a Real line bundle on $W$ with a fixed trivialization on $Y$,
and let $x=c_1(L)\in H^2(W,Y;\Z)$.
The APS formula implies
\[
\mathcal I_{\tf \otimes L}(\mathfrak a) - \mathcal I_{\tf}(\mathfrak a)=
\ind_{\R}^{APS}D^{+,\sr}_{\tf \otimes L}
   -\ind_{\R}^{APS}D^{+,\sr}_{\tf}
 =\ind_{\C}^{APS}D^+_{\tf \otimes L}
   -\ind_{\C}^{APS}D^+_{\tf}
 =\frac{1}{2}\bigl(c_1(\mathfrak{t}) \cdot x + x \cdot x\bigr).
\]
Here the products are evaluated
on the relative fundamental class $[W,Y]$.
Since $c_1(\tf)$ is a characteristic element, the expression in parentheses is
even; nevertheless, the extra factor of $1/2$ means that the above change need not be
even. As a preliminary absolute grading, the quantity $\mathcal{I}_{\tf}$ therefore needs a further correction.\footnote{
	So far the construction is reminiscent of the absolute $\mathbb{Q}$-grading in Floer homology.
	However, we allow non-torsion \text{spin\textsuperscript{c}} structures on $Y$,
	and we seek a correction term to $\mathfrak{t}$ that is mod-two index-theoretic in nature,
	as opposed to $\eta$-type invariants which can be non-integral rational numbers.}

\begin{figure}
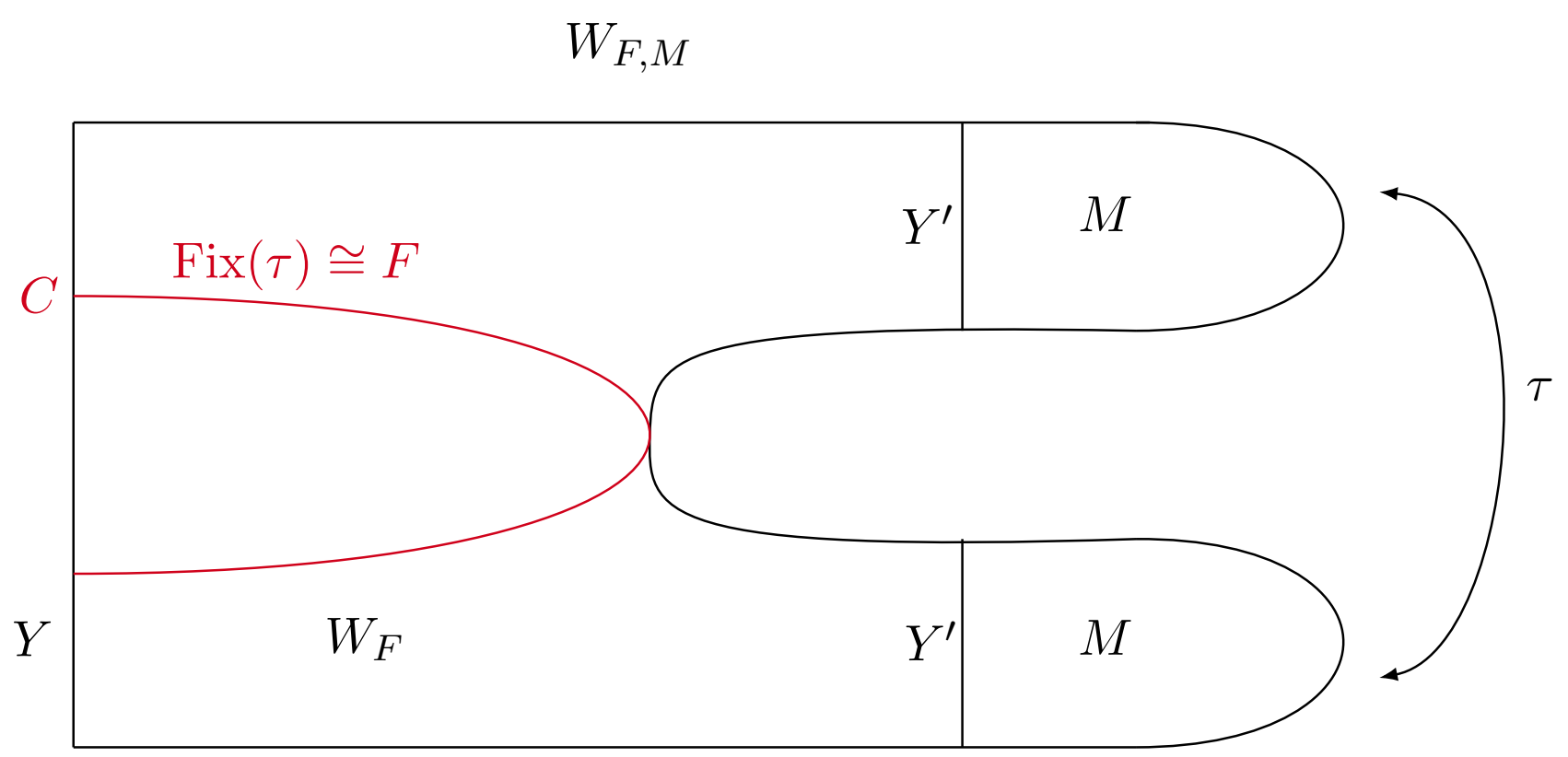

	\adjustbox{max width=0.7\textwidth}{\figUpstairsConnected}
	\caption{Cap of $Y$ constructed from the Seifert cobordism $W_{F} \colon Y \to Y' \sqcup Y'$ and two copies of an outgoing cap $M$. The fixed point set can be identified with a pushed-in Seifert surface $F$.}
	\label{fig:upstairs-only}
\end{figure}

The Seifert surface $F$ organizes the choices of extensions $\mathfrak{t}$,
by relating the real manifold $(Y, \tau)$ to a disconnected one, as follows. We push  $F$ into the interior of $ [0,1] \times Y'$, and take the double branched cover. The result is a cobordism  $W_F$ from $Y$ to $Y' \sqcup Y'$, which we call the {\em Seifert cobordism} from $Y$ to $Y' \sqcup Y'$. We then obtain a cap  $W_{F, M}$ by filling its outgoing ends with two copies of a four-manifold $M$ with $\partial M = -Y'$. See Figure~\ref{fig:upstairs-only}. We ask that $M$ is {\em admissible}, i.e., $H_1(M; \Z)=0$, so that every $\spinc$ structure on $Y'$ extends to $M$. The extensions $\tf$ of $\s$ over $W_F$ rel $Y$ form a torsor for \[\Lambda= H^1(F, \partial F; \Z).\]
Denote the twist of $\mathfrak{t}$ by $\theta\in\Lambda$ as $\tf^\theta$.

Via a slight abuse of notation, we also denote by $\tf$ a compatible extension over the
cap $W_{F, M}$.
Record the index change on $W_{F, M}$ induced by $\tf \to \tf^\theta$ by
\[
 Q_{\tf}(\theta)
 :=\ind_{\C}^{APS}D^+_{\tf^{\theta}}
   -\ind_{\C}^{APS}D^+_{\tf}\pmod2.
\]
This parity is independent of choice of $M$,
because real Seiberg--Witten theory over a pair of disjoint manifolds $M \sqcup M$ with conjugate data reduces to the ordinary theory on $M$. It turns out that  $Q_{\tf}(\theta)$ depends only on $\theta\bmod2$. Thus, we get a quadratic refinement
$Q_{\tf}$ of the intersection form $B$ on $H_1(F; \mathbb{Z}/2)$, satisfying
$Q_{\tf^\theta}=Q_{\tf}+B(\theta,-)$.
For a knot, $B$ is nondegenerate.
The first important idea is to introduce the Arf invariant,
which satisfies the following translation formula:
\[
 \Arf(Q_{\tf^\theta})-\Arf(Q_{\tf})=Q_{\tf}(\theta).
\]
Thus $\Arf(Q_{\tf})$ is a \emph{primitive} of the index change induced by changes in $\tf$.
From here we can introduce the knot-version of the absolute mod-two grading:
\begin{equation}
\label{eq:grading-knot-intro}
\gr^{(2)}(\mathfrak a)
=\gr_{\mathrm{SW}}(\mathfrak a, W_{F,M}, \tf)
+\iota_R(W_{F, M})
+\Arf(Q_{\mathfrak{t}}) \ \bmod2.
\end{equation}

\subsection{The link case}

The preceding approach stops short of the link case,
when $B$ becomes degenerate.
Indeed, the Arf invariant is not well-behaved in this case,
and $Q_{\mathfrak{t}}$ fails to detect the difference between relative and absolute real \text{spin\textsuperscript{c}} structures.
To return to the non-degenerate case, we must use a relative real \text{spin\textsuperscript{c}} structure $\widetilde{\sr}$ on $Y$. We fix the choice of such a structure.

The next important idea is to turn $Q_{\mathfrak{t}}$ into spin-geometric data on the surface $F$.
If $F$ were closed, Johnson's correspondence~\cite{Johnson} would associate
to every quadratic form $Q$ on $H_1(F;\F)$
a spin structure $\mathfrak{s}$ on $F$, satisfying
\begin{equation*}
	\Arf(Q) = \Arf(\mathfrak{s}).
\end{equation*}
Here, $\Arf(\mathfrak{s}) \in \Z/2$ tells us whether the majority of the curves on $F$ bound a spin disk or not. Moreover, we  would have $\Arf(\mathfrak{s}) \equiv \dim \ker(\overline{\partial}_{\mathfrak{s}}) \pmod 2$~\cite{Atiyah71},
where $\overline{\partial}_{\mathfrak{s}}$ denotes the Cauchy--Riemann operator on the Riemann surface $F$
and $\mathfrak{s}$ is interpreted as a square root of the canonical bundle.

In our setting, we construct a spin structure $\sF_{\tf}$ on $F$ directly from the real spinors of $\mathfrak{t}$ on $W_F$,
and prove that this construction recovers $Q_{\mathfrak{t}}$ under Johnson's correspondence. Consider the closure $\widehat{F} = G \cup_C F$ by a planar surface $G$. (Since $\widehat{F}$ is closed, its intersection form is nondegenerate.)
The given orientation of $C$, together with an ordering of the components of $C$ and the relative real \text{spin\textsuperscript{c}} structure $\widetilde{\sr}$,
yields a spin structure $\widehat{\mathfrak{s}}_{\tf}$ on $\widehat{F}$,
for which $\widehat{\mathfrak{s}_{\mathfrak{t}}}|_F = \mathfrak{s}_{\mathfrak{t}}$.
This implies that $\Arf(q_{\widehat{\mathfrak{s}}_{\mathfrak{t}}})$ offsets the index change under $\mathfrak{t} \to \mathfrak{t}^{\theta}$.
This construction produces the desired correction term for any $n = |\pi_0(C)| \geq 1$,
and recovers the correction term~\eqref{eq:grading-knot-intro} in the knot case $(n=1)$.

\begin{theorem}\label{thm:main}
Under the hypotheses above,
an oriented Seifert surface $F$ determines an absolute mod-two grading $\gr^{(2)}$ of the framed real monopole Floer homology $\tHMR(Y, \tau; \tsr)$.
At a configuration $\mathfrak a$:
\begin{equation}
 \gr^{(2)}(\mathfrak a)
 =\gr_{\mathrm{SW}}(\mathfrak a;W_{F,M},\tf)+\iota_R(W_{F,M})
  +\Arf(\, \widehat{\mathfrak{s}}_{\mathfrak{t}})\pmod2.
\end{equation}
Here, $M$ is an admissible outgoing cap, $W_{F,M}$ is the capped cobordism,
and $\gr_{\mathrm{SW}}$ is the real Seiberg--Witten index with APS boundary condition given by $\mathfrak a$.
The grading $\gr^{(2)}$ is independent of the extension $\mathfrak{t}$ and the \text{spin\textsuperscript{c}} cap $M$. However, it may depend on an ordering of the components of $C$, and it may depend on $F$, but only through its class $[F]\in H_2(Y',C;\Z)$.
\end{theorem}

\subsection{Comparison to real Heegaard Floer homology}
In \cite{GM}, Guth and the second author defined real Heegaard Floer homology, a symplectic-geometric counterpart to real monopole Floer homology which is expected to be isomorphic to it. Srivastava \cite{Srivastava} defined an absolute $\Z/2$-grading on the hat version of real Heegaard Floer homology, $\hHFR(Y, \tau, \mathbf{p})$. This depends on $[F]$ and the ordering of the components of $C$, just like in our case.

In fact, we expect our grading $\gr^{(2)}$ to correspond to the one in \cite{Srivastava}, with an additional correction term. (The difference can be seen even in simple examples such as the trefoil.) We define the normalized mod-$2$ grading by the formula
\begin{equation} \label{eq:normalization-gr2} \widehat{\gr}:= \gr^{(2)} + \frac{\sigma(W_F)-b_1(Y)-2b_1(Y')-|C|+1}{2}\pmod2.
\end{equation}
When $Y$ is the double branched cover of a link $L \subset S^3$, we can equivalently write
\begin{equation} \label{eq:seifert-normalization-gr2} \widehat{\gr}:= \gr^{(2)} + \frac{-\sigma(L)-\mathrm{null}(L)-|L|+1}{2}.
\end{equation}
Here, $\sigma(L)$ and $\mathrm{null}(L)$ are the signature and nullity of $L$, respectively.

We phrase the expected graded equivalence between the two theories as follows.

\begin{conjecture}
\label{conj:HFR}
Fix the homology class $[F]$ and the ordering of the components of $C$. Then, there is a grading-preserving isomorphism between the framed real monopole Floer homology $\tHMR(Y, \tau, \mathbf{p})$ equipped with the $\widehat{\gr}$ grading  from \eqref{eq:normalization-gr2}, and the hat real Heegaard Floer homology $\hHFR(Y, \tau, \mathbf{p})$ equipped with the grading from \cite{Srivastava}.
\end{conjecture}

\begin{remark}
	While the mod-two grading as a $\mathbb{Z}/2$-valued function on $\tHMR(Y, \tau, \sr)$ may depend on the ordering of $C$,
	the isomorphism class of $\tHMR(Y, \tau, \sr)$, as a graded $\F$-vector space,  is independent of this ordering (Theorem~\ref{thm:long-absolute-grading}(iii)).
	This is the same case as in Srivastava's grading on $\hHFR(Y, \tau, \mathbf{p})$.
	Intriguingly, the two absolute gradings encounter this dependence through very different routes:
	On the monopole side, the ordering was needed to pin down a consistent assignment of the spin marking data on the planar capping surface,
	whereas on the Heegaard Floer side, the ordering was used for the cancellation of generator signs.
\end{remark}

\subsection{The degree invariant}
In \cite{miyazawa}, Miyazawa defined an invariant of surfaces in four-manifolds, called the degree, which he used to detect an exotic embedding of $\mathbb{RP}^2$ in $S^4$. The exotic embedding is the connected sum of the standard one and the roll spun of the $P(-2,3,7)$ pretzel knot. In general, Miyazawa's invariant for the roll spun of a knot $K \subset S^3$ is a $\mathbb{Z}/{\pm 1}$-valued numerical invariant, denoted $|\deg(K)|$. It equals the absolute value of the Euler characteristic of $\tHMR(\dbcv_2(K), \s_0)$, where $\dbcv_2(K)$ denotes the double branched cover of the knot, and $\s_0$ is the spin structure. One can also define degrees from the other $\Spinc$ structures on $\dbcv_2(K)$.

It has been unknown whether $|\deg(K)|$ is a new knot invariant, or something derived from classical knot invariants.
Using our absolute grading $\widehat{\gr}$ from \eqref{eq:seifert-normalization-gr2}, Miyazawa's degree invariant can be defined over $\mathbb{Z}$. Furthermore, we can define a total degree, as the sum of the degrees over all (relative) real $\Spinc$ structures. The analogue of the total degree in real Heegaard Floer homology was shown to be $\Delta_K(i)$ by Srivastava~\cite{Srivastava}. (Here, $\Delta_K$ is the Alexander polynomial.) Our calculations in this paper support that the total degree in $\tHMR$ equals $\Delta_K(i)$ as well. In fact, it will be shown in \cite{LiNew} that this is true for any knot.

\subsection{Organization of the paper}
In Section~\ref{sec:cobordism} we analyze the topology of the Seifert cobordism $W_F$ and the cap $W_{F, M}$, the real $\Spinc$ structures on them, as well as the (normalized) index of the Dirac operator with an APS boundary condition. In Section~\ref{sec:spin} we relate the index to the Arf invariants of spin structures on the Seifert surface.
Section~\ref{sec:arf} proves the main theorem, and in Section~\ref{sec:examples} we compute some examples. Appendix~\ref{app:cover-indentification} contains the proof of a technical lemma.

\subsection{Notational conventions}
Unless otherwise noted, all manifolds are compact, connected, and oriented.

\subsection{Statement of AI use}
We used Claude Fable 5.0, GPT-5.6 Sol and GPT-6 Astra in the following ways.  AI assisted with the classical knot invariant computations in \S\ref{sec:examples}, such as the Alexander polynomial and the signature. It also assisted with Clifford algebra calculations, which led to the proof of Lemma~\ref{lem:Gamma}, given in Appendix~\ref{app:cover-indentification}. AI produced the lens space grading computations in~Lemma~\ref{lem:gr-Lpq}. The AI results were rewritten by the authors. We used AI for creating the tables, for producing the shading in Figure~\ref{fig:charts}, for checking the paper for mathematical errors, and for proofreading.

\subsection{Acknowledgments} The first author would like to thank Peter Kronheimer for many helpful conversations and arranging his visit to Harvard in March 2026. He would also like to thank Tom Mrowka for teaching the Arf invariant in a  Riemann Surfaces course at MIT.

The second author was partially supported by a Simons Investigator Award, the Simons Collaboration Grant on New Structures in Low-Dimensional Topology, and the NSF AIMing grant DMS-2522743.

\section{The Seifert cobordism and the normalized APS index}\label{sec:cobordism}

\subsection{The local octagon model}\label{sec:octagon}
\begin{figure}[htb]
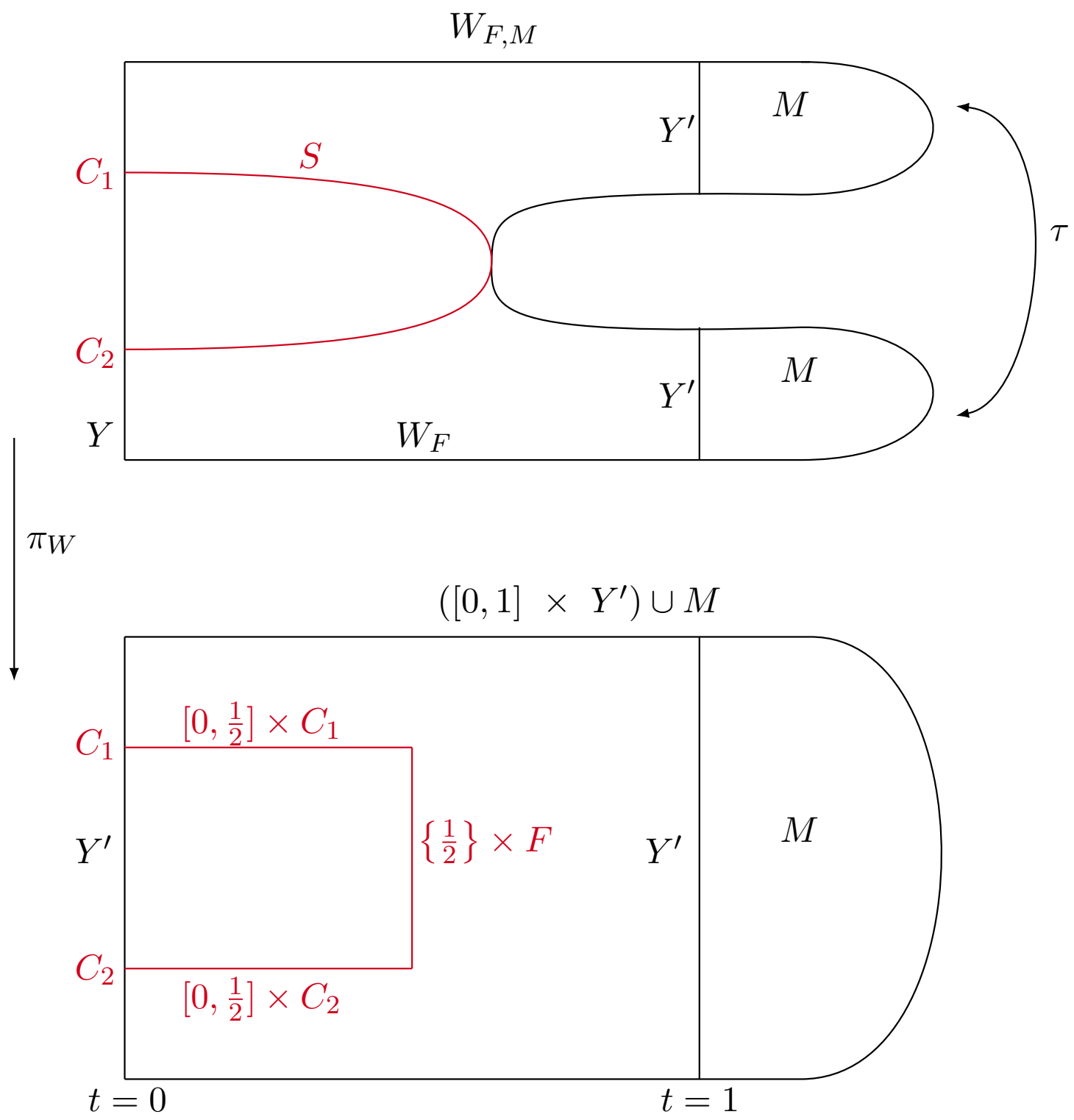

\centering
\adjustbox{max width=0.6\textwidth}{\figQuotient}
\caption{The Seifert cobordism and its quotient by branched covering map.
Below is the cylinder
$[0,1]\times Y'$ capped by $M$, with branch set
$\bigl([0,\tfrac12]\times C\bigr)\cup\bigl(\{\tfrac12\}\times F\bigr)$, rounded to the surface $S$ of Section~\ref{sec:octagon}.
Above is the branched double cover $W_{F,M}$,
obtained by capping the outgoing ends of $W_F\colon Y\to Y'\sqcup Y'$ with two copies of $M$. In our picture, $C$ consists of two components $C_1$ and $C_2$.}
\label{fig:quotient}
\end{figure}

Following the set-up from the introduction, let $F$ be  a Seifert surface for $C \subset Y'$.
Let $S\subset[0,1]\times Y'$ be the push-off of $F$ into the cylinder, i.e., a properly embedded surface obtained from
\[
 \bigl(\{\tfrac12\}\times F\bigr)\cup\bigl([0,\tfrac12]\times\partial F\bigr)
\]
by rounding the corner along $\partial F$.
In particular, $S$ is homeomorphic to $F$,
and $\partial S = \{0\}\times C$ which is disjoint from $\{1\}\times Y'$.
The \emph{Seifert cobordism} is the double cover of $[0,1] \times Y'$ branched along $S$:
\begin{equation}
	\pi_W \colon W_F \longrightarrow [0,1] \times Y',
\end{equation}
specified as follows.
Round the 3-dimensional cut $[0,\tfrac12]\times F$ compatibly with $S$;
the boundary of the cut consists of $S$ and the incoming copy of $F$.
A loop in the complement of $S$ lifts to an open arc in $W_F$ when its mod-two intersection with the cut is one.
Write the deck involution again as $\tau$.
The cut restricts to $F$ at the incoming end and is disjoint from the outgoing end.
Thus restricting $\pi_W$ to the boundary of $W_F$,
we obtain the prescribed $\pi \colon Y \to Y'$ and the trivial double cover of $Y'$.
With boundary orientations:
\begin{equation}
 \partial W_F=-Y\ \sqcup\ \bigl(Y'\sqcup Y'\bigr).
\end{equation}
The deck transformation $\tau$ preserves the incoming $Y$ while interchanging the two outgoing copies of $Y'$. See Figure~\ref{fig:quotient}.

Let us interpret $W_{F}$ as an excision cobordism.
Downstairs, choose a regular neighbourhood $N'(F) \subset Y'$ of $F$,
containing a tubular neighbourhood around $C$.
Let $X' = \overline{Y' \setminus N'(F)}$ be the \emph{$Y'$-cut-along-$F$}.
Then $\partial X' \cong \Sigma$, the double of $F$;
gluing two copies of $X'$ yields $Y$.
The cobordism $W_F$ decomposes as $A \cup B$, where
\begin{equation*}
	A := \pi_W^{-1}([0,1] \times X'), \qquad
	B := \pi_W^{-1}([0,1] \times N'(F)).
\end{equation*}
The subset $B$ contains all of the branch set and cut,
and hence $\pi_W|_A$ is the trivial double cover.
We write
\begin{equation*}
	A = [0,1] \times (X_1 \sqcup X_2), \quad
	A \cap B \cong [0,1] \times (\Sigma_1 \sqcup \Sigma_2),
\end{equation*}
where $X_i$ is a copy of $X'$, and $\Sigma_i := \partial X_i$.

Let $N'(C) \subset N'(F)$ be a tubular neighbourhood of $C$.
Let $$F_0 = F \setminus N'(C),$$ so that $N'(F)$ can be viewed as the union of $N'(C)$
and a thickening of $F_0$.
An interior region of the rounded surface $S \subset ([0,1] \times Y')$ can thus be identified with $F_0$.
The rounding of $S$ takes place in the following tubular region:
\begin{equation*}
	T := \pi_W^{-1}([0,1] \times N'(C)).
\end{equation*}

\begin{figure}[htb]
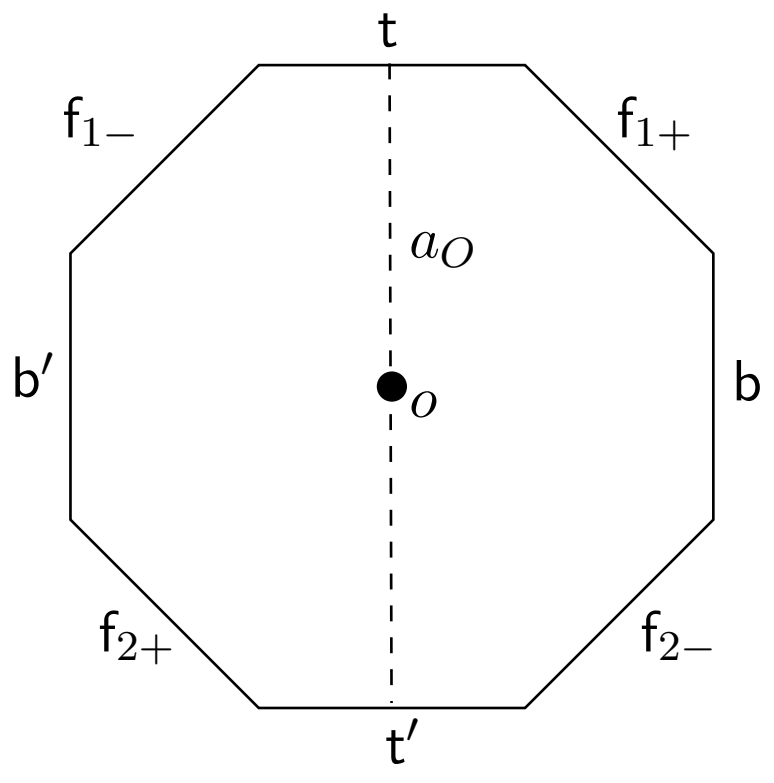

\centering
\figOctagon
\caption{The octagon $O$.
The incoming (bottom) arcs are $\mathsf{b},\mathsf{b}'$,
and the outgoing (top) arcs are $\mathsf{t},\mathsf{t}'$.
The dashed deck-invariant arc $a_O$ joins $\mathsf{t}$ to $\mathsf{t}'$.}
\label{fig:octagon}
\end{figure}

To describe the set $B$,
consider the following model of an octagon $O$:
Let $P = [0,1]_t \times [-1,+1]_s$, and
$O \to P$ be the double cover of $P$ branched at the point $(1/2,0)$.
The surface $O \cong D^2$ can be naturally depicted as an octagon.
The eight edges of $\partial O$ correspond to the preimage of the four edges of the rectangle $P$.
Denote the preimages of the side $t=0$ by $\mathsf{b},\mathsf{b}'$ and the preimages of $t=1$ by $\mathsf{t},\mathsf{t}'$.
The covering involution of $O \to P$ acts by a $180^{\circ}$-rotation,
fixing precisely the point $o$ at the centre of $O$. See Figure~\ref{fig:octagon}.

\begin{lemma}
\label{lem:octagon-model}~
In the notation above,
\begin{equation*}
	B \cong (O \times F_0) \, \cup_{O \times \partial F_0} T, \qquad
	B \cap Y \cong \left((\mathsf{b} \sqcup \mathsf{b}') \times F_0\right) \cup (T \cap Y)
\end{equation*}
Moreover, as pairs:
\[
	(T, T \cap Y) \cong C\times(B^3,D^2).
\]
In particular, $H_*(T, T \cap Y; \mathbb{Z})=0$.
\end{lemma}

\begin{proof}
	Identify an interior region of $S \subset ([0,1]_t \times Y')$ with $F_0$.
	The normal coordinates are the time $t$ and the normal coordinate $s$ of $F_0 \subset Y'$.
	Thus the tubular neighbourhood of $\mathrm{int}(S)$ can be identified with $P \times F_0$,
	where the branch locus is $(1/2,0) \times F_0$,
	with branch cut $[0,1/2] \times \{0\} \times F_0$.
	Hence the covering $\pi_W$ is modelled on $(O \times F_0) \to (P \times F_0)$,
	with deck transformation given by a half turn of the octagon.
	The incoming part of $B$ is precisely the preimage over $\{t=0\} \subset P$ in this model, which is $(\mathsf{b} \sqcup \mathsf{b}') \times F_0$.

	As for the collar, fix a point of $C$,
	and view its normal 3-ball in $[0,1]\times Y'$ as $I \times D^2$.
	Before rounding,
	the branch locus is
	\begin{equation*}
		\{(t,0,0) : 0 \le t < 1/2 \} \cup
		\{(1/2, r, 0) : 0 \le r \le \epsilon\}.
	\end{equation*}
	Rounding turns this arc into a boundary-parallel arc joining $(0,0,0)$ to the side boundary at $t=1/2$.
	By straightening the arc, the branched cover of this 3-ball can be seen to be a 3-ball.
	The preimage of the incoming $D^2$ downstairs of this cover is again a single disc $D^2 \subset B^3$.
	These local models along $C$ fit together as a product globally using the given Seifert framing, which trivializes the normal bundle $N'(C)$.
	Indeed, under this trivialization, the branch cut is constant along $C$.
	This proves $(T, T \cap Y)\cong C \times (B^3, D^2)$.
\end{proof}

The attaching of $B$ onto $A$ can be described more explicitly:
Write $$\Sigma_i = F_{i-} \cup_C F_{i+},  \ \ i = 1,2,$$
where $F_{i\pm}$ is a copy of $F$.
The four remaining sides of $O$ are $\mathsf{f}_{i\pm}$,
with $\mathsf{f}_{1+}$ between $\mathsf{t}$ and $\mathsf{b}$, $\mathsf{f}_{2-}$ between $\mathsf{b}$ and $\mathsf{t}'$, $\mathsf{f}_{2+}$ between $\mathsf{t}'$ and $\mathsf{b}'$, and finally,
$\mathsf{f}_{1-}$ between $\mathsf{b}'$ and $\mathsf{t}$.
Attach the face $\mathsf{f}_{i\pm} \times F_0$ to the corresponding copy of $F_0 \subset F_{i\pm}$.
The remaining collar faces arise from $T$.
Together, these form the attaching region $[0,1] \times (\Sigma_1 \sqcup \Sigma_2)$.

\subsection{Relative topology and extensions}
Let us analyze the relative topology of $(W_F, Y)$, using excision of cohomology.
Let $\mathsf{b},\mathsf{b}'\subset\partial O$ denote the two incoming arcs in the octagon model, and $\mathsf{t},\mathsf{t}'$ denote the two outgoing arcs.
See Figure~\ref{fig:octagon} again.
Write $$\Lambda = H^1(F,\partial F;\mathbb{Z})$$ and $I = [0,1]$.

\begin{proposition}
\label{prop:WF-cohom}
For every $k$,
\begin{equation}\label{eq:WF-cohom}
 H^k(W_F,Y;\Z)\cong H^{k-1}(F,\partial F;\Z).
\end{equation}
In particular,
\[
 H^1(W_F,Y; \mathbb{Z})=0,\qquad H^2(W_F,Y; \mathbb{Z})\cong\Lambda,\qquad
 H^3(W_F,Y; \mathbb{Z})\cong\Z.
\]
The deck involution acts by $-1$ on every group in
\eqref{eq:WF-cohom}.
\end{proposition}

\begin{proof}
Recall $W_F=A\cup B$, where $A=I\times(X_1\sqcup X_2)$, with each $X_i$ a copy of $\overline{Y'\setminus N'(F)}$,
and $B$ is the union of $O \times F_0$ and the tube piece $T$.
Both $(A,A\cap Y)$ and $(A\cap B,A\cap B\cap Y)$ are pairs of the form
$([0,1]\times Z,\{0\}\times Z)$, which has trivial $H^*$.
Relative Mayer--Vietoris therefore gives
\[
 H^*(W_F,Y)\cong H^*(B,B \cap Y).
\]
By Lemma \ref{lem:octagon-model}, the tube part satisfies $H^*(T, T\cap Y) = 0$.
Write $\Bin = B \cap Y$.
We have a long exact sequence of triples:
\begin{equation*}
\cdots \to
H^{k-1}(\Bin \cup T,\Bin)
\to H^k(B, \Bin\cup T)
\to H^k(B,\Bin)
\to H^k(\Bin \cup T,\Bin)
\to \cdots.
\end{equation*}
Since $\Bin \cup T$ deformation retracts onto $\Bin$, we have $H^*(\Bin \cup T,\Bin)=0$ so
$$ H^k(B,\Bin) \cong H^k(B, \Bin\cup T).$$
Excision says:
\begin{equation*}
H^*(B,\Bin\cup T) \cong
H^*\!\left(O \times F_0,\, ((O \times F_0)\cap \Bin)\cup((O \times F_0)\cap T)\right).
\end{equation*}
By Lemma \ref{lem:octagon-model}, the latter pair is
\begin{equation*}
	\left(O \times F_0,\ \left((\mathsf{b}\sqcup \mathsf{b}') \times F_0\right) \, \cup \, O \times \partial F_0 \right) \cong
	(O, \mathsf{b}\sqcup \mathsf{b}') \, \times \, (F, \partial F).
\end{equation*}
Now
\[
 H^k(O,\mathsf{b}\sqcup \mathsf{b}';\Z)=
 \begin{cases}
 \Z,&k=1,\\
 0,&k\ne1.
 \end{cases}
\]
The connecting map $H^0(\mathsf b\sqcup\mathsf b') \to H^1(O,\mathsf b\sqcup\mathsf b')$ identifies $H^1(O,\mathsf b\sqcup\mathsf b';\mathbb Z)$ with $\mathbb Z^2/\langle(1,1)\rangle$,
generated by the class of the constant function $(1,0)$.
The deck transformation $\tau$ interchanges the arcs $\mathsf{b}$ and $\mathsf{b}'$,
and sends the generator to $[(0,1)]=-[(1,0)]$.
Then \eqref{eq:WF-cohom} follows from the relative K\"unneth theorem.
Finally, since $\tau$ acts trivially on the $F$-factor while negating the $H^1$ of the $(O,\mathsf{b} \sqcup \mathsf{b}')$-factor,
it acts by $-1$ on $H^2(W_{F},Y) \cong \Lambda$.
\end{proof}

\begin{remark}
\label{rem:arc-a}
The generator of $H^1(O,\mathsf{b} \, \sqcup \, \mathsf{b}')$ is dual to an arc $a_O$ from $\mathsf{t}$ to $\mathsf{t}'$; see Figure~\ref{fig:octagon}.
Under \eqref{eq:WF-cohom},
every class $\theta\in H^2(W_{F},Y)$ has Poincar\'e--Lefschetz dual represented in
the octagon region $B \subset W_{F}$ by the 2-chain
\[
 a_O\times\gamma_\theta,
\]
where $\gamma_\theta\subset F$ is a multicurve representing $\PD_F(\theta)$.
The deck involution reverses
the orientation of $a_O$ but preserves the orientation of $\gamma_{\theta}$.
This explains geometrically why every element of $H^2(W_F,Y)$ is $\tau^*$-anti-invariant.
\end{remark}

\begin{lemma}
\label{lem:connecting-map}
The connecting map $\partial_2 \colon H^2(Y;\Z) \to H^3(W_F,Y;\Z)$ vanishes.
\end{lemma}

\begin{proof}
Under Poincar\'e--Lefschetz duality, the generator
$\Omega\in H^3(W_F,Y;\Z)\cong\Z$ can be represented by
$a_O\times\{x\}\in H_1(W_F,Y' \sqcup Y';\Z)$, for $x\in\operatorname{int}F$.
Naturality of Poincar\'e--Lefschetz duality identifies the map
\[
 j \colon H^3(W_F,Y;\Z)\longrightarrow H^3(W_F;\Z)
\]
with $H_1(W_F,Y' \sqcup Y';\Z)\to H_1(W_F,\partial W_F;\Z)$.
The image of $(a_O \times \{x\})$ has boundary $[y_2]-[y_1]$,
with $y_1,y_2$ the two endpoints of the arc,
lying in distinct components of $Y' \sqcup Y'$.
But the class $ [y_2] - [y_1] \in H_0(\partial W_F;\Z)$ has infinite order.
Thus $j(\Omega)$ also has infinite order and $j$ has no kernel.
By exactness, $\partial_2 = 0$.
\end{proof}
\begin{lemma}

Every ordinary $\Spinc$ structure $\mathfrak{s}$ on the boundary $Y$ extends over $W_F$.
The set of $\Spinc$ extensions $\mathfrak{t}$ with a fixed boundary identification $\mathfrak{t}|_Y \cong \mathfrak{s}$ is a torsor over
$H^2(W_F,Y;\Z)\cong\Lambda$.
\end{lemma}

\begin{proof}
Since the connecting map vanishes by Lemma \ref{lem:connecting-map}, the restriction map
$H^2(W_F;\Z)\to H^2(Y;\Z)$
is surjective.
This implies that starting from any $\Spinc$ structure $\mathfrak{t}$ on $W_F$,
one can twist it by a line bundle $L$ on $W_F$,
so that $\mathfrak{t} \otimes L$ restricts to the given \text{spin\textsuperscript{c}} structure on $Y$.
The remaining statements follow from Proposition~\ref{prop:WF-cohom}.
\end{proof}

\begin{proposition}
\label{prop:real-extension-WF}
Let $(\mathfrak s,\sr)$ be a real $\Spinc$ structure on $Y$.
Every ordinary extension of $\mathfrak{s}$ over $W_F$ rel $Y$ admits a real structure extending $\sr$. Further, this extended real structure is unique up to isomorphism rel $Y$.
Thus real extensions rel $Y$ form a $\Lambda$-torsor, and the forgetful map is a
bijection:
\[
 \{\text{real \text{spin\textsuperscript{c}} extensions over $(W_{F},\tau)$ rel }(Y,\tau)\}\xrightarrow{\ \cong\ }
 \{\text{ordinary \text{spin\textsuperscript{c}} extensions rel }Y\}.
\]
\end{proposition}

\begin{proof}
Over the $\Lambda$-torsor of ordinary extensions rel $Y$,
the involution $\tau$ on $W_{F}$ defines an affine involution
$\kappa:\tf \mapsto \tau^*\overline{\tf}$.
Its linear part on $H^2(W_F,Y)$ is $\theta\mapsto-\tau^*\theta$,
which by Proposition~\ref{prop:WF-cohom} is the identity.
Thus $\kappa$ is translation by some $\eta \in\Lambda$.
Since pulling back twice is the identity, one has $2\eta=0$.
But then $\eta = 0$ since $\Lambda$ is torsion-free.
Thus every $\mathfrak{t}$ is isomorphic to $\tau^* \bar{\mathfrak{t}}$.

Choose an isomorphism $I \colon \mathfrak{t} \to \tau^*\bar{\mathfrak{t}}$ extending $\sr$.
Its square is not necessarily the identity,
but $I^2$ is always a gauge transformation $g\colon W_F\to S^1$,
for which $g|_Y \equiv 1$ and $g(x)=\overline{g(\tau x)}$.
Since $H^1(W_F,Y)=0$, one can write $g=e^{if}$ for which $f|_Y=0$.
Connectedness gives $f\circ\tau=-f$.
Then $e^{-if/2}I$ defines a genuine real structure.
This proves existence of real \text{spin\textsuperscript{c}} extensions.
Indeed, observe the difference of two real structures (over the same extension)
is a $\tau$-\emph{invariant} map $h \colon W_F\to S^1$ which, by $H^1(W_{F},Y; \mathbb{Z}) =0$, must be of the form $h=e^{if}$ for some $\tau$-invariant $f$ satisfying $f|_Y = 0$.
Thus $e^{-if/2}$ identifies the two real structures rel $Y$.
\end{proof}

\begin{definition}
\label{def:Ltheta}
Given $\theta \in \Lambda = H^1(F, \partial F; \mathbb{Z})$, using the identification of $\Lambda$ with $H^2(W_F, Y; \Z)$ from Proposition~\ref{prop:WF-cohom}, we let $L_{\theta}$ be the corresponding line bundle on $W_F$ rel $Y$. Furthermore, given  an extension $\mathfrak{t}$ on $W_F$,
let $\mathfrak{t}^{\theta}$ be the extension \spinc\ structure on $W_F$ further twisted by $\theta$.
That is, $\mathfrak{t}^{\theta} = \mathfrak{t} \otimes L_{\theta}$.
\end{definition}

\begin{figure}[htb]
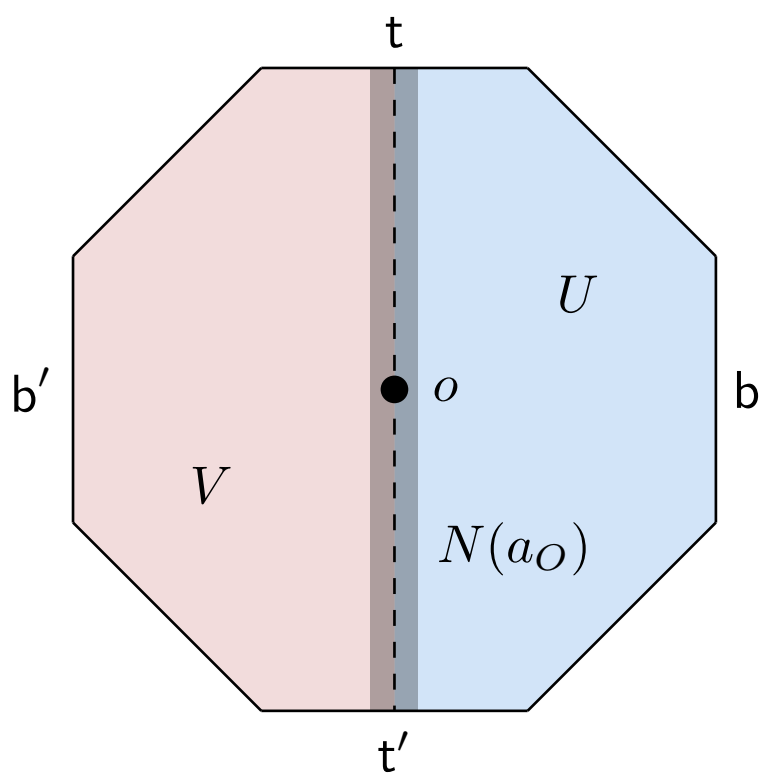

\centering
\figCharts
\caption{
The bundle $L_\theta$ is trivial over the cylinder together with $U\times F$ and over $V\times F$.
The clutching function $g_{\theta}$ is defined on the strip $N(a_O)\times F$.}
\label{fig:charts}
\end{figure}

Let us describe $L_\theta$ more concretely.
Represent $\theta \in H^1(F, \partial F; \mathbb{Z})$ by a map $g_{\theta} \colon F \to S^1$ for which
 $g_{\theta}$ is equal to one on a collar neighbourhood of $\partial F$.
Cut $O$ along the arc $a_O$ into two regions $U$ and $V$, as in~Figure~\ref{fig:charts}.
In particular, $\tau$ reverses the arc $a_O$, fixes the midpoint $o$,
and interchanges $U$ and $V$.
Take trivial complex line bundles on $U \times F$ and $V \times F$,
glued along the overlap $a_O \times F$ via
the clutching function $g_{\theta}$ pulled back from $F$:
\[(x,z)_U\sim(x,g_\theta(x)z)_V\]
whenever $(x,z) \in (a_O \times F) \times \mathbb{C}$.
The assumption $g_\theta=1$ near $\partial F$ ensures $L_{\theta}$ fits with the trivial bundle on the remaining pieces of $W_F$.
Moreover, this construction gives a specified trivialization of $L_\theta$ near the incoming boundary $Y$.
Let $\gamma_{\theta} \in H_1(F;\mathbb{Z})$ be the Poincar\'{e}--Lefschetz dual of $\theta$.
Then
\begin{equation*}
	c_1(L_\theta) = \operatorname{PD}[a_O\times\gamma_\theta] \in H^2(W_{F},Y;\mathbb{Z}),
\end{equation*}
cf.~Remark~\ref{rem:arc-a}.

Further, the complex line bundle $L_{\theta}$ admits a natural real structure.
Indeed, define conjugate-linear lifts $\sr_{\theta}$:
\[
 \sr_\theta(x,z)_U=(\tau x,\overline z)_V,
 \quad
 \sr_\theta(x,z)_V=(\tau x,\overline z)_U.
\]
On the remaining disconnected regions of $W_{F}$ that are freely interchanged by $\tau$,
let $\sr_{\theta}$ be the fibrewise complex conjugation.
Then
\[
 \sr_\theta(x,g_\theta(x)z)_V
 = (\tau x,\overline{g_\theta(x)}\,\overline z)_U
 \sim
 (\tau x,g_\theta(x)\overline{g_\theta(x)}\,\overline z)_V
 = (\tau x,\overline z)_V
 =\sr_\theta(x,z)_U.
\]
The real structure $\sr_{\theta}$ agrees with the fibrewise conjugation real structure in the specified trivialization near $Y$.
Hence tensoring with the real line bundle $(L_{\theta},\sr_\theta)$ preserves the boundary real structures of any real $\Spinc$ structure rel $Y$.

\subsection{Outgoing caps and line bundles}\label{sec:index-form}

\begin{definition}

An \emph{admissible} outgoing cap is a compact connected oriented 4-manifold $M$ with
$\partial M\cong -Y'$ and $H_1(M;\Z)=0$.
The two outgoing components $Y' \sqcup Y'$ of $\partial W_F$ are capped by two copies of $M$,
where we extend the involution $\tau$ on $W_{F,M}$ via exchanging the two components of $M \sqcup M$.
\end{definition}

Given any $Y'$,
an admissible outgoing cap $M$ always exists.
Indeed, by the Lickorish--Wallace theorem~\cite{Lickorish, Wallace},
there exists an integral surgery presentation of $-Y'$.
Let $M$ be the resulting surgery trace of $-Y'$, which satisfies $\partial M=-Y'$, and has no 1-handles.
Then $H_1(M;\Z)=0$.

\begin{lemma}
\label{lem:admissible-cap}
Suppose $M$ is an admissible outgoing cap.
Every $\Spinc$ structure and every line bundle on $Y'$ extend over the outgoing caps $M$.
Moreover, every integral one-cycle in $Y'$ bounds an oriented surface in $M$.
\end{lemma}

\begin{proof}

Consider the following portions of the (co)homology exact sequences
\begin{equation}
\label{eq:LES_outgoing_cap_co}
H^2(M;\Z)\longrightarrow H^2(Y';\Z)
\longrightarrow H^3(M,Y';\Z)\cong H_1(M;\Z)
\end{equation}
and
\begin{equation}
\label{eq:LES_outgoing_cap}
 H_2(M,Y';\Z)\longrightarrow H_1(Y';\Z)
 \longrightarrow H_1(M;\Z).
\end{equation}

Since $H_1(M;\Z)=0$, the first map in~\eqref{eq:LES_outgoing_cap_co} is surjective.
Hence every line bundle on $Y'$ extends to $M$.
This implies that for any $\Spinc$ structure $\widehat{\mathfrak{t}}$ on $M$, there exists a line bundle $\widehat{L}$ such that $\widehat{\mathfrak{t}} \otimes \widehat{L}$ has the correct prescribed $\Spinc$ structure on $Y'$.
Similarly, the first boundary map in~\eqref{eq:LES_outgoing_cap} is surjective and so every one-cycle bounds an element of $H_2(M,Y';\Z)$.
Such a relative homology class in $H_2(M,Y';\Z)$ can be represented by connected embedded oriented surfaces after resolving double points and adding tubes.
\end{proof}

Fix an admissible $M$. Let
\begin{equation*}
	W_{F,M} := W_{F} \cup_{Y' \sqcup Y'} \left(M \sqcup M\right)
\end{equation*}
be obtained by capping the two outgoing components of $W_F$ with the exchanged copies. (See Figure~\ref{fig:quotient}.)
Every real $\Spinc$ structure $\mathfrak{t}$ over $W_F$ rel $Y$ therefore extends over $W_{F,M}$.
Indeed,
choose any $\Spinc$ extension of $\mathfrak{t}|_{Y'}$ on one copy of $M$,
and use the conjugate extension on the other.
This identification with the conjugate copy defines a real structure on $M \sqcup M$.
We refer to such a real $\Spinc$ extension on $M \sqcup M$ as \emph{compatible capping data}.

\begin{corollary}

Under the hypotheses that $Y \to Y'$ is induced by a Seifert surface,
$(Y,\tau;\mathfrak s,\sr)$ is zero in the Real $\Spinc$ bordism group.
\end{corollary}
\begin{proof}
	Proposition~\ref{prop:real-extension-WF} and Lemma~\ref{lem:admissible-cap} together provide a real $\Spinc$ structure $(\widehat{\mathfrak{t}}, \widehat{\sr})$ on $W_{F,M}$, and $\partial W_{F,M} = -Y$.
\end{proof}

Fix a real $\Spinc$ structure $\widehat\tf$ over $W_{F,M}$ (rel $Y$) as above.
Choose extension line bundles $\{\widehat{L}_{\beta}\}$ over $W_{F,M}$ of the relative line bundles $\{L_{\beta}\}$ (as in Definition~\ref{def:Ltheta}),
where $\beta$ is in any integral basis of $\Lambda$.
We arrange it so that $\tau$ pulls back $\widehat{L}_{\beta}|_M$ on one of the outgoing cap $M$ to its conjugate line bundle on the other copy of $M$.
Extend this family from the chosen basis to $\{\widehat{L}_{\theta}\}_{\theta \in \Lambda}$ of line bundles by tensor products.
The resulting family is \emph{coherent}
in the sense that:
\[
 \widehat L_\theta|_{W_F}=L_\theta,\qquad
 \widehat L_{\theta+\eta}\cong\widehat L_\theta\otimes\widehat L_\eta.
\]
Let
\[
 x_\theta=c_1(\widehat L_\theta)\in H^2(W_{F,M},Y;\Z),\qquad
 \widehat\tf^{\theta}=\widehat\tf\otimes\widehat L_\theta.
\]
Then $x_{\theta+\eta}=x_\theta+x_\eta$.

We can choose an integral basis $\{\beta\}$ of $\Lambda$ so that each $\PD_F(\beta)$ can be represented by an oriented embedded curve $\gamma_{\beta}$.
Then the dual of a basis class $x_{\beta}$ can be represented by the two-cycle $a_O\times\gamma_\beta$ capped on $Y'$ by exchanged surfaces in $M \sqcup M$, with orientations reversed by $\tau$.
Here, the oriented capping surfaces exist by the last assertion of Lemma~\ref{lem:admissible-cap}.
The rest of the classes in $H^1(F, \partial F;\mathbb{Z})$ can be obtained additively.

Furthermore, each $L_{\theta}$ is equipped with a real structure $\sr_{\theta}$, unique relative to $Y$.
Using the same basis $\beta$ above, we can extend $\{\sr_{\theta} : \theta \in \Lambda\}$ across $M \sqcup M$ as above.
Thus $\{\widehat{L}_{\theta}\}$ can be upgraded to a coherent family of real line bundles on $W_{F,M}$.

\subsection{Normalizing the APS index}

Let $(W,\tau)$ be a compact oriented 4-manifold with $\partial W = -Y$,
equipped with an orientation-preserving involution $\tau$,
and a $\tau$-invariant metric of product type near the boundary $Y$.
Regard $Y$ as the incoming boundary.

Let $\mathfrak a=(B,0)$ be a reducible configuration on $(Y,\tau,\mathfrak s,\sr)$ for which the Dirac operator $D_B$ is invertible.
Let $\gr_{SW}(\mathfrak a;W, \tf)$ be the index of the linearized Seiberg--Witten operator on $W$,
with a $\Spinc$ structure $\tf$ extending $\mathfrak s$, with fixed boundary data at $\mathfrak a$,
and with spectral boundary conditions as in \cite[\S20.3]{KM} (see~\cite{Li26} for the framed counterpart).
In particular,
the spectral projections are defined using $\mathrm{Hess}_{0,\mathfrak{a}}-\epsilon$, for sufficiently small $\epsilon>0$, where $\mathrm{Hess}_{\mathfrak{a}}$ denotes the Hessian of the CSD functional in the Coulomb slice.

With these spectral boundary conditions on the form part,
let $- \iota_R(W)$ be the index of the ASD operator:
\[
(d^*\oplus d^+)_R\colon
\Omega^1(W;i\R)^{-\tau^*}\longrightarrow
\Omega^0(W;i\R)^{-\tau^*}\oplus\Omega^+(W;i\R)^{-\tau^*}.
\]
This is the real analogue of
\cite[Equation (25.13)]{KM}.
When $Y$ is connected, this recovers the formula $\iota_R(W)=b_R^+(W)-b_R^1(W)$ in~\S\ref{sec:knot_idea}.
Choose a $\Spinc$ connection $A$ on $W$ extending $B$.
Then at the configuration $(A,0)$, the unperturbed linearized Seiberg--Witten operator with gauge fixing splits as
\[
	\mathcal D \mathfrak{F}
	=(d^*\oplus d^+)_R\oplus D^{+,\sr}_{A}.
\]

We now specialize to $W=W_{F, M}$ as in the previous subsections.

\begin{lemma}\label{lem:APS-cap-change}
Fix $\mathfrak{a}$  and a real $\Spinc$ structure $\tf$ over $W_F$ rel $Y$, as above.
Suppose $M_1,M_2$ are  two admissible outgoing caps,
with $\Spinc$ extensions on $M_i \sqcup M_i$ that are real with respect to $\tau$.
Let $\widehat{\mathfrak t}_i$ be the resulting extensions on the caps $W_{F,M_i}$, for $i=1,2$.
Then
\begin{equation}
 \gr_{SW}(\mathfrak a;W_{F,M_1}, \widehat{\tf}_1)
 -\gr_{SW}(\mathfrak a;W_{F,M_2}, \widehat{\tf}_2)
 \equiv\iota_R(W_{F,M_2})-\iota_R(W_{F,M_1})\bmod2.
\end{equation}
\end{lemma}

\begin{proof}
By the splitting of $\mathcal D \mathfrak{F}$ as the de Rham part and Dirac part as above, and the definition of $\iota_R$:
\begin{equation}\label{eq:APS-normalization}
 \gr_{SW}(\mathfrak a;W_{F,M},\tf)=-\iota_R(W_{F,M})
  +\ind_{\C}^{APS}D^+_{\widehat\tf}.
\end{equation}
The differences of the Dirac indices on the right side of \eqref{eq:APS-normalization} can be expressed using the index gluing theorem as the real index of the Dirac operator over two identical closed 4-manifolds with conjugate $\Spinc$ data.
Thus this index must be even.
\end{proof}

The preceding lemma shows that the index of the Dirac operator modulo two is independent of $M$ and the real $\Spinc$ data on $M$:
\[
\ind_{\C}^{APS}D^+_{\widehat\tf} = \gr_{SW}(\mathfrak a;W_{F,M}, \tf)+\iota_R(W_{F,M})\bmod2.\]
\begin{definition}
	For an arbitrary configuration $\mathfrak{a}$, possibly with perturbation $\mathfrak{q}$, the \emph{normalized APS index} is
	\[
	\mathcal{I}_{\mathfrak{t}}(\mathfrak{a})=
	\gr_{SW}(\mathfrak a;W_{F,M}, \widehat{\tf})+\iota_R(W_{F,M})\bmod2.\]
	Here we suppress the perturbation $\mathfrak{q}$ from the notation,
	although $\mathfrak{q}$ is used in the APS boundary condition.
\end{definition}

\subsection{The Dirac index mod-two}\label{ss:index-dirac} We now investigate the dependence of the Dirac index on $\theta$.
Fix the same APS condition for all $\widehat{\mathfrak t}$'s, with invertible boundary Dirac operator.
Let
\[
 P_{\tf}(\theta)=
 \ind_{\C}^{APS}D^+_{\widehat\tf^\theta} - \ind_{\C}^{APS}D^+_{\widehat\tf} \in\Z.
\]
\begin{lemma}\label{lem:Q-cap-independent}
The parity $Q_{\tf}(\theta) := P_{\tf}(\theta)\bmod2$ is independent of the admissible cap $M$, of the
extensions $(\widehat{\tf},\widehat{L}_\theta)$ of $(\tf,L_{\theta})$ from $W_F$ to $W_{F,M}$, and of the common APS boundary condition on $Y$ given by $\mathfrak{a}$.
\end{lemma}

\begin{proof}
Independence of the cap $M$ and of the extensions follows from Lemma~\ref{lem:APS-cap-change}. Changing the common boundary condition changes both $\ind^{APS}D^+_{\widehat\tf^\theta}$ and $\ind^{APS}D^+_{\widehat\tf}$ by the same relative index, since the line bundle $L_{\theta}$ is trivialized on a collar of $Y$.
\end{proof}

\begin{lemma}
For $\theta,\eta\in\Lambda$,
\begin{equation}\label{eq:index-difference}
 \ind_{\C}^{APS}D^+_{\widehat\tf^\theta} - \ind_{\C}^{APS}D^+_{\widehat\tf^\eta}
 =\frac12\left(x_\theta^2-x_\eta^2
       +c_1(\widehat\tf)\cdot(x_\theta-x_\eta)\right).
\end{equation}
In particular,
\begin{equation}\label{eq:P-topological}
 P_{\tf}(\theta)=\frac12\left(x_\theta^2+c_1(\widehat\tf)\cdot x_\theta\right).
\end{equation}
The products above use the relative intersection pairing and the absolute--relative pairing on
$W_{F,M}$, respectively.
\end{lemma}

\begin{proof}
The boundary terms in the APS formula~\cite{APS} cancel.
If $\omega_\theta$ represents $c_1(\widehat\tf^\theta)$,
the remaining term in~\eqref{eq:index-difference} is
\[
 \frac18\int_{W_{F,M}}(\omega_\theta+\omega_\eta)
                         \wedge(\omega_\theta-\omega_\eta).
\]
The first factor above represents
$2c_1(\hat\tf)+2j(x_\theta+x_\eta)$,
where $j \colon H^2(W_{F,M},Y)\to H^2(W_{F,M})$ is the natural map.
The second factor is supported in the interior and represents $2(x_\theta-x_\eta)$.
Thus the integral above is
\[\tfrac12(c_1(\widehat\tf)+x_\theta+x_\eta)\cdot(x_\theta-x_\eta). \qedhere\]
\end{proof}

\begin{remark}
	We have obtained all the ingredients to complete the well-definedness of the absolute mod-two grading in the knot case,
	except for Proposition~\ref{prop:strong-comparison} which relates $Q_{\mathfrak{t}}$ to a quadratic form on $F$.
	The remainder of this paper is devoted to the case when $|C| > 1$.
	This involves a shift of perspective (to spin structures),
	which can be of independent interest,
	and some amount of bookkeeping.
\end{remark}

\begin{remark}
\label{rem:radical-and-Q}
To illustrate the issue in the multi-component case, suppose some $\theta$ dual to $[C_k]$ satisfies $Q_{\mathfrak{t}}(\theta)=1$.
Then
\[Q_{\tf^\theta} - Q_{\tf} = B(\theta, \cdot \,) \equiv 0,\]
by Proposition~\ref{prop:strong-comparison}.
But
\[
 \ind^{APS}_{\C}D^+_{\widehat{\mathfrak{t}}^{\theta}} - \ind^{APS}_{\C}D^+_{\widehat{\mathfrak{t}}} \equiv Q_{\mathfrak{t}}(\theta) = 1 \bmod 2.
\]
Hence the quadratic form $Q_{\mathfrak{t}}$ alone may not sense the index parity change.
This phenomenon appears in~\S\ref{sec:annulus}.
\end{remark}

\section{Spin structures on fixed-point sets}
\label{sec:spin}
In this section, we relate the index calculation in \S\ref{ss:index-dirac} to Arf invariants of spin structures on the Seifert surface $F$.

\subsection{Spin structures on fixed surfaces}
\label{sec:framing}
Suppose $Z$ is \emph{any} embedded oriented surface in an oriented 4-manifold $X$,
with oriented normal bundle $N_Z$.
Since $\mathit{SO}(2) \cong U(1)$, we may regard $TZ$, $N_Z$, and hence $TX|_{Z}$ as complex vector bundles.
Denote the restriction to $Z$ of the anti-canonical bundle of $X$ by
\begin{equation}
 \mathcal A_Z:=(TZ)_{\C}\otimes(N_Z)_{\C}
 =\det{}_{\C}(TX)\bigr|_{Z}.
\end{equation}
Assume $(X,\tau)$ is a real \text{spin\textsuperscript{c}} 4-manifold,
with spinor bundles $S^{\pm}$, Clifford multiplication $\rho \colon TX \to \operatorname{Hom}(S^\pm,S^{\mp})$, and $\rho$-compatible real structure $\sr \colon S^{\pm} \to S^{\pm}$.
\begin{lemma}\label{lem:half4-new}
Suppose the fixed surface $Z=\Fix(\tau)$ is oriented.
Orient $N_{Z}$ so that the induced orientation on $TX|_{Z} = TZ \oplus N_{Z}$ is positive.
Consider the complex structure $J_N=\rho(n_1)\rho(n_2)$,
for any positively oriented frame $(n_1,n_2)$ of $N_Z$.
Then
\[
 S^+|_Z \cong E_- \oplus E_+,
\]
where $E_{\pm}$ is the $(\pm i)$-eigenline bundle of $J_N$. Then:
\begin{enumerate}
\item[(a)]  $E_+ \cong E_-\otimes\mathcal A_Z$;
\item[(b)] The following is an $\mathbb{R}$-linear isometry:
\begin{equation}
	\Phi \colon E_-\longrightarrow(S^+)^{\sr},\quad
	\Phi(\psi)=\frac{\psi+\sr\psi}{\sqrt2}.
\end{equation}
\item[(c)] The real structure $\sr$ defines isomorphisms
\begin{equation}\label{eq:half4-square}
E_+\cong\overline{E_-}\cong E_-^{-1}, \qquad
 E_-^{\otimes2}\cong\mathcal A_Z^{-1}.
\end{equation}
\item[(d)] A normal framing of $Z$ trivializes $(N_Z)_{\C}$ and so it identifies $\mathcal A_Z\cong(TZ)_{\C}$;
further, it defines an isomorphism
\begin{equation*}
	\mu \colon E_-^{\otimes2}\longrightarrow(TZ)_{\C}^{-1}.
\end{equation*}
Hence $E_-$ is a spin structure on $Z$, i.e., a square root of the canonical bundle.
\end{enumerate}
\end{lemma}

\begin{proof}
The map
\[
 (TZ)_{\C}\otimes(N_Z)_{\C}\otimes E_-\longrightarrow E_+,
 \qquad x\otimes n\otimes\psi\longmapsto\rho(x)\rho(n)\psi,
\]
is complex-linear in each factor.
Thus $E_+\cong E_-\otimes\mathcal A_Z$, as claimed in (a).
Since both normal vectors are
anti-invariant, $\sr$ commutes with $J_N$.
Moreover, $\sr$ exchanges the $+i$ and $-i$ eigenlines.
Thus $\Phi$ is an isometry, as claimed in (b).
Moreover, $E_+\cong\overline{E_-}$,
and by the unitary structure $\overline{E_-} \cong E_-^{-1}$.
Comparison of these two descriptions of $E_+$, together with part (a), yields
\eqref{eq:half4-square}. This proves (c). Part (d) follows immediately from the description of $ E_-^{\otimes2}$ in (c).
\end{proof}

Let $F$ be the oriented Seifert surface in $Y'$,
and $e_1$ be the unit tangent of $C$ determined by the orientation of $C$.
Denote by $e_3'$ the outward unit normal along $C$ in $F$,
and choose the unit normal frame $n_F$ of $F$ so that the frame
$(e_1,n_F,e_3')$ is positively oriented in $Y'$.
Choose the positive product normal framing over $F_0$ in the product region $O\times F_0$,
and extend it across the rounding and along $[0,1/2] \times C$ to agree with $(n_F,e_3')$ near $\{0\}\times Y'$.
Denote by $\uptau'$ the resulting normal framing of the branch surface $S$.

To apply Lemma~\ref{lem:half4-new} to $S \subset W_F$,
we construct a normal framing of $S$ upstairs as follows.
Since every positive normal frame downstairs lifts to two positive normal frames upstairs,
let $\mathscr{T} \to S$ be the double cover of $S$,
whose fibre consists of the positive normal frames of $S \subset W_F$
that descend to $\uptau'$ downstairs. The cover $\mathscr T\to S$ is trivial, because the $n_F$-push-off of $F$ can be rotated through the complement of the branch set to the trivially covered region $t > \tfrac12$.

Fix one of two normal framings $\uptau$ of $S$ that lift $\uptau'$,
and let $\mathfrak{t}$ be a real \text{spin\textsuperscript{c}} structure over $W_F$ rel $Y$.
Lemma~\ref{lem:half4-new} gives a spin structure
\begin{equation*}
	\mathfrak{s}_{\mathfrak{t}} \in \mathrm{spin}(F).
\end{equation*}
Here, we are identifying $F$ with $S$.
The two choices $\pm \uptau$ are homotopic through constant rotations of the normal frames,
and define the same normal operator $J_N$ in Lemma~\ref{lem:half4-new}.
Hence $+\uptau$ and $-\uptau$ induce the same $\mathfrak{s}_{\mathfrak{t}}$.

\subsection{Marked real structures and spin lifts in dimension three}
\label{ss:half3}
Recall that an \emph{absolute real $\Spinc$ structure} on $(Y,\tau)$ is a pair $(\mathfrak s,\sr)$ with
$\mathfrak s=(S,\rho)$ a $\Spinc$ structure and $\sr$ a conjugate-linear lift of $\tau$
to $S$ with $\sr^2=1$ and $\sr\rho(\xi)=\rho(\tau_*\xi)\sr$.
An isomorphism of real $\spinc$ structures is a $\Spinc$ isomorphism commuting with the real structures.
A \emph{real gauge transformation} $u \in \mathcal{G}$ is a map
$u\colon Y\to U(1)$ satisfying $u\circ\tau=\bar u$.
In particular, such an $u$ takes values $\pm 1$ on each component of $C$.
The \emph{framed real gauge group} $\mathcal{G}_{\mathbf{p}}$ consists of real gauge transformations equal to $+1$ at the basepoints $\mathbf{p} = \{p_1,\dots,p_n\}$, hence also on $C$ by the assumption that every component of $C$ contains one basepoint.

\begin{definition}
Given a set of basepoints $\mathbf{p} = \{p_1,\dots,p_n\} \subset C$,
a \emph{marked triple} on $(Y, \tau, \mathbf{p})$ is
\[
 \widetilde{\sr}=(\mathfrak s,T,\sr),
\]
where $(\mathfrak s,\sr)$ is a real $\spinc$ structure on $(Y,\tau)$ and $T=(T_k)_{k=1}^n$ are trivializations of the spinor bundle at the basepoints; i.e., isomorphisms of Clifford modules $T_k\colon S_{p_k}\to\C^2$,
satisfying
\[
 T_k\circ \sr=\mathsf{conj}\circ T_k,
\]
where $\mathsf{conj}$ denotes complex conjugation on $\C^2$.
The data of $T$ is a \emph{marking}.
A \emph{marking change} by $a=(a_k)\in \mathbb{Z}_2^n$ is
\[
	a\widetilde{\sr}=(\mathfrak s,a T,\sr), \quad \text{for which} \  (a T)_k=(-1)^{a_k}T_k
\]
Denote the corresponding class of $a$ by $\bar{a}\in \mathscr{A}$, where
\[\mathscr{A}:=\mathbb{Z}_2^n/\mathbb{Z}_2\]
is the quotient by the diagonal $\mathbb{Z}_2$ action (i.e., simultaneous translation of all entries).
An \emph{isomorphism of marked triples}
$u\colon(\mathfrak s,T,\sr)\to(\mathfrak s',T',\sr')$ is an isomorphism of real $\Spinc$ structures with
$T'\circ u=T$.
\end{definition}
\begin{definition}
A \emph{relative real $\Spinc$ structure} (\RRSC\ structure) is an isomorphism class
of marked triples.
Denote by $\RRSC(\mathfrak s,\sr)$
the set of \RRSC\ structures having underlying real $\Spinc$ structure $(\mathfrak s,\sr)$.
\end{definition}

It will be convenient to keep track of the marked triples rather than their equivalence classes as \RRSC\ structures.

In what follows,
we show that a marked triple $(\mathfrak{s},T,\sr)$ on a real 3-manifold $(Y,\tau)$ determines a square root of the unit circle bundle $\mathbb{S}(\nu)$ of the normal bundle $\nu$ of $C \subset Y$.
We work with a convention of Clifford multiplication
$\rho \colon TY \to \operatorname{End}(S)$ that is equivalent to that of~\cite{KM}; denote the Pauli matrices
\begin{equation}
	\label{eq:clifford-in-adapted-frame}
 \sigma_1=\begin{pmatrix}i&0\\0&-i\end{pmatrix},\quad
 \sigma_2=\begin{pmatrix}0&-1\\+1&0\end{pmatrix},\quad
 \sigma_3=\begin{pmatrix}0&i\\i&0\end{pmatrix}.
\end{equation}
Along $C$, let $\nu_F$ be the unit normal of $F$ in $Y$ that pushes forward to $n_F$ under the branched-cover map as in \S\ref{sec:framing}.
Choose $e_3$ so that $(e_1,\nu_F,e_3)$ is a positive orthonormal frame in $Y$.
Then we let
\[\rho(e_1) = \sigma_2, \qquad \rho(\nu_F)=-\sigma_1, \qquad \rho(e_3)=\sigma_3.\]

Over $C$, let $S^{\sr}$ be the real rank-two fixed bundle of $\sr$.
Let $J = -\rho(e_1)$ and, given any unit normal $w\in\mathbb{S}(\nu)$, let
$K_w=\sqrt{-1}\rho(w)$.
By the relations in the Clifford algebra, the $(+1)$-eigenlines of $K_w$'s satisfy the following:
\[
\ker(K_{\cos\varphi\,\nu_F+\sin\varphi\,e_3}-1)
=e^{J\varphi/2}\ker(K_{\nu_F}-1).
\]
Compare with \cite[Lemma~3.1]{Li26}.
Thus each unit vector $u$ of $S^{\sr}$ determines a unique unit normal $\mathsf{w}(u)$ with $K_{\mathsf{w}(u)}u=u$.
Moreover, every normal vector $u \in \mathbb{S}(\nu)$ lifts to two unitary eigenvectors in $S^{\sr}$.
We have obtained the following:

\begin{lemma}\label{lem:half3-new}
	The marked triple $(\mathfrak{s},T,\sr)$ produces a fibrewise double cover of $ \mathbb{S}(\nu)$, consisting of the unit $(+1)$-eigenvectors of $K$:
\[
 U(S^{\sr})\longrightarrow \mathbb S(\nu),\qquad u\longmapsto \mathsf{w}(u).
\]
This defines a square root of $\nu$.
\end{lemma}

\begin{definition}
The set of \emph{spin lifts} of $U(S^{\sr}) \to \mathbb{S}(\nu)$ at the basepoint $p_k \in C_k$ is the two-element set
\[
 \Lift_k=\{u\in U(S^{\sr}_{p_k}):\mathsf{w}(u)=\nu_F(p_k)\}.
\]
\end{definition}

\subsection{Spin markings on the boundary}
\label{ss:spin-marking-on-boundary}
In this section, we incorporate the boundary spin lift $\mathscr{L}$ over $C \subset Y$
from \S\ref{ss:half3} to the spin structure defined in \S\ref{sec:framing}.
To this end, suppose $Z$ is any surface with marked points $(p_1,\dots,p_n)$.
Let $\mathfrak{s}$ be a spin structure on $Z$ and
denote the spin double cover by $P_{\mathfrak{s}} \to \mathbb{S}(TZ)$.
Moreover,
suppose at every $p_k$ we are given a unit tangent vector $v_k \in \mathbb S(T_{p_k}Z)$.
Let $P_{\mathfrak{s},k}$ denote the two-element set of lifts of $v_k$ in $P_{\mathfrak{s}}$.
\begin{definition}

A \emph{spin marking} $\mathsf{m}$ on $(\mathfrak s, (v_k))$ is a choice of lifts of $v_k$ in $P_{\mathfrak{s},k}$ modulo simultaneous sign change, that is, an element of
\[
\mathscr{P}(\mathfrak{s})=
\left(\prod_{k=1}^n P_{\mathfrak s,k}\right)\Big/\{\pm1\}.
\]
A \emph{marked spin structure}  on $(Z,(p_k))$ is a pair $(\mathfrak{s},\mathsf{m})$ consisting of a spin structure $\mathfrak{s}$ and a spin marking $\mathsf{m}$ on $(\mathfrak{s}, (v_k))$.\footnote{The notion of marked spin structure above is related to that of \emph{based $\mathbb F_2$-framings} developed by Kawazumi~\cite{Kawazumi}.
After fixing boundary spin models,
they are equivalent to the relative spin structures in the framework of Randal-Williams~\cite{Randal-Williams}.}
\end{definition}

A marked triple supplies a lift of a unit normal vector rather than a tangent vector as in the definition of a spin marking.
To make $\mathfrak{s}_{\mathfrak{t}}$ a marked spin structure,
we conciliate Lemma~\ref{lem:half3-new} with Lemma~\ref{lem:half4-new}:
Suppose $Z = C \times I$  is a collar of a circle $C$ on a fixed surface in a real 4-manifold $X$, for which $\partial X = Y$.
Lemma~\ref{lem:half3-new} and Lemma~\ref{lem:half4-new} provide two double covers with different base spaces:
\begin{equation}\label{eq:2dbcvs}
	U((S^+)^{\sr}|_{C}) \cong U(S^{\sr}|_{C}) \longrightarrow \mathbb S(\nu), \qquad
	U(E_-) \longrightarrow U(TZ^{-1}|_{C}).
\end{equation}
These two covers are isomorphic by the following lemma,
whose proof is deferred to Appendix~\ref{app:cover-indentification}.

\begin{lemma}
\label{lem:Gamma}
The $\mathbb{R}$-linear isometry $\Phi \colon E_-\to (S^+)^{\sr}$ from Lemma~\ref{lem:half4-new} identifies
the double cover consisting of lifts of $e_1^{-1}$ in $E_-$ in Lemma~\ref{lem:half4-new}
with the double cover consisting of lifts of $\nu_F$ in $S^{\sr}$  in Lemma~\ref{lem:half3-new}.
In particular, the spin lift
$T_k^{-1}(1,0) \in \mathscr{L}_k$ corresponds to a unit vector
$\psi_k\in E_{-,p_k}$ satisfying
$\psi_k^{\otimes2}=e_1(p_k)^{-1}$.
\end{lemma}

Consequently, the two double covers in \eqref{eq:2dbcvs} have the same monodromy around $C_k$,
given by (\cite[\S 3.3]{Li26}):
\[
 w_{1,k} =  \langle w_1\bigl(\nu_F^*U(S^{\sr}|_{C_k})\bigr),[C_k]\rangle.
\]
Indeed, the normal framing $\uptau$ defining $\mathfrak{s}_{\mathfrak{t}}$ and the collar framing $(\nu_F,e_3)$ used in Lemma~\ref{lem:Gamma} are both lifts of $\uptau'$,
so they differ by a constant sign on each $C_k$, which does not affect the spin monodromy.
In particular, $\sum_kw_{1,k} = n \bmod 2$.

Using $U(E_-)\cong U(E_-^{-1})$,
Lemma~\ref{lem:Gamma} identifies the fibre of the fixed-surface spin cover over
$v_k=e_1(p_k)$ with $\Lift_k$.
Define
\begin{equation}\label{eq:rrsc-pointing}
 \varepsilon(\mathfrak{s}_{\mathfrak{t}},\widetilde{\sr})
=[T_1^{-1}(1,0),\ldots,T_n^{-1}(1,0)]\in\mathscr P(\sF_{\tf}).
\end{equation}
The pair $(\mathfrak{s}_{\mathfrak{t}}, \varepsilon(\mathfrak{s}_{\mathfrak{t}},\widetilde{\sr}))$ is a marked spin structure.
Moreover, $\varepsilon(a \widetilde{\sr})=\bar a \cdot\varepsilon(\widetilde{\sr})$,
for any $a\in \mathbb{Z}_2^{n}$ acting by $(-1)^{a_k}$ on each entry.

The two choices of normal framings $\pm \uptau$ on $S$ induce isomorphic marked spin structures.
Indeed, switching from $\uptau$ to $-\uptau$ preserves $E_-$, but turns the square map $\mu$ into $-\mu$.
Hence the map $E_- \to E_-$, $\psi \mapsto i\psi$ induces an isomorphism of spin structures,
which sends the spin marking $[\psi_1,\dots,\psi_n]$ to $[i\psi_1,\dots,i\psi_n]$.
Furthermore, changing the local frame from $(e_1,\nu_F,e_3)$ to $(e_1,-\nu_F,-e_3)$ turns the Clifford trivialization $T_k$ at the basepoint $p_k$ into $-\sigma_2 T_k$.
The resulting spin lift corresponds precisely to $i\psi_k$:
\[(-\sigma_2 T_k)^{-1}(1,0)
=T_k^{-1}(0,1)
=\Phi(i\psi_k).\]
Therefore, $\psi \to i\psi$ is an isomorphism of marked spin structures induced by $\pm \uptau$.

By Proposition~\ref{prop:real-extension-WF},
the \textsc{r}eal line bundles $L_{\theta}$ on $(W_F,Y)$ are parametrized by $\theta \in \Lambda$.
Fix a real \text{spin\textsuperscript{c}} structure $\mathfrak{t}$ rel $Y$.
Consider the spin structure $\mathfrak{s}_{\mathfrak{t}^{\theta}}$ obtained from $\mathfrak{t}^{\theta} = \mathfrak{t} \otimes L_{\theta}$.
By \S\ref{sec:index-form},
we choose a coherent family $\{\widehat{L}_{\theta}\}_{\theta \in \Lambda}$ of \textsc{r}eal line bundles over $W_{F,M} = W_F \, \cup (M \sqcup M)$ rel $Y$ that extend $L_{\theta}$.
The real line bundle $\widehat{L}_\theta$ is canonically trivial on $Y$ by construction in \S\ref{sec:index-form}, with a canonical real section $1$.
Thus for a marked triple $\widetilde{\sr}$ on $Y$, twisting by $\widehat{L}_{\theta}$ does not change the spin marking:
\begin{equation}
 \varepsilon(\mathfrak{s}_{\mathfrak{t}},\widetilde{\sr}) = \varepsilon(\mathfrak{s}_{\tf^\theta},\widetilde{\sr}).
\end{equation}
Let us identify $S \subset W_F$ with $F$.
\begin{lemma}
\label{lem:spin-twist}
On $F$,
\begin{equation}
 \sF_{\tf^\theta}
 =\sF_{\tf}\otimes\ell_\theta,
 \label{eq:spin-twist}
\end{equation}
where $\ell_\theta\to F$ is the real fixed line of $\widehat L_\theta$, with its specified real section $1$ on the boundary collar.
Furthermore, the Stiefel-Whitney class of $\ell_\theta$ relative to this section is
\[
 w_1(\ell_\theta,1)=\theta\bmod2\in H^1(F,\partial F;\F).
\]
\end{lemma}

\begin{proof}
Consider the $\tau$-fixed point set $\{o\}\times F$, which lies in $a_O \times F$. (See~Figure~\ref{fig:charts}.)
At $(o,x)$,
the transition function of $L_{\theta}$ is given by $g_{\theta} \colon (F,\partial F) \to (S^1, 1)$,
while the real structure sends a $U$-fibre $z$ to the conjugate $V$-fibre $\overline z$.
Thus a real fibre $z$ satisfies
\[
 \overline z=g_\theta(x)^{\pm1}z,
\]
and hence
\[
 \ell_\theta|_F=\mathbb R\cdot g_\theta^{\pm1/2}.
\]
Hence $w_1(\ell_\theta,1)=\theta\bmod2$.

To see \eqref{eq:spin-twist},
recall that the spin bundle of $\mathfrak{s}_{\mathfrak{t}}$ is defined (cf.~Lemma~\ref{lem:half4-new})
as the $(-i)$-eigenline bundle $E_-$ of the operator $J_N \otimes 1$. This operator commutes with tensoring by $L_{\theta}$;
therefore
\[
 E_-(\tf\otimes\widehat L_\theta)
 =E_-(\tf)\otimes\widehat L_\theta|_{F}.
\]
Thus $\mathfrak{s}_{\tf^{\theta}} = \mathfrak{s}_{\tf}\otimes \ell_\theta$
and this identification respects the square map $\mu \colon E_{-}^{\otimes 2} \to (TF)^{-1}$.
\end{proof}

\subsection{Index and spin structures}\label{sec:index-geometric}
We begin this subsection with some notation and convention.
Suppose $F$ is an \emph{arbitrary} oriented surface,
possibly with boundary, and equipped with a Riemannian metric.

Let
\[\Lambda = H^1(F,\partial F; \mathbb{Z}) \cong H_1(F;\mathbb{Z}), \qquad V = H_1(F; \F) \cong \Lambda/2\Lambda.\]
Given $\theta \in \Lambda$, let $\gamma_{\theta} \in H_1(F;\mathbb{Z})$ be its Poincar\'{e}--Lefschetz dual.
Denote by $B$ the mod-two intersection form on $F$.
Given $\theta,\eta \in \Lambda$, by conflating them with their dual and using mod-two reduction,
we may write
\[B(\theta, \eta) = \gamma_{\theta} \cdot \gamma_{\eta} \,\bmod 2.\]
The form $B$ is non-degenerate if and only if $F$ has at most one boundary component. In general, let $$\mathrm{rad}(B) = \{\theta  \mid B(\theta,\eta)=0 \text{ for all } \eta \}$$ be the radical of $B$.
Then $\mathrm{rad}(B)$ is generated by the classes of the boundary components of $F$.
Denote $C = C_1 \, \sqcup \cdots \, \sqcup C_n$ the boundary of $F$.

Let $\mathfrak{s}$ be a spin structure on $F$.
Given an embedded oriented circle $\gamma \subset F$,
trivializing $T\gamma$ induces a spin structure $\mathfrak{s}|_{\gamma}$.
Define
\begin{equation}
 q_{\mathfrak s}(\gamma)
 =
\begin{cases}
	0 & \text{if $\mathfrak{s}|\gamma$ is the bounding spin structure},\\
	1 & \text{if $\mathfrak{s}|\gamma$ is the Lie group spin structure.}
\end{cases}
\end{equation}
This extends linearly to disjoint unions.
Johnson~\cite{Johnson} shows that $q_{\mathfrak s}$ is a quadratic refinement of $B$.
That is, for every $x,y \in V$,
\[q_{\mathfrak s}(x+y)=q_{\mathfrak s}(x)+q_{\mathfrak s}(y)+x\cdot y.\]
Assume $B$ is nondegenerate. The \emph{Arf invariant of the spin structure} $\mathfrak{s}$ on $F$ is defined as the Arf invariant
of the quadratic form $q_{\mathfrak{s}}$ on $(V,B)$:
\begin{equation*}
	\Arf(\mathfrak{s}) := \Arf(q_{\mathfrak{s}}).
\end{equation*}
Finally, denote $\lambda = q_{\mathfrak{s}}|_{\rad(B)}$.

\medskip

From now on, let $F$ be a Seifert surface and $W_{F}$ be the Seifert cobordism.
The goal of this subsection is to prove the following identities that relate the index theory on $W_F$ to the spin geometry of $F$.

\begin{proposition}\label{prop:strong-comparison}
For every $\theta\in\Lambda$,
\begin{equation}\label{eq:strong-comparison}
Q_{\mathfrak{t}}(\theta) =
q_{\sF_{\tf}}(\gamma_\theta) \equiv
\frac12\bigl(x_\theta^2+c_1(\widehat\tf)\cdot x_\theta\bigr)\bmod2.
\end{equation}
The equalities hold for all compatible exchanged capping data $M$ and common APS conditions on $Y$.
\end{proposition}

Recall from \S\ref{ss:index-dirac} that $Q_{\mathfrak{t}}$ is defined as the index difference of two Dirac operators over $W_{F,M}$ with identical APS boundary conditions on $Y$, associated to $\mathfrak{t}$ and the twisted $\mathfrak{t}^{\theta}$.
On the other hand, $q_{\mathfrak{s}_{\mathfrak{t}}}$ is the quadratic form defined above,
and the spin structure $\mathfrak{s}_{\mathfrak{t}}$ arises from the $\Spinc$ structure $\mathfrak{t}$.
The third term in \eqref{eq:strong-comparison} relates the two quantities,
where we set $x_{\theta} = c_1(\widehat{L}_{\theta})$.
The proof of this proposition relies on two more lemmas:

\begin{lemma}\label{lem:multicurve-dual}
Given any nonzero $\theta\in\Lambda$,
represent $\PD_F(\theta)$ by a multi-curve $\gamma = \gamma_1 \sqcup \cdots \sqcup \gamma_r$.
Then $\PD(x_\theta)$ can be represented by a connected $\tau$-invariant oriented surface $Z$
for which $\Fix(\tau|_Z)=\gamma$, such that
\[
 [Z]^2\equiv0\bmod2,\qquad
 \frac{\chi(Z)}2=\chi(Z/\tau)\equiv r\bmod2.
\]
Furthermore, along $\gamma$:
\[
	TZ=\langle \mathbf{t},\mathbf{a}\rangle,\qquad N_Z=\langle \mathbf{v},\mathbf{n}\rangle,
\]
where $(\mathbf{v},\mathbf{t})$ is an oriented orthonormal frame of $TF|_{\gamma}$,
and $(\mathbf{n},\mathbf{a})$ is the oriented normal frame of $F$ with $\mathbf{a}$ tangent to
$a_O$.
Over $\operatorname{int} F$, the frame $(\mathbf{n},\mathbf{a})$ differs from the normal frame $\uptau$ by a constant rotation.
\end{lemma}

\begin{lemma}\label{lem:real-eigenline}
Let the eigenline bundle decomposition (Lemma~\ref{lem:half4-new}) of the operator $J_N \colon S^+ \to S^+$ over $Z \subset W_{F,M}$ be
\begin{equation}\label{eq:eigenline-splitting}
 S^+|_Z=E\oplus E',
\end{equation}
for which $E'=E\otimes\mathcal A_Z$.
Then the real structure $\sr$ on $S^+$ restricts to real structures on both $E$ and $E'$.
Moreover,
\begin{align}
	\langle w_1(E^{\sr}), [\gamma_k]\rangle
	&=\langle w_1(E'^{\sr}), [\gamma_k]\rangle,\\
	\langle c_1(E),[Z]\rangle \label{eq:OT-Z-1}
	&\equiv\sum_{k=1}^r\langle w_1(E^{\sr}),[\gamma_k]\rangle\bmod2,\\
	q_{\sF_{\tf}}([\gamma_k]) \label{eq:OT-Z-2}
	&=\langle w_1(E^{\sr}),[\gamma_k]\rangle+1.
\end{align}
\end{lemma}

\begin{proof}[Proof of Proposition~\ref{prop:strong-comparison}, assuming Lemmas \ref{lem:multicurve-dual} and \ref{lem:real-eigenline}]
The second equality in \eqref{eq:strong-comparison} is \eqref{eq:P-topological} modulo two, so it suffices to prove the first.

If $\theta=0$, all three terms in \eqref{eq:strong-comparison} vanish.

If $\theta \neq 0$, choose a dual surface $Z$ of $\theta$ from Lemma~\ref{lem:multicurve-dual}.
Taking the determinant of
\eqref{eq:eigenline-splitting} and $\mathcal A_Z=(TZ)_{\C}\otimes(N_Z)_{\C}$ yields
\[
 c_1(\widehat\tf)\cdot[Z]
 =2\langle c_1(E),[Z]\rangle+\langle c_1(\mathcal A_Z),[Z]\rangle
 =2\langle c_1(E),[Z]\rangle+\chi(Z)+[Z]^2.
\]
Insert the above into \eqref{eq:P-topological} and apply the two preceding lemmas:
\[
Q_{\tf}(\theta) \equiv P_{\tf}(\theta)
 =\langle c_1(E),[Z]\rangle +\frac{\chi(Z)}2+[Z]^2
 \equiv\sum_{k=1}^r
       \bigl(\langle w_1(E^{\sr}),[\gamma_k]\rangle+1\bigr)
=\sum_{k=1}^r q_{\sF_{\tf}}([\gamma_k])
 =q_{\sF_{\tf}}([\gamma]).
\]
The first equality is the definition of $Q_{\mathfrak{t}}$,
and the last equality uses disjointness of $\gamma_k$'s.
\end{proof}

\medskip

\begin{proof}[Proof of Lemma~\ref{lem:multicurve-dual}]
To find the multicurve, choose a regular level of a smooth map $g \colon F \to S^1$ representing $\theta$, at a value $\ne 1$.
Take parallel components to represent multiplicities, if necessary.
Let $\mathbf{t}_k$ be the unit tangent of $\gamma_k$.
Let $a_O$ be a $\tau$-invariant arc from side $\mathsf{t}$ to side $\mathsf{t}'$ as in Figure~\ref{fig:octagon}.
Let $\mathbf{a}$ be the unit tangent of $a_O$ at $o$,
and let $\mathbf{n}$ be a unit vector so that $(\mathbf{n},\mathbf{a})$ is a positively oriented orthonormal frame in $TO$.
In the product region from Lemma~\ref{lem:octagon-model},
$N_F |_{F_0}\cong F_0 \times T_o O$.
The framing $\uptau$ lifts the constant downstairs framing $(n_F,\partial_t)$,
so it is constant in this product trivialization.
Since $(\mathbf{n},\mathbf{a})$ is also constant,
the two framings differ by a constant rotation.

The intersection of the cylinders $a_O\times\gamma_k$'s with an outgoing copy of $Y'$
is an oriented link with $r$ components.
By the assumption  $H_1(M;\Z)=0$,
this link bounds a properly embedded oriented surface in $M$.
One can assume this surface is connected by tubing the components and resolving double points.
Attach also the image of this surface under $\tau$, with reversed orientation; this lives in the other copy of the cap  $M$. The resulting surface is $Z$,
and $Z/\tau$ is a capping surface in $M$ of genus $h$ with $r$ boundary circles.
Hence \[\chi(Z)/2=2-2h-r\equiv r\bmod2.\]

Let $L_Z$ be the complex line bundle dual to the oriented surface $Z$,
with a section $s_Z$ vanishing transversely on $Z$.
Since $Z\cap Y=\emptyset$,
this section trivializes $L_Z$ near $Y$.
As $Z\cap W_F=a_O\times\gamma$,
restricting to $W_F$ yields
\[
 c_1(L_Z|_{W_F},s_Z|_Y)
 =\operatorname{PD}[a_O\times\gamma]
 =c_1(L_\theta,1).
\]
Thus $L_Z|_{W_F}\cong L_\theta$ relative to the boundary trivializations on $Y$.
Furthermore,
$\tau^*L_Z \cong \overline{L}_Z$, because
$\tau$ reverses the orientation of the normal bundle of $Z$.
On the two swapped outgoing caps,
$L_Z$ defines conjugate line bundles.
By Proposition~\ref{prop:real-extension-WF},
each $L_Z$ can be equipped with a unique real structure, compatible with $L_{\theta}$ under the isomorphism.
Let now $\widehat L_\theta=L_Z$ over $W_{F,M}$.
This construction is done separately for each $\theta$, by Lemma~\ref{lem:Q-cap-independent}.

For the self-intersection,
observe that $(NZ)^{\tau}|_{\gamma}=\nu(\gamma \subset F)$.
Since $F$ is oriented,
each $\gamma_k$ is two-sided.
Thus the normal line $(NZ)^{\tau}|_{\gamma}$ is trivial,
and its nonzero section $\mathbf{v}$ extends to an equivariant normal section near $\gamma$.
Extend $\mathbf{v}$ further over the free part of $Z$,
and perturb there equivariantly to make its zeros transverse.
The extended global normal section has no zeros near $\gamma$,
and the remaining zeros appear in $\tau$-pairs.
Hence $[Z]^2$ is even.
\end{proof}

\begin{proof}[Proof of Lemma~\ref{lem:real-eigenline}]
The involution $\tau$ reverses the orientation of $Z$,
and thus reverses the normal orientation.
Consequently
\[\sr_z(J_N)_z=-(J_N)_{\tau z}\sr_z\]
for every $z\in Z$.
Conjugate-linearity implies $\sr(E_z)=E_{\tau z}$ and
	$\sr(E'_z)=E'_{\tau z}$.
Therefore $\sr$ defines real structures on $E$ and $E'$.
Thus \eqref{eq:OT-Z-1} is the rank-one case of \cite[Corollary~5.3]{OT}.
On each $\gamma_k$,
	the operator $I=\rho(\mathbf{v}_k)\rho(\mathbf{t}_k)$ \emph{commutes} with $\sr$,
and satisfies $I^2=-1$.
Hence $I$ defines a complex structure on the real rank-two bundle $(S^+)^{\sr}$.
Since $(S^+)^{\sr}$ is the direct sum $E^{\sr} \oplus E'^{\sr}$,
the $w_1$ classes of the two summands agree.

Let us apply the identification of two spin covers along $\gamma_k$ in the proof of Lemma~\ref{lem:Gamma} along $\gamma_k$,
with $(\partial_t,e_1;\nu_F,e_3)=(\mathbf{v}_k,\mathbf{t}_k;\mathbf{n},\mathbf{a})$.
The constant rotation from the normal framing $\uptau$ to
$(\mathbf{n},\mathbf{a})$ preserves spin monodromy.
The map $\Phi \colon E_- \to (S^+)^{\sr}$ from Lemma~\ref{lem:half4-new}
identifies the spin lifts of $\mathbf{t}_k^{-1}$ with the unit $(-1)$-eigenline bundle of $i\rho(\mathbf{v}_k)\rho(\mathbf{n})$.
Since $\{\mathbf{v}_k,\mathbf{n}\}$ span $N_Z$,
this real eigenline bundle is either $E^{\sr}$ or
$(E')^{\sr}$ along $\gamma_k$,
which have the same $w_1$.
Thus the monodromy of the spin lifts is $\langle w_1(E^{\sr}),[\gamma_k]\rangle$,
and \eqref{eq:OT-Z-2} follows from the definition of $q$.\footnote{Similar relations in Lemma~\ref{lem:real-eigenline} appear for real theta characteristics in \cite[4.3-4.6]{OT}.}
\end{proof}

\subsection{Capping by planar surface}

Subsection~\ref{sec:index-geometric} makes no use of the spin marking on the boundary $C$ of $F$.
This subsection shows that the spin marking is precisely the data needed for gluing two spin structures along $C$.
Let $G$ be a connected planar surface with precisely $n$ boundary components.
Write
\begin{equation*}
	\widehat{F} = F \cup_C G.
\end{equation*}
Orient the planar cap $G$ so that $\Fh$ is closed and oriented.
On each component $C_k$, suppose $p_k$ is a basepoint.

Let $\mathfrak s_F$ be a spin structure on $F$,
with \emph{boundary types} $\lambda_k=q_{\mathfrak s_F}([C_k])$,
and write $\lambda = (\lambda_k)$.
(So $\lambda_k = 1$ if and only if $\mathfrak{s}_F$ is the Lie group spin structure along $C_k$.)
There exists a unique spin structure on $G$ having boundary type $\lambda$ on $G$.
Indeed, since $G$ is planar,
the restriction homomorphism $H^1(G;\F)\to H^1(\partial G;\F)$ is injective
and its image consists of the vectors whose entries sum to zero.
The boundary types of the spin structure induced from the plane are all bounding.
Thus $\sum_k\lambda_k=0$.

Suppose the spin markings on $F$ and $G$ are taken over the same tangent vectors $v_k=e_1(p_k)$ along $C_k$.
For spin markings $\varepsilon\in\mathscr P(\mathfrak s_F)$ and $\mathsf{m} \in\mathscr P(\mathfrak s_G)$,
matching their entries along $C$ yields a spin structure
\[
\widehat{\mathfrak{s}} = \mathfrak s_F\cup_{\mathsf{m}=\varepsilon}\mathfrak s_G\quad\text{on }\Fh.
\]
Indeed, on each $C_k$,
the two restrictions of spin covers of $S(T\Fh)|_{C_k}$ from $F$ and $G$ have the same
monodromy,
each giving rise to two isomorphic double covers over $C_k$.
The gluing of the double covers from $F$ and $G$ is specified precisely by $m_k=\varepsilon_k$ over the basepoint $p_k$.

\begin{lemma}\label{lem:gluing}
The spin structures over $\widehat{F}$ extending $\mathfrak s_F$ form a torsor over $\mathscr{A} = \mathbb{Z}_2^n/\mathbb{Z}_2$.
Changing the spin marking on either $\mathfrak{s}_F$ or $\mathfrak{s}_G$ by $a\in \mathscr{A}$ translates the glued spin structure by $\PD_{\Fh}(u_a)$, where $u_a=\sum_k a_k[C_k] \in H_1(\widehat{F}; \F)$.
\end{lemma}
\begin{proof}
Every extension $\widehat{\mathfrak{s}}$ restricts to the unique $\mathfrak s_G$,
so it is obtained by some gluing $\mathfrak{s}_F \cup \mathfrak{s}_G$.
Automorphisms of the spin covers over the connected surfaces $F$ and $G$ are global signs.
Thus two gluing identifications yield isomorphic extensions precisely when they differ by a global sign.
Moreover,
changing $m_k$ or $\varepsilon_k$ amounts to tensoring $\widehat{\mathfrak{s}}$
by a real line bundle $\ell_k$ obtained from the trivial bundles on $F$ and $G$,
glued by $-1$ along $C_k$, and $+1$ on all other circles.
The monodromy of $\ell_k$ is given by intersection with the class $[C_k]$,
so $w_1(\ell_k) = \PD_{\Fh}[C_k]$.
Any extension spin structure arises from twisting some boundary components this way;
hence the extensions form an $\mathscr{A}$-torsor.
\end{proof}

Using the planar closure of $F$ and an additional spin marking $\mathsf{m}$,
we define a marked version of the Arf invariant,
which applies whenever $|C| \ge 1$.
This is equivalent to a special case of Randal-Williams' generalized Arf invariant \cite[Lemma~2.11]{Randal-Williams}.

\begin{definition}
Given a marked spin structure $(\mathfrak{s}_F, \varepsilon)$ on $F$,
and a spin marking $\mathsf{m}\in\mathscr P(\mathfrak s_G)$,
the \emph{marked Arf invariant} of $(\mathfrak{s}_F, \varepsilon, \mathsf{m})$ is the Arf invariant of
the spin structure on $\widehat{F}$ obtained by combining $\mathfrak{s}_F$ and $\mathfrak{s}_G$,
via $\mathsf{m}=\varepsilon$:
\begin{equation}
 \Arf_{\mathsf{m}}(F,\mathfrak s_F,\varepsilon)
 :=\Arf\!\left(q_{\mathfrak s_F\cup_{m=\varepsilon}\mathfrak s_G}\right).
\end{equation}
\end{definition}

\begin{lemma}\label{lem:pointed-arf}
Let $\ell_\theta$ be a relative real line bundle on $F$ with a fixed trivialization $1$ on the boundary collar for which $w_1(\ell_\theta,1)=\theta\bmod2$.
Let
$(\mathfrak s_F,\varepsilon)\otimes\ell_\theta
=(\mathfrak s_F\otimes\ell_\theta,\varepsilon\otimes1)$.
Then for every $a \in \mathscr{A}$:
\begin{equation}\label{eq:pointed-arf-laws}
 \begin{aligned}
 \Arf_{\mathsf{m}}\bigl(F,(\mathfrak s_F,\varepsilon)\otimes\ell_\theta\bigr)
  -\Arf_{\mathsf{m}}(F,\mathfrak s_F,\varepsilon)
  &=q_{\mathfrak s_F}(\gamma_\theta),\\
 \Arf_{\mathsf{m}}(F,\mathfrak s_F,a\varepsilon)
  -\Arf_{\mathsf{m}}(F,\mathfrak s_F,\varepsilon)&=\lambda(a),\\
 \Arf_{a\mathsf{m}}(F,\mathfrak s_F,\varepsilon)
  -\Arf_{\mathsf{m}}(F,\mathfrak s_F,\varepsilon)&=\lambda(a).
 \end{aligned}
\end{equation}
\end{lemma}
\begin{proof}
Gluing the relative line bundle $\ell_\theta$ to the trivial line bundle on $G$
produces a line bundle $\widehat{\ell}_{\theta}$ over $\widehat{F}$ whose $w_1$ is given by $\PD_{\Fh}(\gamma_\theta)$,
where $\gamma_{\theta} = \PD_F(\theta)$.
Hence tensoring by $\ell_{\theta}$ turns $\widehat{\mathfrak{s}} = \mathfrak{s}_F \cup \mathfrak{s}_G$ over $\Fh$ into $\widehat{\mathfrak{s}}\otimes \widehat{\ell}_{\theta}$.
The quadratic form $q_{\widehat{\mathfrak{s}}}$ satisfies,
for the intersection form $\widehat{B}$ on $\widehat{F}$ and every $\gamma \in H_1(\widehat{F}; \F)$:
\begin{equation}
\label{eq:q-shift-by-gamma}
	q_{\widehat{\mathfrak{s}} \otimes \widehat{\ell}_{\theta}}(\gamma) =
	q_{\widehat{\mathfrak{s}}}(\gamma) + \langle w_1(\widehat{\ell}_{\theta}), \gamma \rangle_{\widehat{F}}
	= q_{\widehat{\mathfrak{s}}}(\gamma) + \widehat{B}(\gamma_{\theta}, \gamma).
\end{equation}
The first equality in~\eqref{eq:pointed-arf-laws} then follows from \eqref{eq:q-shift-by-gamma} and the following elementary Arf translation identity.
For any quadratic refinement $q$ of a nondegenerate skew-bilinear form $B$,
we have
\begin{equation}
\label{eq:aarf}
\Arf(q+B(x,\cdot \, ))-\Arf(q)=q(x).
\end{equation}
This identity follows by computing in a symplectic basis; \cite[Theorem~3A]{Johnson} supplies the affine transformation law \eqref{eq:q-shift-by-gamma}.
By Lemma~\ref{lem:gluing},
a marking change by $a \in \mathscr{A}$, on either $F$ or $G$,  translates $\widehat{\mathfrak{s}}$ by $\PD_{\Fh}(u_a)$,
where $u_a = \sum_k a_k[C_k]$.
Applying \eqref{eq:q-shift-by-gamma} and \eqref{eq:aarf} again shows that the right hand sides of the second and last equalities of \eqref{eq:pointed-arf-laws} are
\begin{equation*}
	q_{\widehat{\mathfrak{s}}}(u_a)= \sum_{k=1}^n a_k \lambda_k = \lambda(a). \qedhere
\end{equation*}
\end{proof}

\begin{remark}
One can show that if $\lambda=0$,
then $\Arf_{\mathsf{m}}$ is independent of $\mathsf{m}$,
and its value coincides with the $\Arf(\bar{q}_{\mathfrak s_F})$, where $\bar{q}_{\mathfrak s_F}$ is the induced nondegenerate form on $H^1(F, \partial F; \F)/\rad(B)$.
This holds for example if $C$ is connected.
Recall from Remark~\ref{rem:radical-and-Q} that an element $\theta$ in the radical can make
\[Q_{\tf^\theta}(\cdot) - Q_{\tf}(\cdot) \equiv 0,
\qquad \text{while} \qquad
Q_{\mathfrak{t}}(\theta)=  \lambda(\theta)=1 \bmod 2.
\]
This means that the missing index variation to keep track of arises precisely from the boundary classes,
through $\lambda = Q|_{\rad(B)}$.
In view of Lemma~\ref{lem:pointed-arf},
the correction can be supplied by the marked Arf invariant.
\end{remark}

\subsection{Spin marking from the ordering of $C$}
In the remainder of this section,
we fix, once and for all, a rule to assign the spin marking $\mathsf{m}$ needed for the definition of the marked Arf invariant.
To this end,
	we fix an ordering $\mathsf{ord}_C$ on the components $(C_i)$ of $$C=C_1 \sqcup \cdots \sqcup C_n.$$

This accounts for the order-dependence of the mod-two grading in the main Theorem~\ref{thm:main}.
By Lemma~\ref{lem:gluing},
the planar cap $G$ has a unique spin structure $\mathfrak s_{\lambda}$ for every boundary type $\lambda$.
The spin markings lie above the tangent vectors $v_k=e_1(p_k)$,
where $e_1$ is determined by the given orientation $o_C$ of $C$.
Let us use the orientation of $G$ induced from $\widehat{F}$, so that $\partial G = -C$.
The following definition is Randal-Williams' $\theta^2$-monodromy for arcs~\cite[Definition~2.3]{Randal-Williams}, expressed in our terms.
See also \cite[Proposition~2.8]{Randal-Williams}.

\begin{definition}
Let $\alpha$ be any embedded arc in $G$ from $p_k$ to $p_j$,
normal to the boundary.
The \emph{transport}
\[ \mathcal{T}_\alpha:P_{\mathfrak s_G,k}\to P_{\mathfrak s_G,j}, \]
is defined by lifting the path of unit tangent vectors in $G$ that begins with the $-\pi/2$-turn from $e_1(p_k)$ to the inward normal,
then stays tangent along $\alpha$,
and ends with the $-\pi/2$-turn from the outward normal to $e_1(p_j)$.
\end{definition}

\begin{figure}[htb]
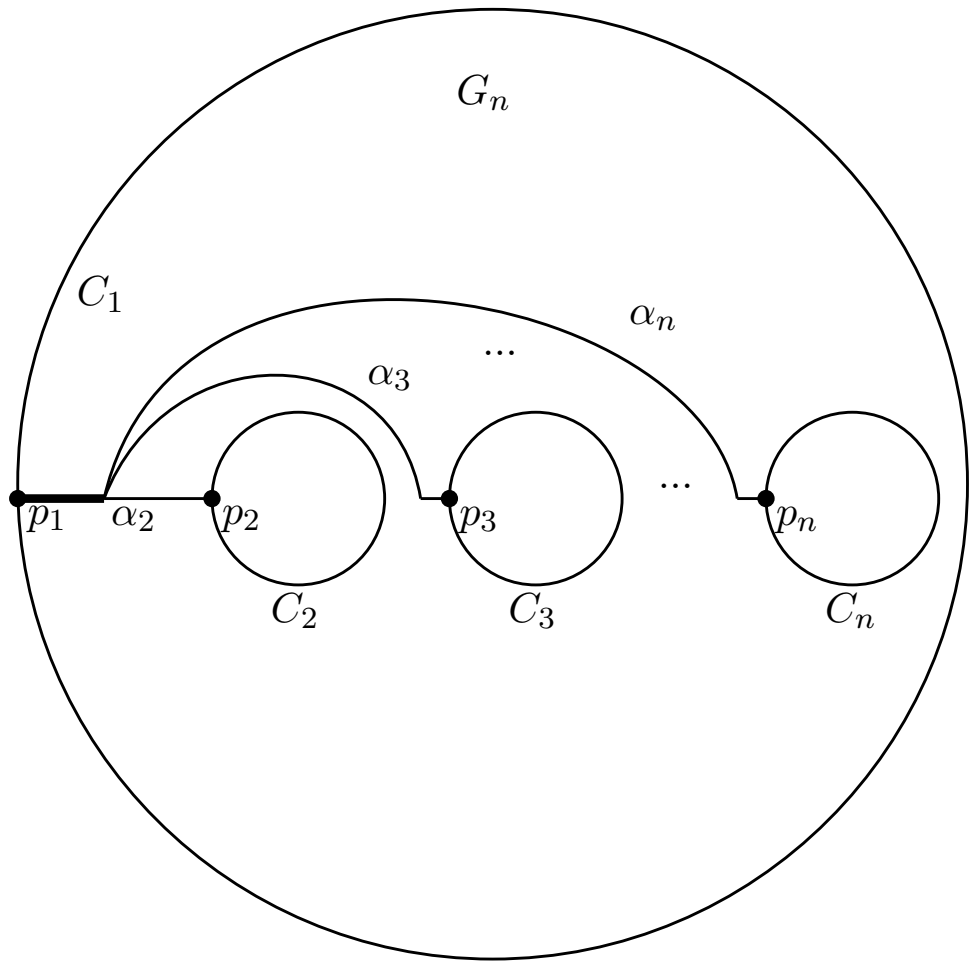

\centering
\adjustbox{max width=0.4\textwidth}{\Gndiagram}
\caption{The round disc model $G_n$.}
\label{fig:Gplanar}
\end{figure}

Let $G_n$ be a round disc in $\mathbb{R}^2$ with $n-1$ small round holes centred at $x_2<\cdots<x_n$ on the real axis.
Label the outer boundary $C_1$ and the inner boundaries $C_2,\ldots,C_n$,
and mark the leftmost point of each circle.
Equip $G_n$ the orientation opposite to the standard orientation on $\mathbb{R}^2$, and orient its boundary circles as $-\partial G_n$.
For $k\ge2$,
let $\alpha_k$ be the arc from $p_1$ to $p_k$ along the real line,
which goes above the intermediate holes by upper semicircles. See Figure~\ref{fig:Gplanar}.

\begin{definition}\label{def:ordered-marking}
Let $\mathfrak s_\lambda$ be the spin structure of boundary type $\lambda$ on $G_n$.
Choose $m_1\in P_{\mathfrak s_\lambda,1}$ over $p_1$ and define
\[
 \mathsf{m}_\lambda=[m_1,\Tcal_{\alpha_2}m_1,\ldots,\Tcal_{\alpha_n}m_1]
 \in\mathscr P(\mathfrak s_\lambda).
\]
Given any model for $G$,
choose an orientation-preserving diffeomorphism $\phi \colon (G, \partial G) \to (G_n, \partial G_n)$
that preserves the ordering of $\partial G$ and the marked points.
The \emph{cap marking} on $G$ is $\mathsf{m}=\phi^*\mathsf{m}_\lambda$.
\end{definition}

Given an ordering $\mathsf{ord}_C$ and the marking $\mathsf{m}$ as above,
we insert $\mathsf{m}$ into the marked Arf invariant.
This is well-defined by the following proposition.

\begin{proposition}

For a fixed ordering $\mathsf{ord}_C$ and a fixed marked spin structure $(\mathfrak s_F,\varepsilon)$ on $F$,
the value $\Arf_{\mathsf{m}}(F,\mathfrak s_F,\varepsilon)$ is independent
of the model identification $\phi:G\to G_n$.
\end{proposition}
\begin{proof}
Let $\phi,\phi' \colon G\to G_n$ be two model identifications,
that induce cap markings $\mathsf m,\mathsf m'$, respectively.
Then $g=\phi^{-1}\circ\phi'$ is an orientation-preserving diffeomorphism of $G$ preserving every marked point.
In particular, $\mathsf m'=g^*\mathsf m$.
Up to isotopy fixing the marked points,
assume $g$ is the identity on a collar of $\partial G$.
Then $g$ can be extended to a diffeomorphism $\widehat g$ of $\widehat F$ via identity across $F$.
Since $g$ preserves the boundary types,
$g^*\mathfrak s_G\cong\mathfrak s_G$ by uniqueness of the cap spin structure.
Pullback by $\widehat g$ preserves $\varepsilon_k$ on the $F$-side,
while carrying this identification to the one matching $g^* m_k= m'_k$ with $\varepsilon_k$.
Thus
\[
 \widehat g^{\,*}
 \bigl(\mathfrak s_F
       \cup_{\varepsilon=\mathsf m}\mathfrak s_G\bigr)
 \cong
 \mathfrak s_F
       \cup_{\varepsilon=\mathsf m'}\mathfrak s_G.
\]
The proposition then follows from the spin diffeomorphism invariance of the Arf invariant.
\end{proof}

\section{Absolute gradings from marked Arf invariants}
\label{sec:arf}
Let us recall the set up.
Let $(Y,\tau)$ be a 3-manifold with involution,
that arises as the double branched cover of an oriented link $C \subset Y' = Y/\tau$,
specified by a Seifert surface $F \subset Y'$,
with compatible orientation.
The link $C$ is viewed as the fixed-point set of $\tau$,
and we fix an ordering $\mathsf{ord}_{C}$ of its components $C_1 \sqcup \cdots \sqcup C_n$.
Each $C_k$ is marked by precisely one basepoint $p_k$.
Finally, let $\widetilde{\sr} = (\mathfrak{s}, T,\sr)$ be a marked triple.

Let $W_F\colon Y \to Y' \sqcup Y'$ be the Seifert cobordism,
and fix an admissible outgoing capping 4-manifold $M$,
so that $W_{F,M} = W_F \cup M \sqcup M$.
Suppose $\mathfrak{t}$ is the extending real $\Spinc$ structure on $W_F$ rel $Y$.

By \S\ref{sec:spin}, let $\mathfrak{s}_{\mathfrak{t}}$ be the spin structure on $F$ determined by the real part of $\mathfrak{t}$ along $F$.
By \S\ref{ss:spin-marking-on-boundary} and \eqref{eq:rrsc-pointing},
the marked triple $\widetilde{\sr}$ determines a spin marking $\varepsilon(\widetilde{\sr})$ on $\mathfrak{s}_{\mathfrak t}$ above the fibres,
and hence a marked spin structure on $(F,\{p_k\})$.
Let $\mathsf{m}$ be the cap marking on $G$ in the sense of Definition~\ref{def:ordered-marking},
which depends on the component ordering $\mathsf{ord}_C$.
Putting everything together and capping $F$ by $G$,
set
\begin{equation}
 \widehat{\mathfrak s}_{\mathfrak t}
 =\mathfrak s_{\mathfrak t}
   \cup_{\varepsilon(\widetilde{\sr})=m}\mathfrak s_G.
\end{equation}
Let
\[
	\Arf(\widehat{\mathfrak s}_{\mathfrak t}) :=
	\Arf_{\mathsf{m}}(F,\mathfrak s_{\mathfrak t},\varepsilon(\widetilde{\sr})).
\]

\subsection{The mod two grading}

\begin{definition}\label{defn:mod-two-final}
Given a framed real configuration $\mathfrak a$ and a perturbation $\mathfrak{q}$, the absolute mod-two grading of $(\mathfrak{a},\mathfrak{q})$ is
\begin{equation}\label{eq:absolute-grading}
	\begin{aligned}
		\gr^{(2)}(\mathfrak a)
		&= \mathcal I_{\widehat{\mathfrak t}}(\mathfrak a)
		+\Arf(\widehat{\mathfrak s}_{\mathfrak t}) \\
		&= \gr_{SW}(\mathfrak a;W_{F,M},\widehat{\tf})
		+\iota_R(W_{F,M})
		+\Arf(\widehat{\mathfrak s}_{\mathfrak t})
		\bmod2.
	\end{aligned}
\end{equation}
Here we suppress the perturbation $\mathfrak{q}$ from the notation,
although $\mathfrak{q}$ enters the APS boundary condition of the formal dimension of the real Seiberg--Witten moduli space $\gr_{SW}(\mathfrak a;W_{F,M},\widehat{\tf})$, and the normalized index $\mathcal{I}_{\widehat{\mathfrak{t}}}(\mathfrak{a})$.
\end{definition}
For $(\mathfrak a,\mathfrak q)$ as above,
choose a reference reducible configuration $\mathfrak{a}_0 = (B_0,0)$ lying in the same framed configuration space as  $\mathfrak a$,
for which $D_{B_0}$ is invertible.
Choose a path $\upgamma_0$ of configurations and a path $\mathfrak p_0$ of perturbations,
from $(\mathfrak a,\mathfrak q)$ to $(\mathfrak{a}_0,0)$.
Let $P_{\upgamma_0,\mathfrak p_0}$ be the Fredholm operator on the finite cylinder as in~\cite[\S5.6]{Li26},
with the corresponding spectral boundary projections.
By the additivity of APS index,
\[
\mathcal I_{\widehat{\mathfrak t}}(\mathfrak a)
=
\mathcal I_{\widehat{\mathfrak t}}(\mathfrak a_0)
+\ind_{\R}P_{\upgamma_0,\mathfrak p_0}
\bmod2.
\]
Consequently,
\begin{equation}\label{eq:relative-grading}
	\begin{aligned}
		\gr^{(2)}(\mathfrak a)
		&= \gr^{(2)}(\mathfrak a_0)
		+\ind_{\R}P_{\upgamma_0,\mathfrak p_0} \\
		&= \ind_{\C}^{APS}D^+_{\widehat{\tf}}
		+\Arf(\widehat{\mathfrak s}_{\mathfrak t})
		+\ind_{\R}P_{\upgamma_0,\mathfrak p_0}
		\bmod2,
	\end{aligned}
\end{equation}
where the APS boundary condition on $D^+_{\widehat{\tf}}$ is determined by $B_0$.
Moreover, \eqref{eq:relative-grading} holds for any reducible reference configuration having invertible Dirac operator.

\begin{proposition}

\label{prop:long-absolute-grading}
~
\begin{enumerate}[label=(\roman*)]
	\item Assume $(\mathfrak{a},\mathfrak{b})$ is a pair of framed configurations
		for which the real monopole Floer relative grading $\gr(\mathfrak{a},\mathfrak{b})$ is defined.
		Then
		\[\gr^{(2)}(\mathfrak{a})-\gr^{(2)}(\mathfrak{b}) = \gr(\mathfrak{a},\mathfrak{b}) \bmod 2.\]
	\item All marked triples representing the same relative real \text{spin\textsuperscript{c}} structure $\widetilde{\sr}$ define the same $\gr^{(2)}$ on $\tHMR_*(Y, \tau, \mathbf{p} ; \, \widetilde{\sr})$.
\end{enumerate}
\end{proposition}
\medskip

\begin{theorem}
~
\label{thm:long-absolute-grading}
\begin{enumerate}[label=(\roman*)]
	\item $\gr^{(2)}(\mathfrak a)$ is independent of the admissible outgoing cap $M$ with exchanged conjugate $\Spinc$ extensions,
		and independent of the real $\Spinc$ extension $\mathfrak{t}$ on $W_F$.
\item The grading $\gr^{(2)}$ may depend on $F$, but only through its class $[F]\in H_2(Y',C;\Z)$.
\item Fix an underlying \emph{absolute} real $\Spinc$ structure $(\mathfrak s,\sr)$. The $\Z/2$-graded isomorphism class of
\[
 \bigoplus_{\widetilde{\sr}\in
              \RRSC(\mathfrak s,\sr)}
 \tHMR_*^{\mathsf{ord}_C}
       (Y, \widetilde{\sr}\,)
\]
is independent of the ordering $\mathsf{ord}_C$ of the components.
Here, $\RRSC(\mathfrak s,\sr)$ denotes the set of relative real \text{spin\textsuperscript{c}} structures whose underlying absolute real $\Spinc$  structure is $(\s,\sr)$.
Thus the \emph{total} framed real monopole Floer homology
\[
 \bigoplus_{\widetilde{\sr}\in
 \underline{\RSpin}^c(Y,\tau,\mathbf{p})}
 \tHMR_*^{\mathsf{ord}_C}
       (Y, \widetilde{\sr} \,)
\]
is independent of $\ord_C$ up to graded isomorphism of $\F$-vector spaces.
Here, $\underline{\RSpin}^c(Y,\tau,\mathbf{p})$ denotes the set of relative real \text{spin\textsuperscript{c}} structures on the marked real 3-manifold $(Y,\tau,\mathbf{p})$.
\end{enumerate}
\end{theorem}

\bigskip
\noindent
Theorem~\ref{thm:main} follows from the combination of the Proposition~\ref{prop:long-absolute-grading} and Theorem~\ref{thm:long-absolute-grading}.
The proofs of these two results occupy \S\ref{ssec:proof-of-long-prop} and \S\ref{ssec:proof-of-long-thm}, respectively.

\subsection{Proof of Proposition~\ref{prop:long-absolute-grading}}
\label{ssec:proof-of-long-prop}

\subsubsection*{(i) From an absolute to a relative grading}
Fix two pairs of configurations and tame perturbations, $(\mathfrak{a},\mathfrak{q})$ and $(\mathfrak{b},\mathfrak{r})$, in the real framed configuration space.
Let $(\upgamma,\mathfrak{p})$ be a path of configurations and perturbations from $(\mathfrak{a},\mathfrak{q})$ to $(\mathfrak{b},\mathfrak{r})$.
The relative grading in real framed monopole Floer homology is defined as
\begin{equation*}
	\gr_{(\upgamma,\mathfrak{p})}(\mathfrak{a},\mathfrak{b}) = \ind(P_{\upgamma,\mathfrak{p}}),
\end{equation*}
where we suppress the perturbation $\mathfrak{q}$ from notation.
By the homotopy invariance of the spectral flow, this integer depends only on the homotopy class of paths of $\upgamma$ with fixed endpoints.
Moreover, by~\cite[Proposition~9.3]{Li26}, the assumptions that every component of $C$ contains a basepoint,
and that the cover $Y \to Y'$ is defined using an oriented Seifert surface, imply that the relative grading mod-two is independent of $(\upgamma, \mathfrak{p})$.

From \eqref{eq:absolute-grading} we obtain
\[\gr^{(2)}(\mathfrak{a})-\gr^{(2)}(\mathfrak{b})
\equiv \gr_{SW}(\mathfrak a;W_{F,M},\widehat{\tf}) - \gr_{SW}(\mathfrak b;W_{F,M},\widehat{\tf}) = \ind_{\R}P_{\upgamma,\mathfrak p} = \gr(\mathfrak{a},\mathfrak{b})
\bmod2. \tag*{\qedsymbol}\]

\subsubsection*{(ii) Marked triple representatives of RRSC}
Let $\widetilde{\sr}=(\mathfrak{s},T,\sr)$ and $ \widetilde{\sr}'=(\mathfrak{s}',T',\sr')$ be two marked triples,
and $u \colon \widetilde{\sr} \to \widetilde{\sr}'$ be an isomorphism of marked triples.
Then $u$ is an unframed real gauge transformation,
which acts on the framed real configuration space of $[\widetilde{\sr}]$ via $\mathfrak{a} \mapsto u \mathfrak{a}$.
Let $\mathfrak{q} \mapsto u_* \mathfrak{q}$ be the corresponding action for perturbations.
Suppose $(\mathfrak{t},\phi)$ is an extension of the $\Spinc$ structure $\mathfrak{s}$ over $W_F$ rel $Y$,
with boundary identification $\phi\colon \mathfrak{t}|_Y \to \mathfrak{s}$.
For $\widetilde{\sr}'$, choose the same $\mathfrak{t}$, but use $\phi' = u \circ \phi$.
Then
\begin{equation*}
	(\phi')^* \mathfrak{a}' = \phi^* \mathfrak{a}, \qquad
	(\phi')^*\mathfrak{q}'=\phi^*\mathfrak{q}.
\end{equation*}
Thus the boundary operators and spectral projections are identical.
Moreover, with the same interior operator,
the APS indices in \eqref{eq:absolute-grading} are identical for $\widetilde{\sr}$ and $\widetilde{\sr}'$.

The Arf term remains the same also.
Indeed, $T'_k \circ u_{p_k} = T_k$.
Thus the fibres at $p_k$ give the same spin lifts under pullback:
\begin{equation*}
	(\phi_{p_k}')^{-1}(T_k')^{-1}(1,0) = \phi_{p_k}^{-1}T_k^{-1}(1,0).
\end{equation*}
The boundary spin markings and spin structures agree.
Denote by $\gr^{(2)}_{\widetilde{\sr},\mathfrak{q}}(\mathfrak{a})$ the grading of $(\mathfrak{a},\mathfrak{q})$ defined with $\widetilde{\sr}$.
Then
\begin{equation*}
	\gr^{(2)}_{\widetilde{\sr},\mathfrak{q}}(\mathfrak{a}) = \gr^{(2)}_{u\widetilde{\sr},u_*\mathfrak{q}}(u\mathfrak{a}).\tag*{\qedsymbol}
\end{equation*}

\subsection{Proof of Theorem~\ref{thm:long-absolute-grading}}
\label{ssec:proof-of-long-thm}

\subsubsection*{(i) Cap $M$ and extension $\mathfrak{t}$}
It suffices to show invariance for a reducible configuration $\mathfrak{a}$ with invertible Dirac operator.
The Arf term depends only on $\mathfrak{t}$ and the data on $Y$.
Hence Lemma~\ref{lem:APS-cap-change} implies
the independence of \eqref{eq:absolute-grading} of the outgoing cap $M$ and $\Spinc$ structure on $M$.
Choose an extending real $\Spinc$ structure $\mathfrak{t}$ rel $Y$.
Any other extension can be written as $\mathfrak{t}^{\theta}$ for some $\theta \in \Lambda$.
Therefore,
\begin{align*}
	&\gr_{SW}(\mathfrak a;W_{F,M},\widehat{\tf^{\theta}})+\iota_R(W_{F,M})+\Arf(\,\widehat{\mathfrak s}_{\mathfrak t^{\theta}})
	-\gr_{SW}(\mathfrak a;W_{F,M},\widehat{\tf})-\iota_R(W_{F,M})-\Arf(\,\widehat{\mathfrak s}_{\mathfrak t})\\
	&=\ind_{\C}^{APS}D^+_{\widehat\tf^\theta} - \ind_{\C}^{APS}D^+_{\widehat\tf}
	+\Arf(\,\widehat{\mathfrak s}_{\mathfrak t^{\theta}})-\Arf(\,\widehat{\mathfrak s}_{\mathfrak t})\\
	&= Q_{\mathfrak{t}}(\theta)+\Arf(\,\widehat{\mathfrak s}_{\mathfrak t^{\theta}})-\Arf(\,\widehat{\mathfrak s}_{\mathfrak t})\\
	&=Q_{\mathfrak{t}}(\theta) - q_{\widehat{\mathfrak{s}}_{\mathfrak{t}}}(\upgamma_{\theta})\\
	&= 0 \bmod 2.
\end{align*}
The first equality uses the assumption that $\mathfrak{a}$ is reducible,
and the second equality uses the definition of $Q_{\mathfrak{t}}$.
The third equality follows from Lemma~\ref{lem:spin-twist} and Lemma~\ref{lem:pointed-arf}.
The last equality is a consequence of Proposition~\ref{prop:strong-comparison}. \qed

\subsubsection*{(ii) Homology class of the Seifert surface}
It is well-known (\cite{GordonLitherland, Murasugi, Trotter}, see also~\cite[Lemma~4.3]{FriedlEtAl}) that all connected Seifert surface representatives of $[F] \in H_2(Y',C ; \mathbb{Z})$ are related by isotopies relative to $C$,
tube attachments and tube removals. In fact, any small isotopy can be viewed as a tube attachment followed by the removal of the same tube, slightly displaced. Thus, it suffices to prove that the grading is preserved under a tube attachment.

Let $\gamma$ be a properly embedded arc with interior in $Y'\setminus F$.
Let $D_1,D_2$ be the disc neighbourhoods of the endpoints of $\gamma$ in $F$.
A {\em tube attachment} amounts to
replacing $D_1 \sqcup D_2$ by the boundary $T$ of a tubular neighbourhood around $\gamma$, with $\partial T = \partial (D_1 \sqcup D_2)$.
Denote
\[
 F'=\bigl(F\setminus\operatorname{int}(D_1\cup D_2)\bigr)\cup T, \qquad
 F^\circ=F\setminus\operatorname{int}(D_1\cup D_2)
        =F'\setminus\operatorname{int}T.
\]

Upstairs, write $W=W_{F,M}$, $W'=W_{F',M}$ for the Seifert cobordisms, using the same outgoing cap $M$.
Then $W'$ is obtained from $W$ by an equivariant surgery (\cite{HildenMontesinos}).
Concretely,
we replace a $\tau$-invariant neighbourhood $N$ of $\pi_W^{-1}(\{\tfrac12\}\times\gamma)$ with $N'\cong D^2\times S^2$.
Denote the corresponding surgery trace as
\[
\mathcal{W}=(W\times[0,1])\cup_{S^1\times D^3}(D^2\times D^3).
\]
Any $\Spinc$ structure $\widehat\tf$ on $W$ extends over $\mathcal{W}$.
Indeed,
pull \(\widehat{\tf}\) back to \(W\times[0,1]\).
There is a unique $\Spinc$ structure on the attaching region (because $H^2(S^1\times D^3;\Z)=0$), and a unique $\Spinc$ structure on $D^2 \times D^3$.
Let $\widehat{\tf}_{\mathcal{W}}$ be such an extension, and let \(\widehat{\tf}'=\widehat{\tf}_{\mathcal W}|_{W'}\).

Let $\mathfrak a$ be a reducible configuration with invertible Dirac operator.
Choose the same metrics and extend the $\Spinc$ connections on $W,W'$ near $Y$,
so that the boundary conditions are identical.
Thus
\begin{equation*}
\mathcal{I}_{W,\widehat{\tf}}(\mathfrak{a})-
\mathcal{I}_{W',\widehat{\tf}'}(\mathfrak{a})
=\ind_{\C}^{APS}D^+_{\widehat\tf}
-\ind_{\C}^{APS}D^+_{\widehat\tf'}
=\ind_{\C}D^+_{\tf_{W\cup_Y(-W')}} = 0.
\end{equation*}
The last term vanishes by the bordism invariance of the Dirac index, and the fact that $\partial \mathcal{W} = W\cup_Y(-W')$.

It remains to prove the equality of the Arf term:
$\Arf(\,\widehat{\mathfrak s}_{\mathfrak t})=\Arf(\,\widehat{\mathfrak s}_{\mathfrak t'})$.
Denote the exterior of the surgery region by
\[
U=W_F\setminus\operatorname{int} N
=W_{F'}\setminus\operatorname{int} N'.
\]
Then the real \text{spin\textsuperscript{c}} structures on $W$ and $W'$ can be identified on $U$.
Furthermore, the normal framings (\S\ref{sec:framing}) $\uptau$ on $F$ and $F'$ upstairs can be identified on $F^\circ$,
as the corresponding downstairs framings coincide in $U/\tau$.
Since the spin structure $\mathfrak{s}_{\tf}$ is induced by $\uptau$,
we have $\mathfrak s_{\tf}|_{F^\circ} \cong \mathfrak s_{\tf'}|_{F^\circ}$.
Finally, the identification of $\widehat{\mathfrak{t}}$ and $\widehat{\mathfrak{t}}'$ on $U$ preserves every spin marking on $F^{\circ}$.
For $F$ and $F'$,
use the same planar cap $(G,\mathfrak{s}_G)$ and ordered marking $\mathsf m$.
Matching the spin markings on $F^{\circ}$ and $G$,
we have an identification:
\[
\widehat{\mathfrak s}_{\tf} \cong \widehat{\mathfrak s}_{\tf'} \quad \text{on} \quad
\widehat F^\circ=F^\circ\cup_C G.
\]

Let $V=H_1(\widehat F;\F)$,
$q=q_{\widehat{\mathfrak s}_{\tf}}$ and $q'=q_{\widehat{\mathfrak s}_{\tf'}}$.
Choose a curve $c \subset F^\circ$ that is parallel to $\partial D_1$, which is isotopic in $F'$ to the meridian of the annulus $T$.
Furthermore, pick a symplectic basis of $V$ represented by curves in the complement of $D_1$ and $D_2$.
The same set of curves in $\widehat F'$ spans a symplectic subspace $V_0$ on which $q'$ agrees with $q$.
Let $e$ represent a curve passing through $T$ once and intersecting $c$ once,
whose class is orthogonal to $V_0$.
Since $c$ bounds a disc in $F$, we have $q(c)=q'(c)=0$.
Thus,
\[
\Arf(\widehat{\mathfrak s}_{\tf'})
=\Arf(\widehat{\mathfrak s}_{\tf})+q'(c)q'(e)
=\Arf(\widehat{\mathfrak s}_{\tf}).
\]

We deduce  that the grading is unchanged by a tube attachment.
\qed

\subsubsection*{(iii). Ordering and isomorphism class of Floer homology}
Let $\mathsf{ord}_C$ and $\mathsf{ord}'_C$ be two orderings.
Use the same planar cap $G$ for both constructions,
and denote the resulting cap markings by $\mathsf m$ and $\mathsf m'$.
Both markings lie above the same $\mathfrak s_G$.
Hence there exists an $a=(a_k)\in \mathbb{Z}_2^n$,
for which $\mathsf m'=a\mathsf m$.

The marked triples $\widetilde{\sr} = (\mathfrak s,T,\sr)$ and $a\widetilde{\sr} = (\mathfrak s,a T,\sr)$ have the same real configuration space and framed gauge group.
Using the same metric and perturbation for both triples,
the critical points and trajectory spaces of $\widetilde{\sr}$ and $a\widetilde{\sr}$ coincide.
The identity map on the common configuration space then induces an isomorphism of ungraded Floer chain complexes associated to $\widetilde{\sr}$ and $a\widetilde{\sr}$.
Let $\varepsilon = (\varepsilon_k) = \varepsilon(\widetilde{\sr})$ be the induced spin marking on $F$.

Changing the spin markings once through $\mathsf{m}' = a\mathsf{m}$ and once through $a T$ preserves $\widehat{\mathfrak{s}}_{\mathfrak{t}}$ and hence the Arf invariant:
\[
 \Arf_{\mathsf m'}
 \bigl(F,\mathfrak s_{\mathfrak t},
           \varepsilon(a\widetilde{\sr})\bigr)
 =
 \Arf_{\mathsf m}
 \bigl(F,\mathfrak s_{\mathfrak t},
           \varepsilon(\widetilde{\sr})\bigr).
\]
The normalized index also remains unchanged, so
the chain complex isomorphism above induces a degree-zero isomorphism:
\[
 \tHMR_*^{\mathsf{ord}_C}
       (Y,\widetilde{\sr})
 \cong
 \tHMR_*^{\mathsf{ord}'_C}
       (Y,a\widetilde{\sr}).
\]
Taking direct sum and reindexing by the permutation $\widetilde{\sr}\mapsto a\widetilde{\sr}$ proves part (iii) of Theorem~\ref{thm:long-absolute-grading}.
\qed

\section{Examples}
\label{sec:examples}
\subsection{The unknot $\boldsymbol{U}$}
Let $U\subset S^3$ be the unknot with one basepoint $p$,
and $F=D^2$ be the standard disc bounding $U$.
The branched double cover is $Y\cong S^3$,
which admits a unique \RRSC\ structure $\widetilde{\sr}$.
Take $M$ to be a 4-ball.
Then $W=W_{F,M}$ is the double cover of a 4-ball branched along the standard 2-disc,
so $W\cong B^4$.
Consider the unique spin structure on $W$ with a compatible real lift.
The induced spin structure on $F$ is the unique spin structure
on $D^2$.
The planar cap $G$ is another disc, and hence
\[
\widehat F=F\cup_U G\cong S^2,\qquad
\Arf(\widehat{\mathfrak s}_{\tf})=0.
\]

Equip $Y$ with the round metric (which is $\tau$-invariant).
There is a unique reducible $\mathfrak a_0=(B_0,0)$ with flat
determinant connection.
By the Weitzenb\"ock formula,
the unperturbed equation admits no irreducible solutions and $\ker D_{B_0}=0$.
Thus $\mathfrak a_0$ is nondegenerate and generates the framed chain complex $\tcmr(S^3,\tau,p)$.
Choose a $\tau$-invariant metric of positive scalar curvature (PSC) on $W$, and
attach a cylindrical end to $W$. By the PSC property of the metric on $W$,
both operators $D^{\pm,\sr}_{\widehat{\tf}}$ have trivial kernel.
Thus
\[
\ind_{\R}^{APS}D^{+,\sr}_{\widehat{\tf}}
=\ind_{\C}^{APS}D^+_{\widehat{\tf}}=0.
\]
Therefore
\[
\gr^{(2)}(\mathfrak a_0)
=\ind_{\C}^{APS}D^+_{\widehat{\tf}}
+\Arf(\widehat{\mathfrak s}_{\tf})
=0\bmod2.
\]
This implies
\[
\tHMR(S^3,\tau,p;\widetilde{\sr})\cong\mathbb{F}_{(0)},\quad
 \chi\bigl(\tHMR(S^3,\tau,p;\widetilde{\sr})\bigr)=1.
\]
\subsection{Hopf links and $\boldsymbol{S}^1 \times \boldsymbol{S}^2$}\label{sec:annulus}
In this subsection we discuss three examples of real manifolds: the two-fold branched covers of the two Hopf links, and  $S^1\times S^2$ with the product involution $\mathrm{id}\times \Rot_\pi$. We treat these examples together because they have similar properties; for example, they all admit  annular Seifert surfaces.

\subsubsection*{Spin structures on the capped $2$-torus}
Fix an ordering of components $(C_1,C_2)$,
which determines the cap marking $\mathsf m$ of $G$ by Definition~\ref{def:ordered-marking}.
Let the Seifert surface $F$ and the planar cap $G$ be annuli,
so \[\Fh=F\cup_CG\cong T^2.\]
Denote
\begin{equation}
\label{eqn:2-torus-notations}
 \Lambda=\Z\theta_1,\quad \theta_1=\delta_F(1_{C_1}),\quad
 V=\F\langle c\rangle,\quad B=0, \quad c=[C_1]=[C_2].
\end{equation}
Here $\delta_F \colon H^0(C;\Z)\to H^1(F,C;\Z)$ is the connecting map.
Then $V=\rad(B)$, and the marking group $\mathscr{A}=\{\pm1\}^2/\{\pm1\}\cong\F$.
Let $e$ be a circle in $\Fh$ that satisfies $\Bh(c,e)=1$,
and $\widehat{\mathfrak s}_{\tf}$ be the glued spin structure.
Write
\[
 \widehat q=q_{\widehat{\mathfrak s}_{\tf}},\qquad
 \lambda=\widehat q(c),\qquad v=\widehat q(e).
\]
Assume $\mathfrak{t}$ is a real $\Spinc$ structure on $W_{F}$,
inducing $\mathfrak{s}$ on $Y$ and a spin structure $\mathfrak{s}_{\mathfrak{t}}$ on $F$.
The two boundary types of $\mathfrak{s}_{\mathfrak{t}}$ are both $\lambda$, and
\begin{equation}\label{eq:arf-product}
 \Arf(\widehat{\mathfrak s}_{\tf})=\lambda v.
\end{equation}
The two spin structures on $\widehat F$ extending $\mathfrak s_{\tf}$ can be specified by
$(\widehat q(c),\widehat q(e))=(\lambda,v)$.
In particular, the Arf invariants of $(0,0)$ and $(0,1)$ vanish,
whereas the Arf invariants of $(1,0)$ and $(1,1)$ differ by one.
By Lemma~\ref{lem:gluing}, changing the marking by the nontrivial element of $\mathscr{A}$ translates the glued spin structure by $\PD_{\widehat F}(c)$.
Consequently, this marking change swaps the extensions $(\lambda,0)$ and $(\lambda,1)$.
For a reducible configuration $\mathfrak a$ with invertible Dirac operator,
inserting \eqref{eq:arf-product} into Definition~\ref{defn:mod-two-final} shows:
\begin{equation}\label{eq:annulus-grading}
 \gr^{(2)}(\mathfrak a)
 =\mathcal I_{\tf}(\mathfrak a)+\lambda v \ \bmod2.
\end{equation}
After changing $\tf$ to $\tf^{\theta_1}$,
the spin structure $\widehat{\mathfrak s}_{\mathfrak t}$ becomes $\widehat{\mathfrak s}_{\mathfrak t}+\PD_{\widehat F}(c)$,
by Lemma~\ref{lem:gluing}.
Hence the $\lambda v$-term in~\eqref{eq:annulus-grading} changes to $\lambda v+ \lambda$.

\subsubsection*{APS indices}
Denote by $\tf$ a real $\Spinc$ structure over $W_F$ rel $Y$ and $\widehat\tf$ a real extension of $\tf$ over $W=W_{F,M}$.
Let $Z=Z_{\theta_1}$ be the dual surface (\S\ref{sec:index-geometric}),
where $\theta_1 \in \Lambda$ as in~\eqref{eqn:2-torus-notations}.
Set $x_{\theta_1}=\PD[Z]$.
Write $h = [Z]^2/2$ and $k=\langle c_1(\widehat\tf),[Z]\rangle/2$.
By Lemma~\ref{lem:multicurve-dual}, $h$ and $k$ are integers.
Equation~\eqref{eq:P-topological} implies
\begin{equation}
	\label{eq:P_t_hkj}
	P_{\tf}(j\theta_1)=hj^2+kj.
\end{equation}
Then by Proposition~\ref{prop:strong-comparison}:
\begin{equation}
	\lambda=q_{\mathfrak s_{\tf}}(c)
	=Q_{\tf}(\theta_1)=h+k.
\end{equation}
Twisting $\widehat{\mathfrak{t}}$ by $j\theta_1$
translates $c_1(\widehat{\mathfrak{t}})$ by $2jx_{\theta_1}$,
and hence translates $k$ to $k+2hj$.
This fixes $\lambda$.
In what follows,
assume $\mathfrak a = (B,0)$ is a reducible configuration for which $D_B$ is invertible and $B^t$ is flat.
Then the normalized index is precisely the index of the Dirac operator on $W = W_{F,M}$ with APS condition:
\[
\mathcal I_{\tf}(\mathfrak a)
=
\gr_{SW}(\mathfrak a;W,\widehat{\tf})+\iota_R(W)
= \ind_{\C}^{APS}D^+_{\widehat\tf}.
\]

\subsubsection*{The Hopf links.}
Let $H_{\epsilon}\subset S^3$ be the Hopf link with linking number $\epsilon =\pm1$,
oriented as the boundary of an annulus $F$.
The double branched cover of $H_{\epsilon}$ is $\mathbb{RP}^3$,
with deck transformation $\tau[z_1,z_2]=[z_1,-z_2]$ in the model $S^3/\{\pm1\}$.
The fixed circles $C_1,C_2$ are $\{[z_1,0]\}$ and $\{[0,z_2]\}$, respectively.
Taking $M=-B^4$ gives a filling $W = W_{F,M}$ with intersection form $[\pm 2]$.
Let $\mathfrak s_{\mathrm{spin}}$ be the real $\Spinc$ structure on $Y$ induced by the spin structure on $W$,
and let $\mathfrak s_{\mathrm{tw}}$ denote the other real $\Spinc$ structure.
Lemma~\ref{prop:hopf-index} shows that $\lambda=1$ for $\mathfrak s_{\mathrm{spin}}$ and $\lambda=0$ for
$\mathfrak s_{\mathrm{tw}}$.
Denote the generator of $H^2(W;\mathbb{Z}) \cong \mathbb{Z}$ by $u$.
Then $\langle u,[Z]\rangle=1$, $h=\epsilon$,
and $[Z]^2= 2\epsilon$.
Let $\widehat\tf$ be the $\Spinc$ structure on $W$ with $c_1(\widehat\tf)=2ku$.
Then $\widehat\tf|_Y = \mathfrak s_{\mathrm{spin}}$ if $k$ is even,
and $\widehat\tf|_Y = \mathfrak s_{\mathrm{tw}}$ if $k$ is odd.
Furthermore, $\lambda=1+k\pmod2$.

\begin{lemma}\label{prop:hopf-index}
	For a PSC metric on $Y \cong \mathbb{RP}^3$,
	\begin{equation}
		\ind_{\C}^{APS}D^+_{\widehat\tf}=
		\begin{cases}
			\epsilon k^2/4, & k\ \text{even},\\
			\epsilon(k^2-1)/4, & k\ \text{odd}.
		\end{cases}
	\end{equation}
	Therefore $\mathcal I_{\tf}=0$ for every extension $\mathfrak{t}$ of $\mathfrak s_{\mathrm{tw}}$,
	and $\mathcal I_{\tf}=k/2\bmod2$ for every extension of $\mathfrak s_{\mathrm{spin}}$.
\end{lemma}

\begin{proof}
	Equip $W$ with a PSC metric that restricts to a PSC metric on a boundary collar of $Y$.
	For example, truncate the cylindrical-end metric in \cite[Fact~4.1(2)]{AB02}.
	For $k=0$, the extension $\widehat\tf$ arises from the spin structure, and its spin connection $A$ on $W$ has flat $U(1)$ part.
	The Weitzenb\"ock formula implies the APS index of $D^+_A$ is zero.
	For $k=2\epsilon j$, Equation~\eqref{eq:P_t_hkj} gives $P_{\tf}(j\theta_1)=\epsilon j^2=\epsilon k^2/4$.

	For $k$ odd, $c_1(\widehat\tf)$ is not torsion.
	To this end, let $W'$ be the disc bundle over $S^2$ of Euler number $-2\epsilon$,
	with the same metric as $W$.
	Let $X = W \cup_Y (-W')$.
	The restriction of the spin structure on $-W'$ onto $Y$ is $\mathfrak{s}_{\mathrm{tw}}$ (via a diffeomorphism $Y \to -Y$).
	Let $\mathfrak{t}_X$ be the $\Spinc$ structure on $X$ obtained by gluing $\widehat{\mathfrak{t}}$ to the unique torsion $\Spinc$ structure on $-W'$.
	By the preceding paragraph, the APS Dirac index on $-W'$ is zero.
	Since the boundary Dirac operator is invertible, APS gluing implies
	\begin{equation*}
		\ind_{\C}^{APS}D^+_{\widehat\tf} + 0 = \frac18 \left(c_1(\tf_X)^2-\sigma(X)\right)
	\end{equation*}
	Moreover, the zero sections $Z$ and $Z'$ of the bundles $W$ and $-W'$
	are orthogonal and both square to $2\epsilon$.
	Thus $\sigma(X)=2\epsilon$.
	Since $\langle c_1(\tf_X), [Z] \rangle = 2k$ and  $\langle c_1(\tf_X), [Z'] \rangle = 0$,
	we obtain $c_1(\tf_X)^2=2\epsilon k^2$.
	Therefore, the Dirac index is $\epsilon(k^2-1)/4$.
\end{proof}

\subsubsection*{A product involution on $S^1 \times S^2$.}
On $Y=S^1\times S^2$, let $\tau=\mathrm{id}\times \Rot_\pi$,
where $\Rot_\pi$ is the $\pi$-rotation about the poles $N,S$ of $S^2$.
The quotient $Y'$ is again $S^1\times S^2$.
Choosing an arc $\alpha$ from $N$ to $S$ in the quotient $S^2$ and a point $z$ in the interior of $\alpha$, we take
\[
 C=S^1\times\{N\}\sqcup S^1\times\{S\},\qquad
 F=S^1\times\alpha, \qquad \gamma = S^1 \times \{z\}.
\]

This example was studied in \cite[Section 14.4]{Li22}.
Only the torsion ordinary $\Spinc$ structure $\mathfrak s_0$ admits real structures.
There are two such real structures, distinguished by the first Stiefel--Whitney class of the real fixed line bundle.
Indeed, the cap $M=D^2\times S^2$ with $\partial M=-Y'$ gives a spin filling $W=W_{F,M}$, with a spin structure $\widehat{\mathfrak{t}}_0$.
In particular,
\begin{equation*}
	W_F = S^1 \times \dbcv_2([0,1]\times S^2,[0,\tfrac12]\times\{N,S\}\cup\{\tfrac12\}\times\alpha).
\end{equation*}
Let $Z_i$ be the $* \times S^2$ in the $i$-th capping $M$, for $i = 1,2$.
In addition, let $Z$ be the annulus $a_O \times \gamma$ capped by the two $D^2\times\{\ast\}$ in the outgoing caps. (Here, $a_O$ is the arc from Figure~\ref{fig:octagon}.)
Using the product normal framing, we have $[Z]^2 = 0$.
Also, $[Z_i]^2 = 0$.
Orient the $Z_i$'s so that $[Z]\cdot [Z_i]=1$.
Then $\tau_*[Z]=-[Z]$ and $\tau_*[Z_1]=-[Z_2]$.
Furthermore,
$H_2(W;\Z)$ has basis $([Z],[Z_1],[Z_2])$.
Thus $W$ is even, and satisfies $\sigma(W) = 0$, $b_R^1(W)=0$, and $b_R^+(W)=1$.
Indeed, the $\tau$-skew-invariant subspace of $H_2$ is generated by $[Z]$ and $[Z_1] + [Z_2]$,
and the intersection form is hyperbolic.

Denote the boundary real lift of the real spin structure on $W$ by $\sr_0$.
Then the other lift is $g\sr_0$, where $g(e^{i\varphi},x)=e^{i\varphi}$.
The boundary types of $(\mathfrak s_0,\sr_0)$ and  $(\mathfrak{s}_0,g\sr_0)$ are $\lambda=0$ and $\lambda=1$, respectively, by Lemma~\ref{prop:s1s2}.
By definition, $\widehat{\tf}_0$ is a real spin structure that extends $(\mathfrak{s}_0,\sr_0)$ and has $k=0$.
The other boundary real $\Spinc$ structure $(\mathfrak{s}_0, g\sr_0)$ extends as $\widehat\tf_0 \otimes E$, with $k=1$.
Here, $E$ is the line bundle satisfying $\langle c_1(E), Z \rangle = 1$ and $\langle c_1(E), Z_i\rangle = 0$.
Indeed, since $\tau^* E \cong \bar{E}$ and $b_1(W) = 0$,
there is a unique real structure $\sr_E$ on $E$.
Since $\tau|_Z$ fixes the core $\gamma$ of the annulus $F$, \cite[Corollary~5.3]{OT} shows
\[
\bigl\langle w_1(E^{\sr_E}|_F),[\gamma]\bigr\rangle
=\langle c_1(E),Z\rangle=1.
\]
In particular, $k=\langle c_1(\widehat{\mathfrak{t}}_0\otimes E), Z\rangle/2 = 1$.
Because the core is homologous to $C_1$ and $C_2$, tensoring $\widehat{\mathfrak{t}}_0$ by $E$ changes $\sr_0$ to $g\sr_0$ on the boundary $Y$.

\begin{lemma}\label{prop:s1s2}\label{prop:s1s2-index}
	We have $\ind_{\C}^{APS}D^+_{\widehat\tf_0}=0$.
	Hence $\mathcal I_{\tf}=0$ for every extension of
	$(\mathfrak s_0,\sr_0)$.
\end{lemma}

\begin{proof}
	Equip $W$ with a $\tau$-invariant metric of positive scalar curvature.
	For the spin connection of $\widehat\tf_0$, the Dirac APS index is zero because of the PSC property.
	Equation~\eqref{eq:P_t_hkj} with $h=k=0$, gives $P_{\tf}(j\theta_1)=0$.
	The extensions of $(\mathfrak s_0,\sr_0)$ are precisely $\widehat\tf_0^{\,j\theta_1}$.
\end{proof}

\subsubsection*{Gradings.}
The two basepoints on $C$ give rise to two distinct $\RRSC$ structures over each absolute real $\Spinc$ structure.
Hence, each example has four $\RRSC$ structures in total.
For each such $\RRSC$ structure,
there is a unique reducible configuration in the framed real monopole Floer chain complex.
Indeed,
use the round metric on $\mathbb{RP}^3$ and a product metric on $S^1\times S^2$,
so that both metrics are $\tau$-metrics and have positive scalar curvature (PSC).
There is a unique real $\Spinc$ connection $B$ with flat $U(1)$ connection,
and hence a unique reducible Seiberg--Witten solution up to gauge.
The Dirac operator $D_B$ is invertible, thanks to the Weitzenb\"ock formula and PSC.
Moreover, since $H^1(Y;\Z)^{-\tau^*}=0$, the de Rham part of the Hessian is also invertible.

By Lemmas~\ref{prop:hopf-index} and~\ref{prop:s1s2-index}, for the two Hopf links $H_{\pm}$ and $S^1\times S^2$ above,
the normalized index $\mathcal{I}_t(B,0)$ vanishes for the two $\lambda=0$ $\RRSC$ structures.
The Arf invariants, which are equal to $\lambda v$, also vanish.
Therefore, the mod-two gradings of corresponding generators of the framed real monopole Floer homology are both zero.
On the other hand,
by~\eqref{eq:annulus-grading}, the gradings for the two $\lambda=1$ $\RRSC$ structures have gradings $\mathcal I_{\tf}+v$ and $\mathcal I_{\tf}+v+1$.
In particular, the $\lambda = 1$ part does not contribute to the Euler characteristic.
We have proved:

\begin{proposition}
For the two Hopf links $H_{\pm}$ and for $S^1\times S^2$ with product involution,
the mod-two gradings of the generators in the $\lambda=0$ $\RRSC$ structures are both zero.
On the other hand, the gradings are $0$ and $1$ for the two $\lambda=1$ $\RRSC$ structures.
We have
\[
	\tHMR=\mathbb{F}_{(0)}\oplus \mathbb{F}_{(0)} \oplus \mathbb{F}_{(0)} \oplus \mathbb{F}_{(1)},
\qquad
\sum_{\widetilde{\sr} \in \underline{\RSpin}^{c}(Y,\tau,\mathbf{p})}
   \chi\bigl(\tHMR(Y,\tau,\mathbf p;\widetilde{\sr})\bigr)=+2.
\]
Here, the subscript denotes the mod-two grading of the $\mathbb{F}$-summand.
\end{proposition}

We summarize the calculations in the following table.
The pair in each ``gradings'' entry of the table corresponds to the two reducibles in the two distinct \RRSC\ structures with the same underlying absolute real $\Spinc$ structure.

\begin{center}
\small
\renewcommand{\arraystretch}{1.2}
\begin{tabular}{l|lcc}
\hline
Example & Underlying $\RSpin^c$ structure & $\lambda$  & Gradings\\
\hline
$H_{\pm}$ & $\mathfrak s_{\mathrm{spin}}$ & $1$
 & $(0,1)$\\
$H_{\pm}$ & $\mathfrak s_{\mathrm{tw}}$ & $0$
 & $(0,0)$\\
 \hline
$S^1\times S^2$ & $(\mathfrak s_0,\sr_0)$ & $0$
 & $(0,0)$\\
$S^1\times S^2$ & $(\mathfrak s_0,g\sr_0)$ & $1$
 & $(0,1)$\\
\hline
\end{tabular}
\end{center}

\bigskip

For comparison, Srivastava's calculation \cite[Example~3.5]{Srivastava} gives the following gradings in real
Heegaard Floer homology:
\begin{center}
  \begin{tabular}{c|cc}
	\hline
           & $\mathfrak s_1^\pm$ & $\mathfrak s_2^\pm$ \\ \hline
    $H^+$  & $(0,0)$             & $(0,1)$             \\
    $H^-$  & $(0,1)$             & $(1,1)$ \\\hline
  \end{tabular}
\end{center}
The labels $H^{\pm}$ refer to the respective Heegaard diagrams,
while columns refer to the underlying real $\Spinc$ structures.
This table has the same pattern of an equal pair and an opposite pair as our table above.
However, our gradings do not depend on whether the Hopf link is $H_+$ or $H_-$,
in contrast with Srivastava's Euler characteristics for $\hHFR$, which are $+2$ for $H^+$
and $-2$ for $H^-$~\cite[Remark~4.2]{Srivastava}.
Hence our groups for $H_{\pm}$ agree with $\hHFR$ for $H^+$,
and differ by an overall shift for $H^-$ from $\hHFR$. This is consistent with Conjecture~\ref{conj:HFR}.

\subsection{Classical Arf invariant of knots and grading for spin structures}
\label{ssec:classical-arf-K}
Let $K \subset S^3$ be a knot and let $(Y,\tau)$ be the double branched cover $\dbcv_2(K)$ equipped with its deck transformation $\tau$.
Put one basepoint on $K$.
Let $\vartheta = (B,0)$ be a reducible configuration, with reducible Dirac operator $D_B$.
Since $K$ is a knot, \eqref{eq:grading-knot-intro} reads:
\begin{equation*}
	\gr^{(2)}(\vartheta) = \ind_{\C}^{\mathrm{APS}} D^+_{\widehat{\mathfrak{t}}} + \Arf(Q_{\mathfrak{t}}).
\end{equation*}
For $M = B^4$, we have $W_{F,M} = -\dbcv_2(B^4, F)$, where $F$ is the pushed-in Seifert surface.
It suffices to work with the orientation-reversed copy $X = \dbcv_2(B^4, F)$,
as the boundary Dirac operator has no kernel, and grading is mod two.
The 4-manifold $X$ is spin, with intersection form given by the matrix $V+V^T$,
where $V$ is the Seifert matrix~\cite{KauffmanTaylor}.
The \emph{Arf invariant} $\Arf(K)$ of $K$ is the Arf invariant of the quadratic form $Q$ on $H_2(X;\F)$, computed using integral representatives by
\begin{equation*}
	Q(x) = \frac{1}{2} x^T (V+V^T)x = x^T V x \ \bmod 2.
\end{equation*}
By \eqref{eq:index-difference}, for all $x \in H^2(X,Y;\mathbb{Z})$,
\begin{equation*}
Q_{\widehat{\mathfrak{t}}}(x) = \frac{x \cdot x + \langle c_1(\widehat{\mathfrak{t}}), x \rangle}{2} \bmod 2.
\end{equation*}

Furthermore, $Y$ admits a unique real spin structure $\mathfrak{s}_0$; see~\cite{KMT:Ktheory}.
Let $\vartheta_{\mathrm{spin}} = (B_{\mathrm{spin}},0)$ be the unique reducible configuration arising from the canonical spin connection $B_{\mathrm{spin}}$.
In particular, $D_{B_{\mathrm{spin}}}$ is quaternionic-linear.
\begin{lemma}\label{lem:spin-red-grading}
The mod-two grading of the spin reducible configuration $\vartheta_{\mathrm{spin}}$ is given by the \emph{classical} Arf invariant:
\begin{equation*}
	\gr^{(2)}(\vartheta_{\mathrm{spin}}) = \Arf(K).
\end{equation*}
Moreover,
\begin{equation}
	\label{eq:spin-determinant}
\gr^{(2)}(\vartheta_{\mathrm{spin}})=\frac{(\det K)^2 - 1}{8} \bmod 2
=\begin{cases}0,&\det K \equiv1,7\pmod8,\\1,&\det K\equiv3,5\pmod8.\end{cases}
\end{equation}
\end{lemma}

\begin{proof}
	Let $\widehat{\mathfrak{t}}$ be the spin structure which extends $\mathfrak{s}_0$.
	Then the 4-dimensional Dirac operator $D^+_{\widehat{\mathfrak{t}}}$ and the APS boundary problem are quaternionic and hence the real part still has even index.
	Hence $Q_{\widehat{\mathfrak{t}}} = Q$.
	The Arf to determinant assertion in~\eqref{eq:spin-determinant} is classical~\cite[\S3.5]{Levine}.
\end{proof}

\subsection{Two-bridge knots}
Let $Y = L(p,q)$ be a lens space. This is the branched cover of $S^3$ along the 2-bridge knot $S(p,q)$,
with deck transformation $\tau$.
Again, pick one basepoint.
Assume $p > q > 0$, $p$ is odd, and $\gcd(p,q)=1$.
By \cite[\S14.2]{Li22} and with a $\tau$-invariant PSC metric, the real monopole Floer chain complex is generated by non-degenerate reducibles, one for each real \text{spin\textsuperscript{c}} structure.
Since every \text{spin\textsuperscript{c}} structure admits precisely one real structure,
there are $p$ \RRSC\ structures in total.
It suffices to compute the mod-two grading of the reducible generators.
This can be done by appealing to the existing calculations of the Heegaard Floer $d$-invariants~\cite{OzSz} as follows.
Label the \text{spin\textsuperscript{c}} structures on $L(p,q)$ as $\mathfrak{s}_j$, $0 \le j < p$ using the convention of \cite{OzSz}.
In particular, the spin structure corresponds to index $j_0$ for which $2 j_0 \equiv q - 1 \pmod p$.
Let $\vartheta_j$ be the corresponding unique reducible monopole Floer generator for $\mathfrak{s}_j$.

Denote $d(p,q,j) = d(L(p,q),\mathfrak{s}_j)$.
By \cite[Proposition~4.8]{OzSz} and $d(-Y,\mathfrak{s})=-d(Y,\mathfrak{s})$, the correction term $d$ can be computed recursively:
\begin{equation}\label{eq:dpqj}
d(p,q,j)
=
-\frac14 + \frac{(2j+1-p-q)^2}{4pq} - d(q, p \bmod q, j \bmod q).
\end{equation}

\begin{lemma}
	\label{lem:gr-Lpq}
	Choose $r,h_j$ for which $qr \equiv 1 \pmod p$ and $2h_j \equiv j - j_0 \pmod p$.
	The mod-two gradings for $L(p,q)$ are:
\begin{equation}\label{eq:Lpq-grading}
	\gr^{(2)}(\vartheta_j) = \frac{p^2-1}{8} + \frac{2r h_j^2}{p} - \frac{d(p,q,j)-d(p,q,j_0)}{2} \bmod 2.
\end{equation}
\end{lemma}

\begin{proof}
	Let $(X, \mathfrak{t}_{\mathrm{spin}})$ be the spin filling of $Y$ as in~\S\ref{ssec:classical-arf-K},
	and $\mathfrak{s}_{j_0}$ be the spin structure on $Y$.
	Since $p$ is odd, there is a unique $u \in H^2(Y;\mathbb{Z})$ for which $\mathfrak{s}_j = \mathfrak{s}_{j_0} + 2u$.
	Extend $u$ to a class $x \in H^2(X; \mathbb{Z})$, representing a line bundle $L$,
	and let $\mathfrak{t} = \mathfrak{t}_{\mathrm{spin}} \otimes L^2$.
	Then $c_1(\mathfrak{t}) = 4x$ and thus $Q_{\mathfrak{t}} = Q_{\mathfrak{t}_{\mathrm{spin}}} = Q$, and
	\begin{equation*}
		\Arf(Q_{\mathfrak{t}}) = \Arf(Q)= \frac{p^2-1}{8} \bmod 2.
	\end{equation*}

	Equip $Y$ with a PSC metric, so that the boundary Dirac operator is invertible. By the work of Ue \cite{Ue}, for spherical manifolds the Heegaard Floer invariant $d(Y,\mathfrak{s}_j)$ agrees with the Fukumoto-Furuta invariant, which implies that
	\begin{equation*}
		\ind^{\mathrm{APS}}_{\C} D^+_{\mathfrak{t}} = \frac{c_1(\mathfrak{t})^2 - \sigma(X)}{8} - \frac{d(Y, \mathfrak{s}_j)}{2}.
	\end{equation*}
	Hence
	\begin{equation*}
		\ind^{\mathrm{APS}}_{\C} D^+_{\mathfrak{t}}-\ind^{\mathrm{APS}}_{\C} D^+_{\mathfrak{t}_{\mathrm{spin}}}=
		2x^2 - \frac{d(Y,\mathfrak{s}_j) - d(Y,\mathfrak{s}_{j_0})}{2}.
	\end{equation*}
	Our choice of $r$ and $h_j$ implies $x^2 \equiv rh_j^2/p \pmod{\mathbb{Z}}$.
	To see this, express $2u$ as $(j-j_0)r\mu$, with $r,\mu$ as in \cite[Appendix~B]{CH15}.
With our orientation convention,
the boundary discriminant pairing $b_X$ satisfies $b_X(\mu,\mu)=q/p \pmod{\mathbb{Z}}$.
Since $p$ is odd and $2h_j\equiv j-j_0\pmod p$, we have $u=rh_j\mu$.
Hence
\[
x^2\equiv b_X(rh_j\mu,rh_j\mu)
=\frac{qr^2h_j^2}{p}
\equiv\frac{rh_j^2}{p}\pmod{\mathbb Z}.
\]
Finally, since the Dirac index of $\mathfrak{t}_{\mathrm{spin}}$ is even,
	\begin{equation*}
	\gr^{(2)}(\vartheta_j) = \ind^{\mathrm{APS}}_{\C} D^+_{\mathfrak{t}} + \Arf(Q_{\mathfrak{t}}) \equiv
	\frac{p^2-1}{8} + \frac{2rh_j^2}{p} - \frac{d(p,q,j) - d(p,q,j_0)}{2}. \qedhere
	\end{equation*}
\end{proof}

The following table gives some concrete examples.
\begin{center}
\small
\renewcommand{\arraystretch}{1.2}
\setlength{\tabcolsep}{4pt}
\begin{tabular}{lcrlllllr}
\hline
Knot
&$(p,q)$
&$j_0$
&$(\gr^{(2)}(\vartheta_j))_{0 \le j < p}$
&$\Delta_K(i)$
&$(-1)^{\sigma(K)/2}$
&$\mathrm{rk} \thmr_0$
&$\mathrm{rk} \thmr_1$ &$\hat{\chi}$\\
\hline
Trefoil
&$(3,1)$&0&$(\underline{1},0,0)$
&$-1$&$-1$&2&1 & $-1$\\
$T(2,5)$
&$(5,1)$&0&$(\underline{1},1,0,0,1)$
&$-1$&$+1$&2&3 & $-1$\\
Figure-eight
&$(5,2)$&3&$(0,0,0,\underline{1},0)$
&$3$&$+1$&4&1 & 3\\
$5_2$
&$(7,2)$&4&$(0,0,1,0,\underline{0},0,1)$
&$-3$&$-1$&5&2 &$-3$\\
$6_1$
&$(9,2)$&5&$(0,0,0,1,\underline{0},0,0,1,0)$
&$5$&$+1$&7&2 & $5$\\
\hline
\end{tabular}
\end{center}

\subsection{Odd torus knots}
Let $K = T(p,q)$, with $p,q > 1$ coprime.
Put one basepoint on $K$.
The double branched cover $Y$ is $\Sigma(2,p,q) = \{(z_1,z_2,z_3) \in \mathbb{S}(\mathbb{C}^3):z_1^2 + z_2^p + z_3^q = 0\}$, with deck transformation $\tau (z_1,z_2,z_3) = (-z_1,z_2,z_3)$.
Let $\vartheta_{\mathrm{spin}}$ be the unique spin reducible for the unique spin structure.
By Lemma~\ref{lem:spin-red-grading}:
\begin{equation}\label{eq:torus-spin}
\begin{array}{c|c|c|c}
\hline
\text{Parity of }(p,q)&\det T(p,q) & \sigma(K) \bmod 4 &\gr^{(2)}(\vartheta_{\mathrm{spin}})\\ \hline
p,q\text{ odd}&1&0 &0\\
p\text{ even},\ q\text{ odd}& q &-(q-1) &(q^2-1)/8\pmod2.\\
\hline
\end{array}
\end{equation}
If $p,q$ are both odd, then \cite[\S14.5]{Li22} shows that the real monopole chain complex contains no real irreducible generators. The unique real structure acts freely on irreducibles and fixes $\vartheta_{\mathrm{spin}}$.
In this case $\Sigma(2,p,q)$ is a $\mathbb{Z}HS^3$, so
\begin{equation*}
	\thmr(\Sigma_2(T(p,q))) = \mathbb{F}_{(0)}, \quad \text{$p,q$ odd}.
\end{equation*}
By comparison,
let $\Delta_K(t)$ be the Alexander polynomial with convention $\Delta_K(t)=\Delta_K(t^{-1})$ and $\Delta_K(1)=1$.
Then
\[\Delta_{T(p,q)}(t) = t^{-(p-1)(q-1)/2}\frac{(t^{pq}-1)(t-1)}{(t^p-1)(t^q-1)}, \qquad
\Delta_{T(p,q)}(i) = 1= \chi\bigl(\thmr(\dbcv_2(T(p,q)))\bigr).
\]

\subsection{$(2,p)$-torus knots}
Consider the 2-bridge torus knots $T(2,p)$, and use the convention $\dbcv_2(T(2,p))=L(p,1)$.
The real \text{spin\textsuperscript{c}} structures are $\{\mathfrak{s}_j: 0 \le j < p\}$,
with $j = 0$ corresponding to the spin structure.
Take $r=1$ and use any $h_j$ for which $2h_j \equiv j \bmod p$.
Formula \eqref{eq:Lpq-grading} becomes
\begin{equation*}
\gr^{(2)}(\vartheta_j)
=
\frac{p^2-1}{8}
+\frac{2h_j^2}{p}
-\frac{d(p,1,j)-d(p,1,0)}2
\bmod2,
\end{equation*}
Apply \eqref{eq:dpqj} to compute $d(p,1,j) = ((2j-p)^2-p)/4p$.
Plug in $h_j = j/2$ for $j$ even and $(j+p)/2$ for $j$ odd.
Then
\begin{equation*}
	d(p,1,j) - d(p,1,0) = -\frac{j(p-j)}{p}, \quad
	\frac{4h_j^2 + j(p-j)}{2p} = \begin{cases} j/2, & j \ \text{even}\\ (3j+p)/2 & j \ \text{odd}\end{cases}
\end{equation*}
Modulo two, we obtain
\begin{equation*}
\gr^{(2)}(\vartheta_j)
=
\frac{p^2-1}{8}
+\frac{j(p-j)}{2} \bmod 2.
\end{equation*}
In particular,
\begin{equation*}
	\chi (\thmr(\dbcv_2(T(2,p))))
=\sum_{j=0}^{p-1}
(-1)^{(p^2-1)/8+j(p-j)/2}
=(-1)^{(p^2-1)/8+(p-1)/2}.
\end{equation*}
On the other hand,
\[\Delta_{T(2,p)}(t) = t^{-(p-1)/2} \frac{t^p+1}{t+1}, \qquad
\Delta_{T(2,p)}(i) = (-1)^{(p^2-1)/8}.
\]
In other words,
\begin{equation*}
	(-1)^{(p-1)/2}\chi( \thmr(\dbcv_2(T(2,p))))
	 = \Delta_{T(2,p)}(i) = \widehat{\chi}\,(\hHFR(\dbcv_2(T(2,p)))).
\end{equation*}
Indeed, the normalization~\eqref{eq:seifert-normalization-gr2} for $T(2,p)$ is $-(p-1)/2$.
Hence the normalized Euler characteristics of the real monopole Floer homology agree with the Euler characteristics for real Heegaard Floer homology as in~\cite{Srivastava}.

\subsection{Montesinos knots}
Let $Y = \Sigma(p,q,r)$ be the Brieskorn sphere with $p,q,r \ge 2$ coprime,
and $\tau(z_1,z_2,z_3) =(\bar z_1, \bar z_2, \bar z_3)$ be the conjugation involution.
Then $Y/\tau \cong S^3$ and the branch set is a Montesinos knot $K_{p,q,r}$.
Since $Y$ is an integral homology sphere, there is a unique real \text{spin\textsuperscript{c}} structure arising from the spin structure.
We use the MOY setup~\cite{MOY} again, as adapted to the real case in \cite[\S14.6]{Li22}.
In particular, we use the Seifert metric connection rather than the typical Levi--Civita connection.
In \cite[\S14.6]{Li22} it is shown that the real structure preserves all the irreducibles and the unique spin reducible $\vartheta_{\mathrm{spin}}$ in the MOY description.
The same argument as in Lemma~\ref{lem:spin-red-grading} applies to this setting;
since $\det(K)=1$,
\begin{equation*}
\gr^{(2)}(\vartheta_{\mathrm{spin}}) = \Arf(K_{p,q,r}) = \frac{1-1}{8} = 0.
\end{equation*}
The gradings of all the critical points in the MOY description can be computed by combining the index calculation in~\cite[\S9]{MOY}, with the fact that $\sr$ acts conjugate-linearly on the moduli spaces of holomorphic divisors.
Normalizing $\vartheta_{\mathrm{spin}}$ to have grading zero mod two, this yields a complete calculation of the Euler characteristics.
Alternatively, the framed real monopole homology of $\Sigma(p,q,r)$ can be computed using the lattice techniques from~\cite{KPT:cables}.

\appendix
\section{Identification of the spin lifts}
\label{app:cover-indentification}
This appendix proves Lemma~\ref{lem:Gamma}.
Use the frame $(e_1,\nu_F,e_3)$ along $C$ on $Y$,
and the frame $(\partial_t,e_1,\nu_F,e_3)$ on the collar $Z = I \times C \subset I_t \times Y$ and bundle identification $S \cong S^+|_Y$ via $\rho_Y(v) = \rho(\partial_t)^{-1}\rho(v)$.
Denote
\[
D=\rho(\partial_t),\qquad T=\rho(e_1),\qquad
N=\rho(\nu_F),\qquad B=\rho(e_3).
\]
In particular,
\[
	1=\rho_Y(e_1)\rho_Y(\nu_F)\rho_Y(e_3)
	=(D^{-1}T)(D^{-1}N)(D^{-1}B)
	=-DTNB.
\]
For a unit normal vector $w\in\operatorname{span}(\nu_F,e_3)$
and a unit tangent vector $x\in\operatorname{span}(\partial_t,e_1)$,
consider the operators
\[
K_w=i\rho(\partial_t)^{-1}\rho(w),\qquad
H_{x^{-1}}=i\rho(e_3)^{-1}\rho(x).
\]
Here, $x^{-1}$ is the unit $\mathbb{C}$-covector for which $x^{-1}(x) = 1$.
Both $K$ and $H$ preserve $(S^+)^{\sr}$ and square to $1$.
Furthermore,
\begin{equation}\label{eq:collar-HK}
	H_{e_1^{-1}}=-iBT=iD^{-1}N=K_{\nu_F}.
\end{equation}
In fact, if $a^2+b^2=1$:
\begin{equation}\label{eq:collar-HK-ab}
H_{(ae_1-b\partial_t)^{-1}}=K_{a\nu_F+be_3}.
\end{equation}
The square map $\mu \colon E_-^{\otimes2}\to (TZ)_{\C}^{-1}$ from Lemma~\ref{lem:half4-new}
can be normalized by requiring,
for every $\psi \in U(E_-)$ and $x\in TZ$, that
\begin{equation}
	\label{eq:mu-normalization}
	-\rho(x)\rho(\nu_F)\psi
	=\mu(\psi^{\otimes2})(x)\,\sr\psi,
\end{equation}

\begin{proof}[Proof of Lemma~\ref{lem:Gamma}]
For $\psi\in E_-$, the relation $NB\psi=-i\psi$ gives $B\psi=iN\psi$.
Therefore, for unit $\psi$ and unit tangent $x$,
\[
H_{x^{-1}}\psi
=i\rho(x)B\psi
=-\rho(x)N\psi
=\mu(\psi^{\otimes2})(x)\,\sr\psi.
\]
The operator $H_{x^{-1}}$ commutes with $\sr$,
because $\tau_*$ fixes $x$ and negates $e_3$.
Writing $c=\mu(\psi^{\otimes2})(x)$, we obtain
\[
H_{x^{-1}}\Phi(\psi)
=\frac{\overline c\,\psi+c\,\sr\psi}{\sqrt2}.
\]
Since $\psi$ and $\sr\psi$ lie in distinct eigenlines,
the above equals $\Phi(\psi)$ if and only if $c=1$, or equivalently if $\mu(\psi^{\otimes2})=x^{-1}$.

Thus $\Phi$ identifies the spin lifts of $x^{-1}$ with the unit real
$+1$-eigenvectors of $H_{x^{-1}}$.
For $x=e_1$, Equation~\eqref{eq:collar-HK} says $H_{x^{-1}} = K_{\nu_F}$.
Hence $\Phi$ identifies the double cover $\nu_F^*U(S^{\sr})$ of $C$,
which consists of unit real $+1$-eigenvectors of $K_{\nu_F}$,
with the double cover consisting of square roots of $e_1(p)^{-1}$:
\[
\left\{\psi\in U(E_{-,p}): p \in C, \ \mu(\psi^{\otimes2})=e_1(p)^{-1}\right\}.
\]
It remains to show that the spin lift $T_k^{-1}(1,0)$ corresponds to a unit vector $\psi_k\in E_{-,p_k}$ satisfying
$\psi_k^{\otimes2}=e_1(p_k)^{-1}$.

Choose a local chart around $p_k$ for which $\sr$ is complex conjugation on $\mathbb{C}^2$,
and $\rho$ satisfies~\eqref{eq:clifford-in-adapted-frame}.
We have
\[
J_N=NB=-\sigma_2,\qquad
E_-=\C\psi_0,\quad
\psi_0=\frac1{\sqrt2}\binom{1}{-i};
\qquad
K_{\nu_F} = \begin{pmatrix}1 & 0 \\ 0 & -1\end{pmatrix}.
\]
Moreover,
\[
-TN\psi_0=\overline{\psi_0},\qquad
\mu(\psi_0^{\otimes2})=e_1^{-1},\qquad
\Phi(\psi_0)=\binom10.
\]
The two lifts of $e_1^{-1}$ are thus $\pm\psi_0$.
Under $\Phi$, they map respectively to $\pm(1,0)$,
the two unit real $+1$-eigenvectors of $K_{\nu_F}$.
Any transition map between two such charts commutes with complex conjugation and Clifford multiplication,
forcing it to lie in $\{\pm1\}$.
Since $\{\pm 1\}$ acts the same way on $\psi_0$ and on $(1,0)$,
the two double covers above have the same transition functions along $C$.
\end{proof}

\bibliography{arf}
\bibliographystyle{custom}
\end{document}